\documentclass[11pt]{article}
\usepackage{amsmath,amssymb,amsthm,bm}
\usepackage{graphicx,color}
\usepackage{subfigure}
\usepackage{epsfig}

\newcommand{\beq}{\begin{equation}}
\newcommand{\eeq}{\end{equation}}
\newtheorem{theorem}{Theorem}[section]
\newtheorem{lemma}[theorem]{Lemma}
\newtheorem{corollary}[theorem]{Corollary}
\newtheorem{proposition}[theorem]{Proposition}
\newtheorem{definition}[theorem]{Definition}

\newcommand{\goto}{\rightarrow}

\newcommand{\<}{\langle}
\renewcommand{\>}{\rangle}

\newcommand{\vf}{\varphi}

\newcommand{\R}{{\mathbb R}}
\newcommand{\C}{{\mathbb C}}
\newcommand{\bZ}{{\mathbb Z}}
\newcommand{\Z}{{\mathbb Z}}

\newcommand{\cR}{{\cal R}}

\newcommand{\jlk}{{j,\ell,k}}

\newcommand{\bitem}{\begin{itemize}}
\newcommand{\eitem}{\end{itemize}}

\title{The Curvelet Representation of Wave Propagators is Optimally Sparse}

\author{Emmanuel J.~Cand\`es and Laurent Demanet\\
  Applied and Computational Mathematics\\
  California Institute of Technology\\
  Pasadena, California 91125 }
\date{December 2003} 

\begin{document}

\maketitle
\begin{abstract}
  This paper argues that {\em curvelets} provide a powerful tool for
  representing very general linear symmetric systems of hyperbolic
  differential equations. Curvelets are a recently developed
  multiscale system \cite{TVSynthesis,CurveEdge} in which the elements
  are highly anisotropic at fine scales, with effective support shaped
  according to the parabolic scaling principle $width \approx
  length^2$ at fine scales. We prove that for a wide class of linear
  hyperbolic differential equations, the curvelet representation of
  the solution operator is both optimally sparse and well organized.
  \begin{itemize}
  \item It is sparse in the sense that the matrix entries decay nearly
    exponentially fast (i.e. faster than any negative polynomial),
  \item and well-organized in the sense that the very few
    nonnegligible entries occur near a shifted diagonal. 
  \end{itemize}
  Indeed, we actually show that the action of the wave-group on a
  curvelet is well-approximated by simply translating the center of
  the curvelet along the Hamiltonian flow---hence the diagonal shift
  in the curvelet representation.  A physical interpretation of this
  result is that curvelets may be viewed as coherent waveforms with
  enough frequency localization so that they behave like waves but at
  the same time, with enough spatial localization so that they 
  simultaneously behave like particles.
\end{abstract}


{\bf Acknowledgments.}  This research was supported by a National
Science Foundation grant DMS 01-40698 (FRG), by a DOE grant
DE-FG03-02ER25529, and by an Alfred P.~Sloan Fellowship.

\pagebreak

\section{Introduction}
\label{sec:introduction}
\setcounter{theorem}{0}
\setcounter{equation}{0}

\renewcommand{\d}{{\partial}}
\newcommand{\bu}{{\bf u}}
\newcommand{\bx}{{\bf x}}

This paper is concerned with the representation of symmetric systems
of linear hyperbolic differential equations of the form
\begin{equation}
  \label{eq:hyperbolic-system}
 \frac{\d u}{\d t} + \sum_k A_k(x) \frac{\d u}{\d x_k} + B(x) u = 0, 
\qquad u(0,x) = u_0(x),    
\end{equation}
where $u$ is an $m$-dimensional vector and $x \in \R^n$.  The matrices
$A_k$ and $B$ may depend on the spatial variable $x$, and the $A_k$'s
are symmetric. Linear hyperbolic systems are ubiquitous in the
sciences and a classical example are the equations for acoustic waves
which read
\begin{equation}
 \label{eq:system-wave}
  \begin{array}{ll}
      \frac{\d \rho}{\d t} + \nabla\cdot (\rho_0 u) & = 0 \\ 
      \rho_0 \frac{\d u}{\d t} + \nabla (c^2_0 \rho) & = 0,
  \end{array}
\end{equation}
where $u$ and $\rho$ are the velocity and density disturbances
respectively. (Here, $\rho_0 = \rho_0(x)$ is the density and $c_0 =
c_0(x)$ is the speed of sound.) Other well-known examples include
Maxwell's equations of electrodynamics and the equations of linear
elasticity. All the results presented in this paper equally apply to
higher-order scalar wave equations, e.g., of the form
\[
\frac{\d^2 u}{\d t^2} - \sum_{ij}a_{ij}(x) \frac{\d^2 u}{\d x_j \d x_k}
  = 0, \qquad u(0,x) = u_0(x), \quad \frac{\d u}{\d t}(0,x) = u_1(x),
\]
(where $u$ is now a scalar and $a_{ij}(x)$ is taken to be symmetric
and positive definite) as it is well-known that such single
second-order equations can be reduced to a symmetric system of
first-order equations (\ref{eq:hyperbolic-system}) by appropriate
changes of variables.


\subsection{About representations}

We are interested in representations of the solution operator $E(t)$
to the system (\ref{eq:hyperbolic-system}) 
\[
u(t, \cdot) = E(t) u_0, 
\]
which may be expressed as an integral involving the so-called
Green's function
\[
u(t,x) = \int E(t; x,y) u_0(y) \, dy. 
\]
To introduce the concept of representation, suppose that the
coefficient matrices do not depend on $x$. As is well-known, the
Fourier transform is a powerful tool for studying differential
equations in this setting.  Indeed, in the Fourier domain,
(\ref{eq:hyperbolic-system}) takes the form
\[
(\d_t + i \sum_k A_k \xi_k + B)\hat{u}(t,\xi) = 0; 
\]
in short, (\ref{eq:hyperbolic-system}) reduces to a system of ordinary
differential equations which can be solved analytically. This shows
the power of the representation; in the frequency domain, the solution
operator is diagonal and the study becomes ridiculously simple.

These desirable properties are very fragile, however. Both
mathematicians and computational scientists know that Fourier methods
are not really amenable to differential equations with variable
coefficients, and that we need to find alternatives.  Instead of
considering the evolution of Fourier coefficients, we may want to
think, instead, of the action of the propagation operator $E(t)$ on
other types of basis elements. This connects with the viewpoint of
modern harmonic analysis whose goal is to develop representations,
e.g. an orthonormal basis $(f_n)$ of $L^2(\R^n)$, say, in which the
solution operator
\begin{equation}
\label{eq:sol-operator-basis}
E(t; n, n') = \<f_n, E(t) f_{n'} \> 
\end{equation}
is as simple as possible; that is, such that $E(t) f_{n}$ is a {\em
  sparse} superposition of those elements $f_{n'}$. Such sparse
representations are extremely significant both in mathematical
analysis, where sparsity allows for sharper inequalities and in
numerical applications where sparsity allows for faster algorithms. 
\begin{itemize}
\item In the field of mathematical analysis, for example, Calder\'on
  introduced what one would nowadays call the Continuous Wavelet
  Transform (CWT) in which objects are represented as a superpositions
  of simple elements of the form $\psi((x-b)/a)$, with $a > 0$ and $b
  \in \R^n$; i.e, objects are represented as a superposition of
  dilates and translates of a single function $\psi$. These elements
  proved to be {\em almost} eigenfunctions of large classes of
  operators, the Calder\'on-Zygmund operators which are special types
  of singular integrals, some of which arising in connection with
  elliptic problems. It was later gradually realized that tools like
  atomic decompositions of Hardy spaces \cite{FJW,Ste} and orthonormal
  bases of Wavelets \cite{Mey1,Mey2} provide a setting in which some
  aspects of the mathematical analysis of these operators is
  dramatically eased.
  
\item Clever representation of scientific and engineering computations
  can make previously intractable computations tractable.  Here,
  sparsity may allow the design of fast matrix multiplication and/or
  fast matrix inversion algorithms. For example, Beylkin, Coifman and
  Rokhlin \cite{BCR} exploited the sparsity of those singular
  integrals mentioned above, and showed how to use wavelet bases to
  compute such integrals with very low complexity algorithms.
\end{itemize}
In short, a single representation, namely, the wavelet transform
provides sparse decompositions of large classes of operators
simultaneously. 

\subsection{Limitations of Classical Multiscale Ideas}

Our goal in this paper is to find a representation which provides
sparse representations of the solution operators to fairly general
classes of systems of hyperbolic differential equations.  Now the last
two decades have seen the widespread development of multiscale ideas
such as Multigrid, Fast Multipole Methods, Wavelets, Finite Elements with or
without adaptive refinement, etc. All these representations propose
dictionaries of roughly isotropic elements occurring at all scales and
locations; the templates are rescaled treating all directions in
essentially the same way.  Isotropic scaling may be successful when
the object under study does not exhibit any special features along
selected orientations. This is the exception rather than the rule. 

Tools from traditional multiscale analysis are very powerful for
representing certain elliptic problems but unfortunately, they are
definitely ill-adapted to hyperbolic problems such as
\eqref{eq:hyperbolic-system}.  Indeed,
\begin{enumerate}
\item they fail to sparsify the wave propagation, i.e. the solution
  operator $E(t)$, 
\item and they fail to provide a sparse representation of oscillatory
  signals which are the solutions of those equations.
\end{enumerate}
To make things concrete, consider the problem of propagating elastic
waves as in geophysics. Consider the scalar wave equation in two dimensions 
\begin{equation}
  \label{eq:wave}
  \d_{tt} u = c^2(x) \Delta u,   
\end{equation}
where $\Delta$ is the Laplacian defined by $\Delta =
\partial^2/\partial x_1^2 + \partial^2/\partial x_2^2$ (we may take
the velocity field to be constant). To describe the action of the
wave group, we assume that the initial condition takes the form of a
wavelet with vanishing initial velocity, say.  Then it is clear that
at a later time, the wavefield is composed of large concentric rings
(imagine throwing a stone in a lake). Now, it is also clear that many
wavelets are needed to represent the wavefield.  In other words, the
wavefield is a rather dense superposition of wavelets. Note that this
may be quantified.  Suppose that the velocity field is identically
equal to one, say, and that the initial condition is a wavelet at
scale $2^{-j}$; that is, of the form $2^{j/2} \psi(2^j x)$ so that in
frequency, the energy is concentrated near the dyadic annulus $|\xi|
\sim 2^j$. Then one would need at least $O(2^j)$ wavelets to
reconstruct the wavefield at time $t = 1$ to within reasonable
accuracy.

Our simple example above shows wave-like flows do not preserve the
geometry and characteristics of classical multiscale systems.  To
achieve sparsity, we need to rethink the geometry of multiscale
representations.

\subsection{A New Form of Multiscale Analysis}

As we will see in section \ref{sec:curvelets}, curvelets are waveforms
which are highly anisotropic at fine scales, with effective support
obeying the parabolic principle $length \approx width^2$. Just as for
wavelets, there is both a continuous and a discrete curvelet
transform.  A curvelet is indexed by three parameters which---adopting
a continuous description of the parameter space---are: a scale $a$, $0
< a < 1$; an orientation $\theta$, $\theta \in [-\pi/2, \pi/2)$ and a
location $b$, $b \in \R^2$. At scale $a$, the family of curvelets is
generated by translation and rotation of a basic element $\vf_a$
\[
\vf_{a,b,\theta}(x) = \vf_a(R_\theta(x - b)).
\]
Here, $\vf_{a}(x)$ is some kind of directional wavelet with spatial
width $\sim a$ and spatial length $\sim \sqrt{a}$, and with minor axis
pointing in the horizontal direction
\[
\vf_a(x) \approx \vf(D_a x), \quad D_a = \begin{pmatrix} 1/a & 0 \\
  0 & 1/\sqrt{a}
\end{pmatrix}; 
\]
$D_a$ is a parabolic scaling matrix, $R_\theta$ is a rotation by
$\theta$ radians. 

An important property is that curvelets obey the principle of harmonic
of analysis which says that it is possible to analyze and reconstruct
an arbitrary function $f(x_1,x_2)$ as a superposition of such
templates.  It is possible to construct tight-frames of curvelets and
one can, indeed, easily expand an arbitrary function $f(x_1, x_2)$ as
a series of curvelets, much like in an orthonormal basis.  Continuing
at an informal level of exposition, there is a sampling of the plane
$(a,b,\theta)$
\[
a_j = 2^{-j}, \quad 
\theta_{j,\ell} = 2\pi \ell \cdot 2^{-\lfloor j/2 \rfloor}, \quad
 R_{\theta_{j,\ell}}  b_k^{(j,\ell)} = (k_1 2^{-j}, k_2 2^{-j/2}),  
\]
such that with $\mu$ indexing the triples $(a_j, \theta_{j,\ell},
b_k^{(j,\ell)})$ the collection $\vf_\mu$ is a tight-frame:
\beq\label{eq:tight-frame}
f = \sum_{\mu} \langle f, \vf_\mu  \rangle \vf_\mu, \quad 
\| f \|_2^2 = \sum_\mu |\langle f, \vf_\mu \rangle|^2. 
\eeq
(Note that these formulae allow us to analyze and synthesize arbitrary
functions in $L^2(\R^2)$ as a superposition of curvelets in a stable
and concrete way.) 

\begin{figure}
  \centering 
  \includegraphics[scale=.55]{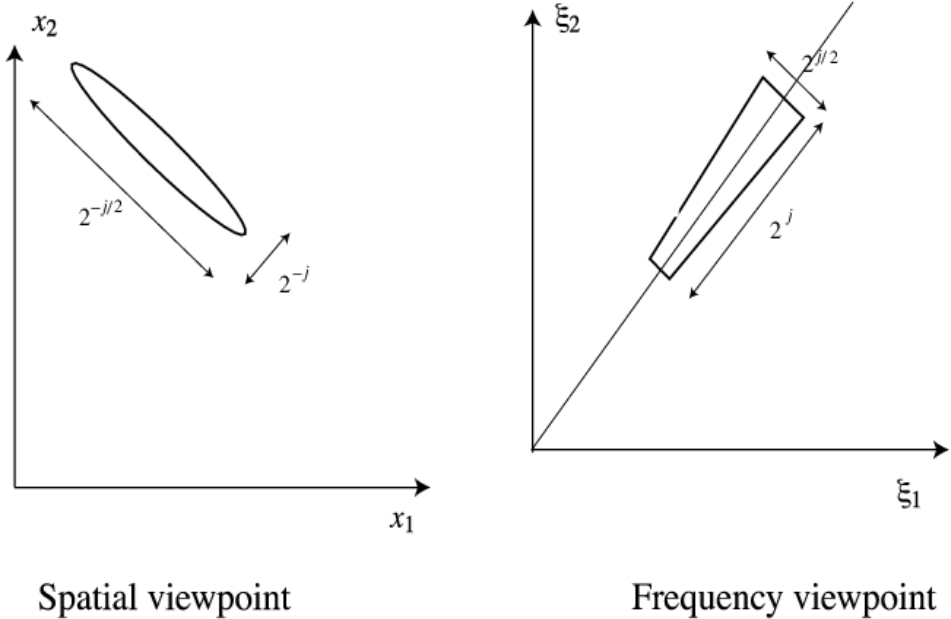}
  \caption{Schematic representation of the support of a 
    curvelet in both space and frequency.  In the spatial domain, a
    curvelet has an envelope strongly aligned along a specified
    `ridge' while in the frequency domain, it is supported near a box
    whose orientation is aligned with the co-direction of the ridge.}
  \label{fig:SpaceFreqView}
\end{figure}

As we have seen, a curvelet is well-localized in space but it is also
well-localized in frequency. Recall that a given scale, curvelets
$\vf_\mu$ are obtained by applying shifts and rotations to a `mother'
curvelet $\vf_{j,0,0}$. In the frequency domain then
\begin{equation*}
  \label{eq:vfj1}
  \hat{\vf}_{j,0,0}(\xi) = 
2^{-3j/4} W(2^{-j}|\xi|) \, V(2^{\lfloor j/2 \rfloor} \theta).  
\end{equation*}
Here, $W, V$ are smooth windows compactly supported near the intervals
$[1, 2]$ and $[-1/2, 1/2]$ respectively. Whereas in the spatial domain
curvelets live near an oriented rectangle $R$ of length $2^{-j/2}$ and
width $2^{-j}$, in the frequency domain, they are located in a parabolic
wedge of length $2^j$ and width $2^{j/2}$ and whose orientation is
orthogonal to that of $R$. The joint localization in both space allows
us to think about curvelets as occupying a `Heisenberg cell' in
phase-space with parabolic scaling in both domains. Figure
\ref{fig:SpaceFreqView} offers a schematic representation of this
joint localization. As we shall see, this microlocal behavior is key
to understanding the properties of curvelet-propagation. Additional
details are given in section \ref{sec:curvelets}.

\subsection{Curvelets and Geometrical Optics}\label{sec:curvgeomoptics}

A hyperbolic system can typically be considered in the approximation
of high-frequency waves, also known as {\em geometrical optics}. In
order to best describe our main result, it is perhaps suitable first
to exhibit the connections between curvelets and geometrical optics. In
that setting it is not necessary to describe the dynamics in terms of
the wavefield $u(t,x)$. Only its prominent features are studied: wave
fronts, or equivalently rays. The latter are trajectories $(x(t),
\xi(t))$ in phase-space $\R^2 \times \R^2$, and are the solutions to
the $m$ Hamiltonian flows (indexed by $\nu$)
\begin{equation}\label{eq:HS}
\left\{ \begin{array}{llll}
\dot{x}(t) &= \nabla_\xi \lambda^0_\nu(x,\xi), 
\qquad &x(0) &= x_0,\\
\dot{\xi}(t) &= -\nabla_x \lambda^0_\nu(x,\xi), 
\qquad & \xi(0) &= \xi_0.
           \end{array}
      \right.
\end{equation}
The system (\ref{eq:HS}) is also called the {\em bicharacteristic
  flow} and the rays $(x(t), \xi(t))$ the {\em bicharacteristics.}  To
see how this system arises, consider the classical high-frequency
wave-propagation approximation to the wavefield $u(x,t)$ of the form
\[
u(x,t) = \sum_\nu e^{i \omega \Phi_\nu(x,t)} \left( \sigma_{\nu,0}(x,t) +
  \frac{\sigma_{\nu,1}(x,t)}{\omega} + \frac{\sigma_{\nu,2}(x,t)}{\omega^2} +
  \ldots\right),
\]
where $\omega$ is a large parameter; it then follows (after
substituting the approximation in the wave equation
\eqref{eq:hyperbolic-system}) that the phase functions $\Phi_\nu$ must
obey the eikonal equations
\begin{equation}
  \label{eq:eikonal}
  \partial_t \Phi_\nu + \lambda^0_\nu(x,\nabla_x \Phi) = 0, 
\end{equation}
and that the amplitudes must obey `transport' equations we shall not
detail here (see section \ref{sec:lax}). In the above expression (and
in the Hamiltonian equations \eqref{eq:HS}), the
$\lambda^0_\nu(x,\xi)$ are the eigenvalues of the dispersion matrix
\begin{equation}
  \label{eq:dispersion}
  a^0(x,\xi) = \sum_k A_k(x) \xi_j. 
\end{equation}
It is well-known that the Hamiltonian equations describe the evolution
of the wavefront set of the solution as $\Phi_\nu(t,x(t))$ is actually
constant along the $\nu$th Hamiltonian flow \eqref{eq:HS}. 

With this background, we are now in a position to qualitatively
describe the behavior of the wave-propagation operator $E(t)$ acting
on a curvelet $\vf_\mu$. However, we first need to introduce a notion
of vector-valued curvelet since $E(t)$ is acting on vector fields. 
Let $r^0_\nu(x,\xi)$ be the eigenvector of the dispersion matrix
associated with the eigenvalue $\lambda^0_\nu(x,\xi)$. We then define
hyper-curvelets by
\begin{equation}
  \label{eq:hyper1}
\bm{\vf}^{(0)}_{\mu \nu}(x) = \frac{1}{(2 \pi)^2} 
\int e^{i x \cdot \xi} r^0_\nu(x, \xi) \hat{\vf}_{\mu}(\xi) \, d\xi. 
\end{equation}
Later in this section, we will motivate this special choice but for
now simply observe that $\bm{\vf}^{(0)}_{\mu \nu}$ is a vector-valued waveform.

\begin{figure}
  \centering \includegraphics[scale=.45]{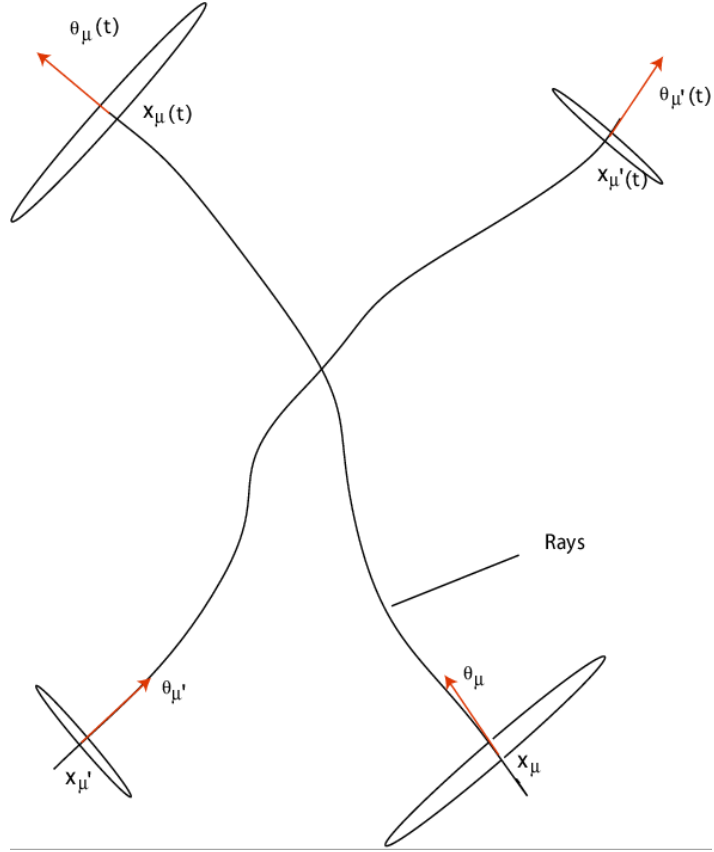}
  \caption{Schematic representation of the action of the wave group 
    on a hyper-curvelet. The solution is well-approximated by rigid
    motion along the Hamiltonian flow.}
  \label{fig:Hflow} 
\end{figure}

Consider then the solution to the wave equation $\bm{\vf}^{(0)}_{\mu \nu}(t,x)$
with initial value $\bm{\vf}^{(0)}_{\mu \nu}(x)$. Then we claim that  
\begin{center}
  {\em the action of the wave group on a hyper-curvelet is
    well-approximated by simply translating the center of the curvelet
    along the corresponding Hamiltonian flow.}
\end{center}
To examine this claim, set $n(t) = \xi(t)/|\xi(t)|$ and consider the
reduced Hamiltonian flow obtained by adding the equation 
\begin{equation}
  \label{eq:Hamiltonian-flow}
  \dot{U}(t)  = - U(n \otimes \nabla_x \lambda^0_\nu(x,n) - 
                    \lambda^0_\nu(x,n) \otimes n)
\end{equation}
to the system \eqref{eq:HS}. Here, $U(t)$ is the rotation matrix
tracking the orientation rotation of $\xi(t)$ in the sense that at all
times $U(t) n(t) = n(0) = \xi(0)/|\xi(0)|$, or $n(t) = U(t)^* n(0)$.
Our claim says that the solution to the wave equation nearly follows
the dynamics of the reduced Hamiltonian flow, i.e.
\begin{equation}
  \label{eq:approx-flow}
    \bm{\vf}^{(0)}_{\mu \nu}(t,x) \approx \bm{\vf}^{(0)}_{\mu \nu}(U_\mu(t)(x - x_\mu(t)) + x_\mu) 
\end{equation}
where $x_\mu(t)$ and $U_\mu(t)$ are the rays of
\eqref{eq:Hamiltonian-flow} with initial starting points $x_\mu(0) =
x_\mu$ and $U_\mu(0) = Id$. Figure \ref{fig:Hflow} illustrates this
fact.

We now return to the interpretation of a hyper-curvelet. Suppose that
$r^0_\nu$ only depends on $\xi$ as in the case of the acoustic
system \eqref{eq:system-wave} 
\begin{displaymath}
  r^0_0(\xi)= 
  \begin{pmatrix}
    \xi^\perp/|\xi| \\ 0
  \end{pmatrix},\qquad 
  r^0_\pm(\xi) = \frac1{\sqrt{2}} 
  \begin{pmatrix}
   \pm \xi /|\xi| \\ 1
  \end{pmatrix}.
\end{displaymath}
(Here and below, $\xi^\perp$ denotes the vector obtained from $\xi$
after applying a rotation by 90 degrees).  In this special case, we
see that the hyper-curvelet is obtained by multiplying---in the frequency
domain---a scalar-valued curvelet with the eigenvectors of the
dispersion matrix
\[
\hat{\bm{\vf}}^{(0)}_{\mu \nu}(\xi) = r^0_\nu(\xi) \hat{\varphi}_\mu(\xi),
\quad \nu \in \{+, -, 0\}. 
\]
This is useful for the curvelet $\hat{\bm{\vf}}^{(0)}_{\mu \nu}$ will
essentially follow only one flow, namely, the $\nu$th flow. Suppose we
had started, instead, with an initial value of the form
$\bm{\vf}_{\mu\nu} = \vf_\mu \mathbf{e}_\nu$, where $\mathbf{e}_\nu$
is the canonical basis of $\R^3$, say. Then our curvelet would have
interacted with the three eigenvectors of the dispersion matrix, and
would have `split' and followed the three distinct flows.  By forcing
$\hat{\bm{\vf}}^{(0)}_{\mu\nu}(\xi)$ to be aligned with $r^0_\nu(\xi)$, we
essentially removed the components associated with the other flows.
In the general case \eqref{eq:hyper1}, we build hyper-curvelets by
applying $R^0_\nu$, which is now a pseudo-differential operator with
symbol $r^0_\nu(x,\xi)$, mapping scalars to $m$-dimensional vectors,
and independent of time.  The effect is of course the same.

Note that when $r^0_\nu$ is independent of $x$, hyper-curvelets build-up
a (vector-valued) tight-frame; letting $[F,G]$ be the usual inner
product over three-dimensional vector fields in $L^2(\R^2)$, the
family $(\bm{\vf}^{(0)}_{\mu \nu})_{\mu \nu}$ obeys the reconstruction formula
\begin{equation}
  \label{eq:reproducing-formula}
u = \sum_{\mu, \nu} [u, \bm{\vf}^{(0)}_{\mu \nu}] \bm{\vf}^{(0)}_{\mu \nu}
\end{equation}
and the Parseval relation 
\begin{equation}
  \label{eq:Parseval}
 \|u\|_{L^2}^2 = \sum_{\mu, \nu} |[u, \bm{\vf}^{(0)}_{\mu \nu}]|^2.  
\end{equation}
Just as one can decompose a scalar field as a superposition of scalar
curvelets, one can analyze and synthesize any wavefield as a
superposition of hyper-curvelets in a stable and concrete way. For
arbitrary $r^0_\nu(x,\xi)$, this is, however, in general not true.

We would like to emphasize that although the Eikonal equations only
have solutions for small times, the approximation
\eqref{eq:approx-flow} and, more generally, all of our results are
valid for all times since the rays \eqref{eq:HS} are always
well-defined, see section 1.6 below for a more detailed discussion.

\subsection{Curvelets and Hyperbolic Systems}

The previous section gave a qualitative description of the action of
the wave group on a curvelet and we we shall now quantify this fact.
The evolution operator $E(t)$ acting on a curvelet
$\bm{\vf}^{(0)}_{\mu_0 \nu_0}$ is of course not exactly another
curvelet $\bm{\vf}^{(0)}_{\mu_0(t) \nu_0}$ which occurs at a displaced
location and orientation. Instead, it is a superposition of curvelets
$\sum_{\mu, \mu} \alpha_{\mu \nu} \bm{\vf}^{(0)}_{\mu \nu}$ such that
\begin{enumerate}
\item the coefficients $(\alpha_{\mu\nu})$ decay nearly exponentially,
  \item and the significant coefficients of this expansion are all
    located at indices $(\mu, \nu)$ `near' $(\mu_0(t), \nu_0)$. By
    near, we mean nearby scales, orientations and locations.
\end{enumerate}
  
To state the key result of this paper, we need a notion of distance
$\omega$ between curvelet indices which will be formally introduced in
section \ref{sec:curvelets}. Crudely, $\omega(\mu,\mu')$ is small if and only if
both curvelets are at roughly the same scale, have similar orientation
and are at nearby spatial locations. In the same spirit, the distance
$\omega(\mu,\mu')$ increases as the distance between the scale, angular,
and location parameters increases.

For each $\mu = (j,k,\ell)$ and $\nu = 1, \ldots, m$, define the
vector-valued curvelets
\begin{equation}
  \label{eq:hyper2}
\bm{\vf}_{\mu\nu} = \mathbf{e}_\nu \vf_\mu,
\end{equation}
where $\mathbf{e}_\nu$ is the $\nu$th canonical basis vector in
$\R^m$. The $\bm{\vf}_{\mu\nu}$'s inherit the tight-frame property
\eqref{eq:reproducing-formula}--\eqref{eq:Parseval}. We would like to
again remind the reader that these vector-valued curvelets are simpler
and different from the hyper-curvelets $\bm{\vf}_{\mu\nu}^{(0)}$
defined in the previous section.  Consider now the representing the
operator $E(t)$ in a tight-frame of vector-valued curvelets, namely,
\begin{equation}
  \label{eq:rep-op}
  E(t;\mu ,\nu; \mu', \nu') = \<\bm{\vf}_{\mu\nu},  E(t) \bm{\vf}_{\mu' \nu'}\>. 
\end{equation}
We will refer to $E(t; \mu,\nu; \mu',\nu')$ or simply $E$ as the curvelet matrix of
$E(t)$, with row index $\mu,\nu$ and column index $\mu',\nu'$. Decompose 
the initial wavefield $u_0 = \sum_{\mu, \nu}
c_{\mu\nu} \bm{\vf}_{\mu\nu}$. Then one can express the action of
$E(t)$ on $u_0$ in the curvelet domain as
\[
E(t)u_0 = \sum_{\mu \nu} c_{\mu\nu}(t) \bm{\vf}_{\mu\nu}, \quad
c_{\mu\nu}(t) = \sum_{\mu',\nu'} E(t; \mu',\nu'; \mu,\nu)  
c_{\mu'\nu'}
\]
with convergence in $L^2(\R^2, \C^m)$.  In short, the curvelet matrix
maps the curvelet coefficients of the initial wavefield $u_0(\cdot)$
into those of the solution $u(t, \cdot)$ at time $t$.  

\begin{theorem}
\label{teo:main}
Suppose that the coefficients $A_k(x)$ and $B(x)$ of the hyperbolic
system are $C^\infty$ and that the multiplicity of the eigenvalues of
the dispersion matrix $\sum_k A_k(x) \xi_k$ is constant in $x$ and
$\xi$. Then
  \begin{itemize}
  \item The matrix $E$ is sparse. Suppose $a$ is either a row or a
    column of $E$, and let $|a|_{(n)}$ be the $n$-largest entry of the
    sequence $|a|$, then for each $M > 0$, $|a|_{(n)}$ obeys
\beq\label{eq:main1}
|a|_{(n)} \le C_{tM} \cdot n^{-M}.
\eeq

\item The matrix $E$ is well-organized. For each $N > 0$, the
  coefficients obey
\beq\label{eq:main2}
|E(t;\mu,\nu; \mu',\nu')| \le C_{tN} \cdot \sum_{\nu'' = 1}^m
\omega(\mu,\mu'_{\nu''}(t))^{-N}.
\eeq
Here $\mu_\nu(t)$ is the
curvelet index $\mu$ flown along the $\nu$th Hamiltonian system.
\end{itemize}
Both constants $C_{tM}$ and $C_{tN}$ grow in time at most like $C_1 e^{C_2t}$
for some $C_1, C_2 > 0$ depending on $M$, resp. $N$.

\end{theorem}
In effect, the curvelet matrix of the solution operator resembles a
sum of $m$ permutation matrices where $m$ is the order of the
hyperbolic system; first, there are significant coefficients along $m$
shifted diagonal and second, coefficients away from these diagonals
decay nearly exponentially; i.e. faster than any negative polynomial.
Now just as wavelets provide sparse representations to the solution
operators to certain elliptic differential equations, our theorem
shows that curvelets provide an optimally sparse representation of
solution operators to systems of symmetric hyperbolic equations.

We can also resort to hyper-curvelets as defined in the previous
section and formulate a related result where the curvelet matrix is
sparse around a \emph{single} shifted diagonal. This refinement
approximately decouples the evolution into polarized components and
will be made precise later.

To grasp the implications of Theorem \ref{teo:main}, consider the
following corollary:
\begin{corollary}
\label{teo:maincoro}
  Consider the truncated operator $A_B$ obtained by keeping $m \cdot B$
  elements per row---the $B$ closest to each shifted diagonal in the
  sense of the pseudo-distance $\omega$. Then the truncated matrix obeys
\begin{equation}
\label{eq:approx_T}
\|A - A_B\|_{L^2 \goto L^2} \le C_M \cdot B^{-M}, 
\end{equation}
for each $M > 0$. 
\end{corollary}
The proof follows from that of Theorem \ref{teo:main} by an
application of Schur's lemma and is omitted.  Hence, whereas the
Fourier or wavelet representations are dense, curvelets faithfully
model the geometry of wave propagation as only a few terms are needed
to represent the action of the wave group accurately.

\subsection{Strategy}\label{sec:strategy}

In his seminal paper \cite{Lax}, Lax constructed approximate
solution operators to linear and symmetric hyperbolic systems, also
known as {\em parametrices.} He showed that these parametrices are
oscillatory integrals in the frequency domain which are commonly
referred to as Fourier integral operators (FIO) (the development and
study of FIOs is motivated by the connection). An operator $T$ is said
to be an FIO if it is of the form
\begin{equation}\label{eq:FIO}
  Tf(x) = \int e^{i \Phi(x,\xi)} \sigma(x,\xi) \hat{f}(\xi) \, d\xi.
\end{equation}
We suppose the phase function $\Phi$ and the amplitude $\sigma$ obey
the following standard assumptions \cite{Ste}:
\begin{itemize}
\item the phase $\Phi(x,\xi)$ is $C^\infty$, homogeneous of degree 1
  in $\xi$, i.e. $\Phi(x,\lambda \xi) = \lambda \Phi(x,\xi)$ for
  $\lambda > 0$, and with $\Phi_{x \xi} = \nabla_x \nabla_{\xi} \Phi$,
  obeys the nondegeneracy condition
\begin{equation}
\label{eq:nondegenerate}
|\mbox{det}\, \Phi_{x \xi}(x, \xi)| > c > 0, 
\end{equation}
uniformly in $x$ and $\xi$;
\item the amplitude $\sigma$ is a symbol of order $m$, which means that $\sigma$
  is $C^\infty$, and obeys
\begin{equation}
\label{eq:symbol}
|\partial^\alpha_\xi \partial^\beta_x \sigma(x, \xi)| \le C_{\alpha\beta}(1
+ |\xi|)^{m - |\alpha|}.
\end{equation}
\end{itemize}

Lax's insight is that the solution of the initial value problem for a
variable coefficient hyperbolic system can be well approximated by a
superposition of integrals of the form \eqref{eq:FIO} with
matrix-valued amplitudes of order 0.  The phases of these FIO's are, of course,
those solving the Eikonal equations \eqref{eq:eikonal}. Hence, a
substantial part of our argument will be about proving that curvelets
sparsify FIO's.  Now an important aspect of this construction is that
this approximation is only valid for small times whereas our theorem
is valid for {\em all} times.  The reason is that the solutions to the
Eikonal equations \eqref{eq:eikonal} are not expected to be global in
time, because $\Phi_\nu$ would become multi-valued when rays
originating from the same point $x_0$ cross at a later time.  This
typically happens at cusp points, when caustics start developing.  We
refer the reader to \cite{Gar, Whi}. Because, we are interested in a
statement valid for all times, we need to bootstrap the construction
of the FIO parametrix by composing the small time FIO parametrix with
itself.  Now this creates an additional difficulty.  Each parametrix
convects a curvelet along $m$ flows, and we see that after each
composition, the number of curvelets would be multiplied by $m$, see
section 4.1 for a proper discussion.  This would lead to matrices with
poor concentration properties. Therefore, the other part of the
argument consists in decoupling the equations so that this phenomenon
does not occur. In summary, the general architecture of the proof of
Theorem \ref{teo:main} is as follows:
\begin{itemize}
\item We first decompose the wave-field into $m$-one way components,
  i.e. components which essentially travel along only one flow. We show
  that this decomposition is sparse in tight-frames of curvelets.
\item Second, we show that curvelet representations of FIO's are
  optimally sparse in tight-frame of curvelets, a result of
  independent interest. 
\end{itemize}

As a side remark, we would like to point out that the result about
optimally sparse representations of FIO's was announced without a
proof in the companion paper \cite{CurveFIO}. This paper, however,
gives the first proof of this optimality result.

\subsection{Inspiration and Relation to Other Work}

Underlying our results is a mathematical insight concerning the
central role for the analysis of hyperbolic differential equations,
played by the {\em parabolic scaling}, in which analysis elements are
supported in elongated regions obeying the relation $width \approx
length^2$. In fact, curvelets imply the same tiling of the frequency
plane as the Second Dyadic Decomposition (SDD), a construction
introduced in the seventies by Stein and Fefferman
\cite{Fef,Ste}, originally for the purpose of
understanding boundedness of Riesz spherical means, and later widely
adapted to the study of various Fourier Integral Operators. More
specifically, we would like to single out the work of Hart Smith with
which we became familiar while working on this project.  Smith
\cite{Sm1} used parabolic scaling to define function spaces preserved
by Fourier Integral Operators \cite{Sm1}, and to analyze the behavior
of wave equations with low-regularity coefficients \cite{SmithWave}. The
latter reference actually develops curvelet-like systems which provide
a powerful tool to derive so-called sharp Strichartz estimates for
solutions to such equations in space dimensions $d = 2, 3$ (a
Strichartz estimate is a bound on the norm of the solution in some
appropriate functional space, e.g.  $L^p$). We find the connection
with the work of Smith especially stimulating. From a broader
viewpoint, the literature on the subject indicates that curvelets are
in some sense compatible with a long tradition in harmonic analysis.

The fact that the action of a FIO should be seen as a `rearrangement
of wave packets' was discovered by C\'{o}rdoba and Fefferman in their
visionary paper \cite{CF}. They show how simple proofs of $L^2$
boundedness and the Garding inequality follow in a straightforward way
from the decomposition into wave packets. Without rotations and the
parabolic scaling, however, their estimates are not sharp and do not
completely capture the geometry of FIO's.

Next, there is of course the inspiration of modern computational
harmonic analysis (CHA) whose agenda is the development of orthobases,
tight-frames, which are `optimal' for representing objects (operators,
functions) of scientific interest together with rapid algorithms to
compute such representations. The point of view here is to develop new
mathematical ideas and and turn these ideas into effective algorithms
and these effective algorithms into effective and targeted
applications. At the beginning of this introduction, we mentioned an
instance of this scientific vision: (1) wavelets provide sparse
representations of objects with punctuated smoothness and of large
classes of singular integrals and other pseudo-differential operators
\cite{Mey1, Mey2}; (2) there are fast discrete wavelet transforms
operating in $O(N)$ for a signal of length $N$ \cite{MallatBook}; (3)
this creates an opportunity for targeted applications in signal
processing where wavelets allow for better compression
\cite{Donoho-IEEE}, scientific computing where wavelets allow for
faster algorithms \cite{BCR}, and for statistical estimation where
wavelets allow for sharper reconstructions \cite{DonohoSoft}. This
vision was perhaps championed in \cite{BCR} where (1)--(3) were
combined to demonstrate how one can use the wavelet transform to
compute certain types of singular integrals in a number of operations
of the order of $C(\epsilon) \cdot N \log N$ where $C(\epsilon)$ is a
constant depending upon the desired accuracy $\epsilon$.

\subsection{Significance}

We would like to mention how we see our work fit with the vision
described above.

\begin{itemize}
\item {\em Curvelets and wavefronts.} Curvelets are ideal for
  representing wavefront phenomena \cite{DWTPSI}, or objects which
  display {\it curve-punctuated smoothness} ---smoothness except for
  discontinuity along a general curve with bounded curvature
  \cite{Curvelets-StMalo,CurveEdge}.  This fact originally motivated
  their construction \cite{Curvelets-StMalo,CurveEdge}.  For example,
  \cite{CurveEdge} established that curvelets provide the sparsest
  representations of functions which are $C^2$ away from piecewise
  $C^2$ edges. Such representations are nearly as sparse as if the
  object were not singular and turn out to be far more sparse than the
  wavelet decomposition of the object.
  
  Hence, we see that curvelets provide the unique opportunity for
  having a representation giving enhanced sparsity of wave groups, and
  {\em simultaneously} of the solution space.  We believe that this
  will eventually be of great practical significance for applications
  in fields which are great consumers of these mathematical models,
  e.g.  seismics.
  
\item {\em New ideas for new numerical solvers.} Clearly, Theorem
  \ref{teo:main} may serve as a basis for faster geometric multiscale
  PDE solvers. In fact, this paper is the first of a projected series
  showing how one can exploit the structure of the curvelet transform
  and the enhanced sparsity of wave groups to derive new numerical
  low-complexity algorithms for accurately computing the solution to
  large classes of differential equations, see the concluding section
  for a discussion.
  
\item {\em Digital curvelet transforms.} In order to deploy
  curvelet-like ideas in practical applications, one would need a
  digital notion of curvelet transform which (1) would be rapidly
  computable and (2) would be geometrically faithful in the sense that
  one would want an accurate digital analog of the corresponding
  geometric ideas defined at the level of the continuum. There
  actually is progress on this front. Donoho and one of the authors
  recently proposed a Digital Curvelet Transform via Unequispaced Fast
  Fourier Transforms (DCTvUSFFT) \cite{DCTvUSFFT}.  The DCTvUSFFT is a
  fast algorithm which allows analysis and synthesis of Cartesian
  arrays as superpositions of discrete curvelets; for practical
  purposes, the algorithm runs in $O(N\log N)$ operations for input
  array of size $N$. Digital curvelets obey sharp frequency and
  spatial localization.
\end{itemize}
In short, this paper is an essential piece of a much larger body of
work.

\subsection{Contents}

Section 2 below reviews the construction of Curvelets.  Section 3
below examines second-order scalar hyperbolic equations and gives a
heuristic indicating why the sparsity may be expected to hold. Section
4 links our main result with properties of FIOs. Section 5 proves that
FIO's are optimally sparse in scalar curvelet tight-frames. Section 6
discusses implications of this work, namely, in the area of scientific
computing. Finally, proofs of key estimates supporting our main result
are given in Section 7.

\section{Curvelets}
\label{sec:curvelets}
\setcounter{theorem}{0}
\setcounter{equation}{0}

This section briefly introduces tight frames of curvelets, see
\cite{CurveEdge} for more details. 

\subsection{Definition}

We work throughout in $\R^2$, with spatial variable $x$, with $\xi$ a
frequency-domain variable, and with $r$ and $\theta$ polar coordinates
in the frequency-domain.  We start with a pair of windows $W(r)$ and
$V(t)$, which we will call the `radial window' and `angular window',
respectively. These are both smooth, nonnegative and real-valued, with
$W$ taking positive real arguments and supported on $r \in [1/2,2]$
and $V$ taking real arguments and supported on $t \in [-1,1]$. These
windows will always obey the admissibility conditions:
\begin{equation} \label{def-W-disc}
    \sum_{j=-\infty}^\infty W^2(2^j r) = 1, \qquad r > 0;
\end{equation}
\begin{equation} \label{def-V-disc}
    \sum_{\ell=-\infty}^\infty V^2(t - \ell) = 1, \qquad t \in \R.
\end{equation}
Now, for each $j \ge j_0$, we introduce $\vf_j(x_1,x_2)$ defined in
the Fourier domain by
\begin{equation}
  \label{eq:vfj}
  \hat{\vf}_j(\xi) = 2^{-3j/4} W(2^{-j}|\xi|)  \cdot 
V(2^{\lfloor j/2 \rfloor} \theta)  
\end{equation}
Thus the support of $\hat{\vf_j}$ is a polar `wedge' defined by the
support of $W$ and $V$, the radial and angular windows, applied with
scale-dependent window widths in each direction.  

We may think of $\vf_j$ as a ``mother'' curvelet at scale $2^{-j}$ in
the sense that all curvelets at that scaled are obtained by rotations
and translations of $\vf_j$. Introduce
\begin{itemize}
\item the equispaced sequence of {\em rotation angles}
  $\theta_{j,\ell} = 2\pi \cdot 2^{-\lfloor j/2 \rfloor} \cdot \ell$,
  $0 \le \ell < L_j = 2^{\lfloor j/2 \rfloor}$,
\item and the sequence of {\em translation parameters} $k = (k_1,
  k_2) \in \Z^2$.
\end{itemize}
With these notations, we define curvelets (as function of $x = (x_1,
x_2)$) at scale $2^{-j}$, orientation $\theta_{j,\ell}$ and position
$b^{(j,\ell)}_k = R_{\theta_{j,\ell}}(k_1 \cdot 2^{-j}/\delta_1, k_2 \cdot
2^{-j/2}/\delta_2)$ for some adequate constants $\delta_1, \delta_2$ by
\[
\vf_{j,k,\ell}(x) = 
\vf_j\left(R_{-\theta_{j,\ell}}(x - b_k^{(j,\ell)})\right). 
\] 

As in wavelet theory, we also have coarse scale elements. We introduce
the low-pass window $W_0$ obeying
\[
|W_0(r)|^2 + \sum_{j \ge 0} |W(2^{-j} r)|^2 = 1,  
\]
and for $k_1, k_2 \in \bZ$, define coarse scale curvelets as
\[
\Phi_{j_0,k}(x) = \Phi_{j_0}(x - 2^{-j_0} k), \quad 
\hat{\Phi}_{j_0}(\xi) = 2^{-j_0} W_0(2^{-j_0} |\xi|). 
\]
Hence, coarse scale curvelets are nondirectional. The `full' curvelet
transform consists of the fine-scale directional elements
$(\vf_\jlk)_{j \ge j_0, \ell, k}$ and of the coarse-scale isotropic
father wavelets $(\Phi_{j_0,k})_{k}$. For our purposes, it is the behavior
of the fine-scale directional elements that matters.

In the remainder of the paper, we will use the generic notation
$(\vf_\mu)_{\mu \in M}$ to index the elements of the curvelet tight-frame. The dyadic-parabolic subscript $\mu$ stands for the triplet $(j,k,\ell)$. We will also make use of the convenient notations
\begin{itemize}
\item $x_\mu = b_k^{(j,\ell)}$ is the center of $\vf_\mu$ in space.
\item $\theta_\mu = \theta_{j,\ell}$ is the orientation of $\vf_\mu$ with respect to the vertical axis in $x$.
\item $\xi_\mu = (2^j \cos\theta_\mu, 2^j \sin\theta_\mu)$ is the center of $\hat{\vf}_\mu$ in frequency.
\item $e_\mu = \xi_\mu/|\xi_\mu|$ indicates the codirection of $\vf_\mu$.
\end{itemize}

\subsection{Properties}

\begin{figure}
\includegraphics[width=8cm]{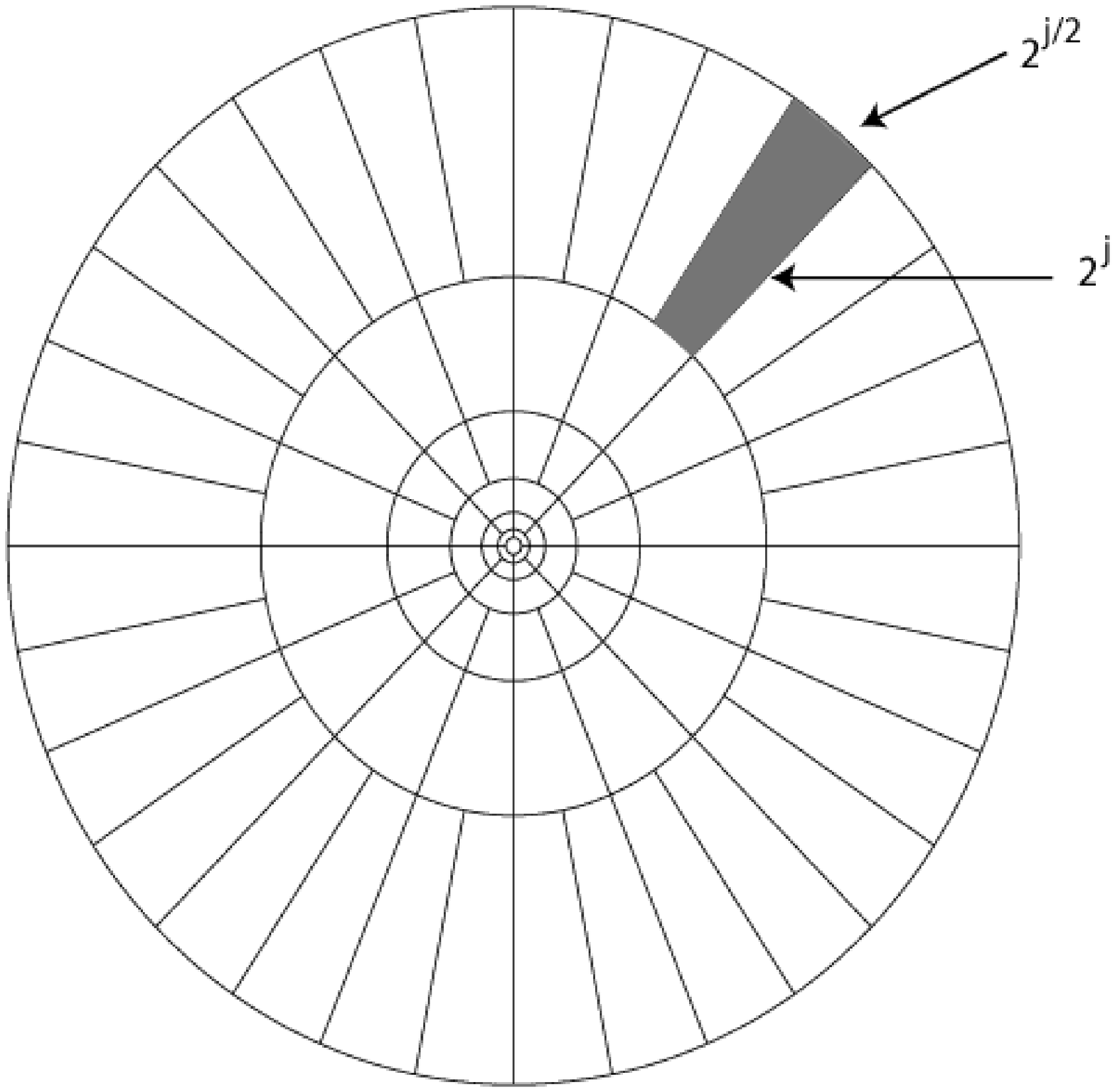}\quad
\includegraphics[width=8cm]{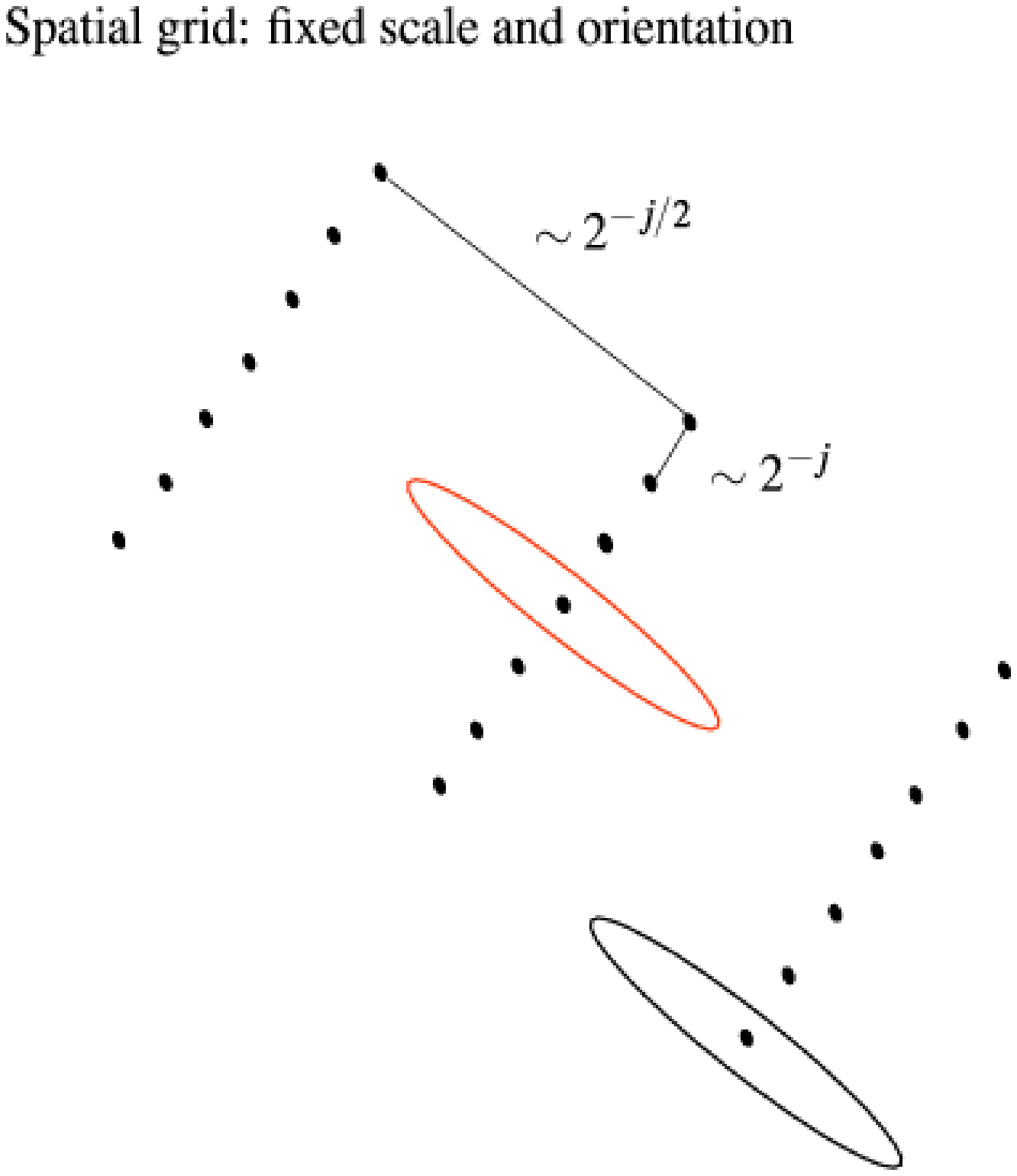}
\caption{Curvelet tiling of Phase-Space. 
  The figure on the left represents the sampling in the frequency
  plane.  In the frequency domain, curvelets are supported near a
  `parabolic' wedge.  The shaded area represents such a generic wedge.
  The figure on the right schematically represents the spatial
  Cartesian grid associated with a given scale and orientation.}
  \label{fig:phase-space} 
\end{figure}

We now list a few properties of the curvelet transform which will play
an important role throughout the remainder of this paper. 
\begin{enumerate}
\item {\bf Tight-frame}. Much like in an orthonormal basis, we can easily
  expand  an arbitrary function $f(x_1, x_2) \in L^2(\R^2)$ as a
  series of curvelets: we have a reconstruction formula 
\begin{equation*}
  \label{eq:tight-frame_sec2}
  f = \sum_{\mu} \<f,\vf_{\mu}\>  \vf_{\mu}, 
\end{equation*}
with equality holding in an $L^2$ sense; and a Parseval relation 
\begin{equation*}
\label{eq:Parseval2}
\sum_{\mu} |\<f,\vf_{\mu}\>|^2  = 
\|f\|_{L^2(\R^2)}^2 , \quad \forall f \in L^2(\R^2).  
\end{equation*}
\item {\bf Parabolic scaling}. The frequency localization of $\vf_j$
  implies the following spatial structure: $\vf_j(x)$ is of rapid
  decay away from an $2^{-j}$ by $2^{-j/2}$ rectangle with minor axis
  pointing in the horizontal direction. In short, the effective length
  and width obey the anisotropy scaling relation
  \begin{equation}
    \label{eq:anisotropy}
length \approx 2^{-j/2}, \quad width \approx 2^{-j} \quad 
\Rightarrow  \quad width \approx length^2.
  \end{equation}
\item {\bf Oscillatory behavior}. As is apparent from its definition,
  $\hat{\vf}_j$ is actually supported away from the vertical axis
  $\xi_1 = 0$ but near the horizontal $\xi_2 = 0$ axis. In a nutshell,
  this says that $\vf_j(x)$ is oscillatory in the $x_1$-direction and
  lowpass in the $x_2$-direction. Hence, at scale $2^{-j}$, a curvelet
  is a little needle whose envelope is a specified `ridge' of
  effective length $2^{-j/2}$ and width $2^{-j}$, and which displays
  an oscillatory behavior across the main `ridge'.
\item {\bf Phase-Space Tiling/Sampling}. We can really think about
  curvelets as Heisenberg tiles of minimum volume in phase-space. In
  $x$, the essential support of $\vf_\mu$ has size $O(2^{-j} \times
  2^{-j/2})$. In frequency, the support of $\hat{\vf}_\mu$ has size
  $O(2^{j/2} \times 2^{j})$. The net volume in phase-space is
  therefore
\[
O(2^{-j} \times 2^{-j/2}) \cdot O(2^{j/2} \times 2^{j}) = O(1),
\]
which is in accordance with the uncertainty principle. The parameters
$(j,k,\ell)$ of the curvelet transform induce a new non-trivial
sampling of phase-space, Cartesian in $x$, polar in $\xi$, and based
on the parabolic scaling.
\item{\bf Complex-valuedness}. Since curvelets do not obey the
  symmetry $\hat \vf_\mu(-\xi) = \overline{\hat \vf_\mu(\xi)}$,
  $\vf_\mu$ is complex-valued. There exists a related construction for
  real-valued curvelets by simply symmetrizing the construction, see
  \cite{CurveEdge}.  The complex-valued transform is better adapted to
  the purpose of this paper.
\end{enumerate}
Figure \ref{fig:phase-space} summarizes the key components of the
construction. 

%
%
%

\subsection{Curvelet Molecules}\label{sec:molecules}

We introduce the notion of {\em curvelet molecule}; our objective,
here, is to encompass under this name a wide collection of systems
which share the same essential properties as the curvelets we have
just introduced. Our formulation is inspired by the notion of
`vaguelettes' in wavelet analysis \cite{Mey2}.  Our motivation for
introducing this concept is the fact that operators of interest do not
map curvelets into curvelets, but rather into these molecules. Note
that the terminology `molecule' is somewhat standard in the literature
of harmonic analysis \cite{FJW}.

\begin{definition} A family of functions 
  $(m_{\mu})_\mu$ is said to be a family of curvelet molecules with
  regularity $R$ if (for $j > 0$) they may be expressed as
\[
m_{\mu}(x) = 2^{3j/4} a^{(\mu)}\left(D_{2^{-j}} R_{\theta_\mu} x - k' \right),
\]
where $k' = (\frac{k_1}{\delta_1},\frac{k_2}{\delta_2})$ and where for all
$\mu$, the $a^{(\mu)}$'s verify the following properties:
\begin{itemize}
\item Smoothness and spatial localization: for each $|\beta| \le R$,
  and each $M = 0, 1, 2 \ldots$, there is a constant $C_M > 0$ such
  that
  \begin{equation}
    \label{eq:localization}
    |\partial_{x}^{\beta} a^{(\mu)}(x)| \leq C_{M} \cdot (1 + |x|)^{-M}.
  \end{equation}
\item Nearly vanishing moments: for each $N = 0, 1, \ldots, R$, there
  is a constant $C_N > 0$ such that
  \begin{equation}
    \label{eq:bandpass}
    |\hat{a}^{(\mu)}(\xi)| \leq 
C_N \cdot  \min (1, 2^{-j} + |\xi_1| + 2^{-j/2} |\xi_2|)^{N}. 
  \end{equation}
\end{itemize}
Here, the constants may be chosen independently of $\mu$ so that the
above inequalities hold uniformly over $\mu$.  There is of course an
obvious modification for the coarse scale molecules which are of the
form $a^{(\mu)}(x - k')$ with $a^{(\mu)}$ as in
(\ref{eq:localization}).
\end{definition}


This definition implies a series of useful estimates. For instance,
consider $\theta_\mu = 0$ so that $R_{\theta_\mu}$ is the identity
(arbitrary molecules are obtained by rotations). Then, $m_\mu$ obeys
\begin{equation}\label{E:locspace}
|m_{\mu}(x)| \leq C_M \cdot 2^{3j/4} \cdot 
\left(1 + |2^j x_1 - \frac{k_1}{\delta_1}| + |2^{j/2} x_2 - \frac{k_2}{\delta_2}|\right)^{-M}
\end{equation}
for each $M > 0$ and $|\beta| \le R$, and similarly for its derivatives
\begin{equation}\label{E:smospace}
|\partial_x^{\beta} m_{\mu}(x)| \leq C_M \cdot   
2^{3j/4} \cdot 2^{(\beta_1 + \beta_2/2)j} 
\cdot \left(1 + |2^j x_1 - \frac{k_1}{\delta_1}| + |2^{j/2} x_2 - \frac{k_2}{\delta_2}|\right)^{-M}. 
\end{equation}
Another useful property is the almost vanishing moments property which
says that in the frequency plane, a molecule is localized near the
dyadic corona $\{2^j \le |\xi| \le 2^{j+1}\}$; $|\hat{m}_{\mu}(\xi)|$
obeys
\begin{equation}
  \label{eq:bp2}
|\hat{m}_{\mu}(\xi)| \leq C_N \cdot  
2^{-3j/4} \cdot \min(1, 2^{-j}(1+|\xi|))^N,
\end{equation}
which is valid for every $N \le R$, which gives the frequency localization 
\begin{equation}
  \label{eq:freq-local}
  |\hat{m}_{\mu}(\xi)| \leq C_N \cdot  2^{-3j/4} \cdot |S_\mu(\xi)|^N
\end{equation}
where for $\mu_0 = (j,0,0)$,
\begin{equation}
\label{eq:Smu0}
S_{\mu_0}(\xi) = \min(1, 2^{-j}(1+|\xi|)) \cdot 
(1 + |2^{-j} \xi_1| + |2^{-j/2} \xi_2|)^{-1}. 
\end{equation}
For arbitrary $\mu$, $S_\mu$ is obtained from $S_{\mu_0}$ by a simple
rotation of angle $\theta_\mu$, i.e.  $S_{\mu_0}(R_{\theta_\mu} \xi)$.
Similar estimates are available for the derivatives of
$\hat{\vf}_\mu$.

In short, a curvelet molecule is a needle whose envelope is supported
near a ridge of length about $2^{-j/2}$ and width $2^{-j}$ and which
displays an oscillatory behavior across the ridge. It is easy to show
that curvelets as introduced in the previous section are indeed
curvelet molecules for arbitrary degrees $R$ of regularity.

\subsection{Near Orthogonality of Curvelet Molecules}

Curvelets are not necessarily orthogonal to each other\footnote{It is
  an open problem whether orthobases of curvelets exist or not.}, but
in some sense they are almost orthogonal. As we show below, the inner
product between two molecules $m_\mu$ and $p_{\mu'}$ decays nearly
exponentially as a function of the `distance' between the subscripts
$\mu$ and $\mu'$.

This notion of distance in phase-space, tailored to curvelet analysis,
is to be understood as follows. Given a pair of indices
$\mu=(j,k,\ell), \mu'=(j',k',\ell')$, define the
\emph{dyadic-parabolic pseudo-distance}
\begin{equation}\label{eq:omega}
\omega(\mu, \mu') = 2^{|j - j'|} \cdot 
\left(1+\min(2^{j},2^{j'}) \, d(\mu, \mu') \right),
\end{equation}
where 
\[
d(\mu,\mu') = |\theta_\mu - \theta_{\mu'}|^2 + |x_\mu - x_{\mu'}|^2 + 
|\< e_\mu, x_\mu - x_{\mu'} \>|.
\]
Angle differences like $\theta_\mu - \theta_{\mu'}$ are understood
modulo $\pi$. As introduced earlier, $e_\mu$ is the codirection of the
first molecule, i.e.,  $e_\mu = (\cos \theta_\mu, \sin \theta_\mu)$.

The pseudo-distance (\ref{eq:omega}) is a slight variation on that
introduced by Smith \cite{SmithWave}. We see that $\omega$ increases
by at most a constant factor every time the distance between the scale,
angular, and location parameters increases. The extension of the
definition of $\omega$ to arbitrary points $(x,\xi)$ and $(x', \xi')$
is straightforward. Observe that the extra term $|\<
e_\mu, x_\mu - x_{\mu'} \>|$ induces a
non-Euclidean notion of distance between $x_{\mu}$ and $x_{\mu'}$. The
following properties of $\omega$ are proved in section
\ref{sec:appendix_curvelets}. (The notation $A \asymp B$ means that $C_1 \le
  A/B \le C_2$ for some constants $C_1, C_2 > 0$.)
\begin{proposition}\label{teo:omega}

\begin{enumerate}
\item Symmetry: $\omega(\mu,\mu') \asymp \omega(\mu',\mu)$.
\item Triangle inequality: $d(\mu,\mu') \leq C \cdot (d(\mu,\mu'') +
  d(\mu'',\mu'))$ for some constant $C > 0$.
\item Composition: for every integer $N > 0$, and some positive constant 
$C_N$ 
\[
\sum_{\mu''} \omega(\mu,\mu'')^{-N} \cdot \omega(\mu'', \mu')^{-N}
\leq C_N \cdot \omega(\mu, \mu')^{-(N-1)}.
\]
\item Invariance under Hamiltonian flows: $\omega(\mu,\mu') \asymp
  \omega(\mu_\nu(t),\mu'_\nu(t))$.
\end{enumerate}
\end{proposition}

We can now state the almost orthogonality result

\begin{lemma}
\label{teo:almost_orthogonality}
Let $(m_{\mu})_{\mu}$ and $(p_{\mu'})_{\mu'}$ be two families of
curvelet molecules with regularity $R$. Then for $j, j' \ge 0$,
\begin{equation}\label{eq:decay}
|\< m_{\mu}, p_{\mu'} \>| \leq C_N \cdot \omega(\mu,\mu')^{-N}.
\end{equation}
for every $N \le f(R)$ where $f(R)$ goes to infinity as $R$ goes to
infinity.
\end{lemma}
\begin{proof}
  Throughout the proof of (\ref{eq:decay}), it will be useful to keep
  in mind that $A \leq C \cdot (1+|B|)^{-M}$ for every $M \leq 2M'$ is
  equivalent to $A \leq C \cdot (1+B^2)^{-M}$ for every $M \leq M'$.
  Similarly, if $A \leq C \cdot (1+|B_1|)^{-M}$ and $A \leq C \cdot
  (1+|B_2|)^{-M}$ for every $M \leq 2M'$, then $A \leq C \cdot
  (1+|B_1| + |B_2|)^{-M}$ for every $M \leq M'$. Here and throughout,
  the constants $C$ may vary from expression to expression.
   
  For notational convenience put $\Delta \theta = \theta_\mu -
  \theta_{\mu'}$ and $\Delta x = x_{\mu} - x_{\mu'}$.  We abuse
  notation by letting $m_{\mu_0}$ be the molecule
  $a^{(\mu)}(D_{2^{-j}} R_{\theta_\mu} x)$, i.e., $m_{\mu_0}$ is
  obtained from $m_\mu$ by translation to that it is centered near the
  origin.  Put $I_{\mu\mu'} = \< m_{\mu}, p_{\mu'} \>$. In the
  frequency domain, $I_{\mu\mu'}$ is given by
\[
I_{\mu\mu'} = \frac{1}{(2 \pi)^2} \int \hat{m}_{\mu_0}(\xi) \overline{\hat{p}_{\mu'_0}(\xi)} \,
e^{-i (\Delta x) \cdot \xi} \, d\xi.
\]
Put $j_0$ to be the minimum of $j$ and $j'$. The Appendix shows that 
\begin{equation}
  \label{eq:freq-loc}
   \int | S_{\mu_0}(\xi) \, S_{\mu'_0}(\xi)|^N \, d\xi \le 
C \cdot  2^{3j/4 + 3j'/4} \cdot 
 2^{-|j - j'| N} \cdot \left(1 + 2^{j_0} 
|\Delta \theta|^2\right)^{-N},
\end{equation}
where $S_{\mu_0}$ is defined in equation (\ref{eq:Smu0}).  Therefore,
the frequency localization of the curvelet molecules
(\ref{eq:freq-local}) gives
\begin{eqnarray}
\nonumber
\label{E:decay-angle}
\int |\hat{m}_{\mu_0}(\xi)| \, |\hat{p}_{\mu'_0}(\xi)| \, d\xi & \le & 
C \cdot 2^{-3j/4 - 3j'/4} \cdot  \int |S_{\mu_0}(\xi) \, S_{\mu'_0}(\xi)|^N \, d\xi \\
& \leq & C \cdot 2^{-|j - j'| N} \cdot \left(1 + 2^{j_0} 
|\Delta \theta|^2\right)^{-N}. 
\end{eqnarray}
This inequality explains the angular decay. A series of integrations by
parts will introduce the spatial decay, as we now show.

The partial derivatives of $\hat{m}_{\mu}$ obey 
\[
|\partial_{\xi}^{\alpha} \hat{m}_{\mu}(\xi)| \le C \cdot 2^{-3j/4} \cdot
2^{-j (\alpha_1 + \frac{\alpha_2}{2})} \cdot |S_\mu(\xi)|^N.
\]
Put $\Delta_{\xi}$ to be the Laplacian in $\xi$. Because $\hat{p}_{\mu'}$ is
misoriented with respect to $e_\mu$, simple calculations show that 
\begin{align*}
  |\Delta_{\xi} \hat{p}_{\mu'}(\xi)| & \le C \cdot 2^{-3j'/4} \cdot
  2^{-j'} \cdot
  |S_{\mu'}(\xi)|^N, \\
  |\frac{\partial^2}{\partial \xi_1^2} \hat{p}_{\mu'}(\xi)| & \le C
  \cdot 2^{-3j'/4} \cdot (2^{-2j'} + 2^{-j'}|\sin(\Delta \theta)|^2)
  \cdot |S_{\mu'}(\xi)|^N.
\end{align*}
Recall that for $t \in [-\pi/2 , \pi/2]$, $2/\pi \cdot |t| \le |\sin
t| \le |t|$, so we may just as well replace $|\sin(\Delta \theta)|$ by
$|\Delta \theta|$ in the above inequality.
Set 
\[
L = I - 2^{j_0} \Delta_{\xi} - \frac{2^{2 j_0}}{1+2^{j_0}|\Delta
  \theta|^2} \frac{\partial^2}{\partial \xi_1^2},
\]
On the one hand, for each $k$, 
$L^k (\hat{m}_{\mu} \overline{\hat{p}_{\mu'}})$ obeys 
\[
|L^k (\hat{m}_{\mu} \overline{\hat{p}_{\mu'}})(\xi)| \le 
C \cdot 2^{-3j/4 - 3j'/4} \cdot |S_{\mu}(\xi)|^N \cdot |S_{\mu'}(\xi)|^N.
\]
On the other hand
\[
L^k e^{- i (\Delta x) \cdot \xi} = [1 + 2^{j_0} |\Delta
x|^2 + \frac{2^{2 j_0}}{1+2^{j_0}|\Delta \theta|^2} |\< e_\mu,
\Delta x \>|^2]^k e^{- i (\Delta x) \cdot \xi}.
\]
Therefore, a few integrations by parts give
\begin{equation*}\label{eq:decay1}
|I_{\mu\mu'}| \leq 
C \cdot 2^{-|j - j'| N} \cdot \left(1 + 2^{j_0} 
|\theta_\mu - \theta_{\mu'}|^2\right)^{-N} \cdot 
\left(1 + 2^{j_0} |\Delta x|^2 + 
\frac{2^{2 j_0}}{1+2^{j_0}|\Delta \theta|^2} 
|\< e_\mu, \Delta x \>|^2\right)^{-N}, 
\end{equation*}
and then 
\[
|I_{\mu\mu'}| \leq C \cdot 2^{-|j - j'| M} \cdot \left(1 + 2^{j_0}
  (|\Delta \theta|^2 + |\Delta x|^2) +
  \frac{2^{2j_0}}{1+2^{j_0}|\Delta \theta|^2} |\< e_\mu, \Delta x
  \>|^2\right)^{-N}.
\]
One can simplify this expression by noticing that
\[
(1 + 2^{j_0} |\Delta \theta|^2) + \frac{2^{2j_0} |\< e_\mu, \Delta x
  \>|^2}{1+2^{j_0}|\Delta \theta|^2} \gtrsim \sqrt{1 + 2^{j_0} |\Delta
  \theta|^2} \frac{2^{j_0} |\< e_\mu, \Delta x
  \>|}{\sqrt{1+2^{j_0}|\Delta \theta|^2}} = 2^{j_0} |\< e_\mu, \Delta x
\>|.
\]
This yields equation (\ref{eq:decay}) as required.
\end{proof}

{\em Remark}. Assume that one of the two terms or both terms are
coarse scale molecules, e.g. $p_{\mu'}$, then the decay estimate is of
the form 
\[
|\< m_{\mu}, p_{\mu'} \>| \leq C \cdot 2^{-j N} \cdot \left(1 +
  |x_{\mu} - x_{\mu'}|^2 + |\<e_\mu, x_{\mu} - x_{\mu'} \>|\right)^{-N}.
\]
For instance, if they are both coarse scale molecules, this would give
\[
|\< m_{\mu}, p_{\mu'} \>| \leq C \cdot 
\left(1 +  |x_{\mu} - x_{\mu'}|\right)^{-N}.
\]

The following result is a different expression for the
almost-orthogonality, and will be at the heart of the sparsity
estimates for FIO's.
\begin{lemma} 
\label{teo:summable}
Let $(m_{\mu})_{\mu}$ and $(p_{\mu})_{\mu}$ be two families of
  curvelet molecules with regularity $R$. Then for each $p > p^*$,
\[
\sup_\mu \sum_{\mu'} |\< m_{\mu}, p_{\mu'} \>|^p \le C_p. 
\] 
Here $p^* \goto 0$ as $R \goto \infty$.  In other words, for $p >
p^*$, the matrix $I_{\mu\mu'} = (\< m_{\mu}, p_{\mu'} \>)_{\mu, \mu'}$
acting on sequences $(\alpha_\mu)$ obeys
\[
\|I \alpha\|_{\ell_p} \le C_p \cdot \|\alpha\|_{\ell_p}. 
\]
\end{lemma}
\begin{proof}
Put as before $j_0 = \min(j,j')$. The appendix shows that 
\begin{equation}\label{eq:schurtest}
\sum_{\mu \in M_{j'}} 
\left(1 + 2^{j_0}(d(\mu, \mu')\right)^{-Np} \le C \cdot 2^{2|j - j'|}
\end{equation}
provided that $Np > 2$. We then have 
\[
\sum_{\mu'} |I_{\mu\mu'}|^{p} \le C \cdot \sum_{j' \in \mathbb{Z}}
2^{-2|j - j'|Np} \cdot 2^{2|j - j'|}\le C_p,
\]
provided again that $Np > 2$. 

Hence we proved that for $p \le 1$, $I$ is a bounded operator from
$\ell_p$ to $\ell_p$. We can of course interchange the role of the two
molecules and obtain
\[
\sup_{\mu'} \sum_{\mu} |\< m_{\mu}, p_{\mu'} \>|^p \le C_p. 
\]
For $p = 1$, the above expression says that $I$ is a bounded operator
from $\ell_\infty$ to $\ell_\infty$. By interpolation, we then
conclude that $I$ is a bounded operator from $\ell_p$ to $\ell_p$ for
every $p$.
\end{proof}

\section{Heuristics}
\label{sec:intuition}
\setcounter{theorem}{0}
\setcounter{equation}{0}

We would like to explain why curvelets are special and why we may
expect sparsity by considering the situation
of the classical acoustic wave propagation equation \eqref{eq:wave}
\begin{equation*}
 \frac{\d^2 u}{\d t^2} = c^2(x) \Delta u,   
\end{equation*}
with initial conditions $u(x,0) = u_0(x)$, $\frac{\d u}{\d t}(x,0) =
u_1(x)$.  To picture the action of the solution operator on a
curvelet, consider the case where the velocity is constant, e.g. $c =
1$.  Put $\hat{u}(\xi,t)$ to be the spatial Fourier transform of the
solution at time $t$ (and similarly for $\hat{u}_0, \hat{u}_1$). In
the Fourier domain, the PDE (\ref{eq:wave}) is then transformed into
an ODE whose solution is given by the relation
\begin{equation}
\hat{u}(\xi,t) = \cos (|\xi| t) \, \hat{u}_0(\xi) + 
\frac{\sin (|\xi| t)}{|\xi|} \, \hat{u}_1(\xi).  
\end{equation} 
To understand the sparsity we need only to consider the sparsity of
the multiplication by $ \cos (|\xi| t)$, as the exact same arguments apply to 
the other term. In addition, since 
\[
\cos (|\xi| t) = (e^{i |\xi| t} + e^{-i |\xi| t})/2, 
\]
we may just as well study the sparsity of the the multiplication by
$e^{i |\xi| t}$. We will now explain how the parabolic scaling plays an
essential role.

\subsection{Frequency localization}

Because the frequency support of a curvelet $\vf_\mu$ is effectively a
thin dyadic rectangle of sidelengths about $2^{j}$ in the direction
$e_\mu$ and about $2^{j/2}$ in the orthogonal direction, the modulus
of $|\xi|$ is close to the component of $\xi$ along $e_\mu$, and we
may want to linearize the phase as follows
\[
|\xi| = \xi \cdot e_\mu + (|\xi| - \xi \cdot e_\mu).
\]
The point is that the ``phase perturbation'' $\delta(\xi) = |\xi| -
\xi \cdot e_\mu$ is bounded over rectangles which support $\hat \vf_\mu(\xi)$.
With these notations, the phase multiplication takes the form
\[
e^{i |\xi| t} \hat \vf_\mu(\xi) = e^{i \delta(\xi) t} (e^{i (\xi \cdot e_\mu) t}
\hat \vf_\mu(\xi));
\]
now the multiplier $e^{i \delta(\xi) t}$ is not really oscillatory (at
least for $t$ not too large) thanks to $\delta(\xi) = O(1)$.
Therefore, the phase multiplication amounts to a linear modulation
followed by a multiplication with a smooth term. The frequency picture
is given in figure \ref{fig:phase_linearization}. In the spatial
domain, this simply corresponds to a shift by $t$ along the
co-direction of the curvelet followed by a convolution with a gentle
kernel.


As it seems reasonable to assume that curvelets sparsify both
translations and convolutions, one might expect sparse representations
of the solution operator, at least in the case where the coefficients
of the differential equations are constant.

\subsection{Spatial localization}

On the one hand, curvelets have enough frequency localization so that
they approximately behave like waves. But on the other, they have
enough spatial localization so that the flow will essentially preserve
their shape, even when the velocity is varying.

\begin{figure}
   \begin{minipage}[c]{.58\linewidth}
     \includegraphics[width=\linewidth]{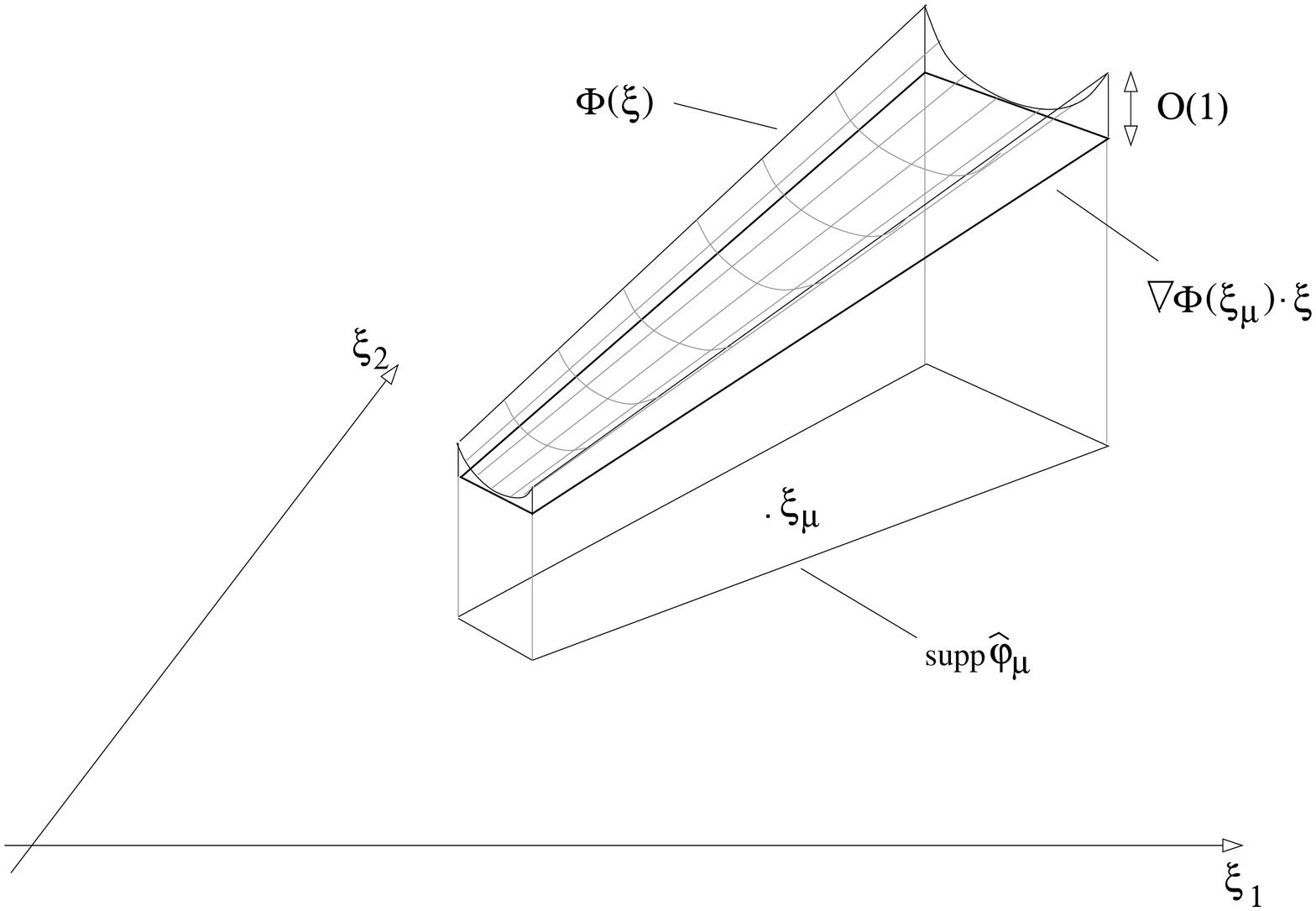}
      \caption{The phase function $\Phi(\xi)$ is well-approximated by its linearization $\nabla \Phi(\xi_\mu) \cdot \xi$ on the support of $\hat{\vf}_\mu$. In the example discussed in the text, $\Phi(\xi) = t|\xi|$.}
      \label{fig:phase_linearization}
   \end{minipage} \hfill
   \begin{minipage}[c]{.34\linewidth}
     \includegraphics[width=\linewidth]{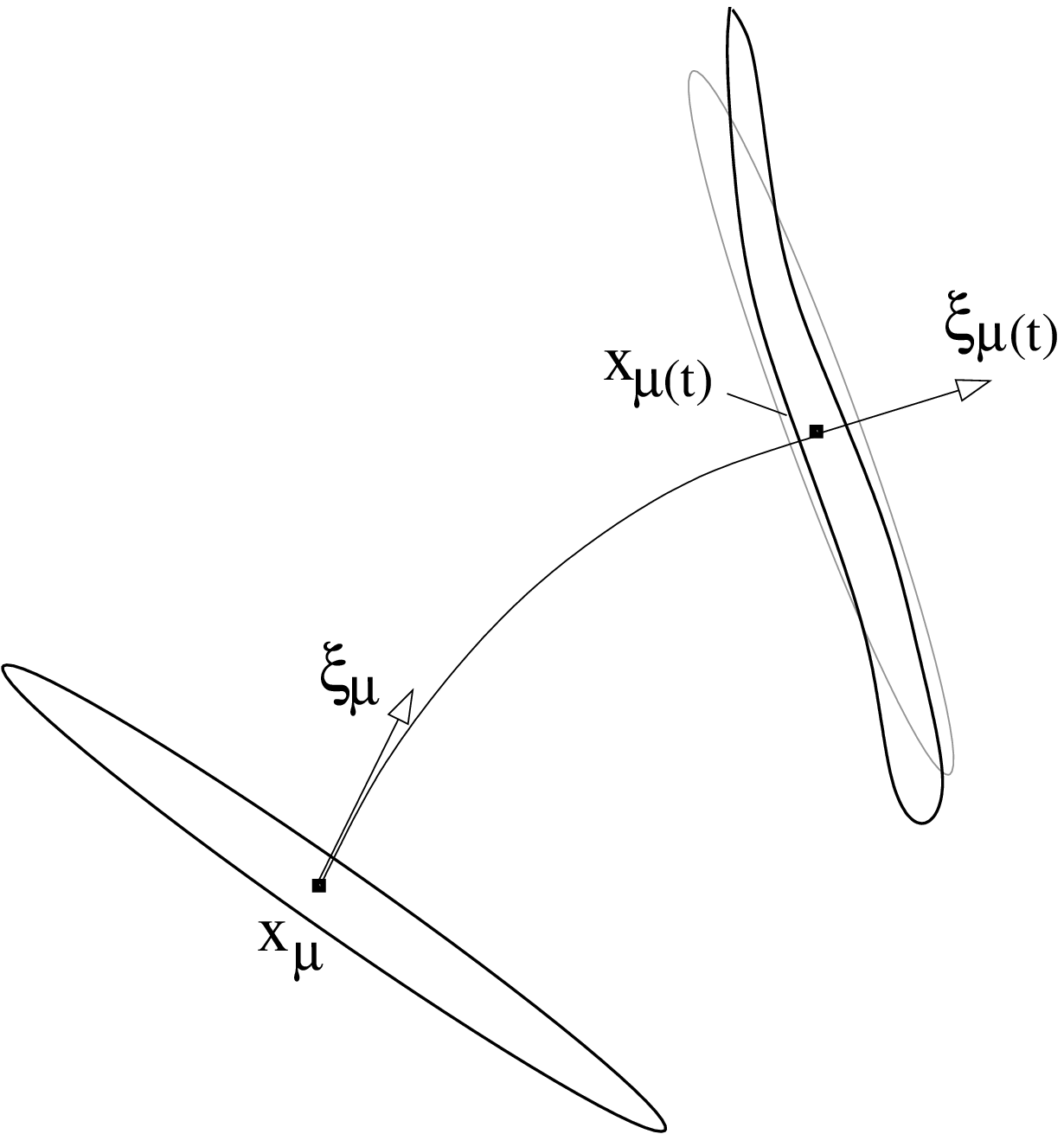}
      \caption{The bottom ellipse indicates the essential support of $\vf_\mu$, initially at $(x_\mu,\xi_\mu)$. After propagation, the essential support of $\vf_{\mu(t)}$ is warped but not by much.}
      \label{fig:harmless_warping}
   \end{minipage}
\end{figure}

At high frequencies or equivalently at small scales, a curvelet has an
envelope strongly aligned along a specified `ridge' of length
$2^{-j/2}$ and width $2^{-j}$ and is oscillatory across the main
`ridge' with a frequency roughly inversely proportional to its width.
Now the key-point is that the velocity field does not vary very much
over the support of a curvelet. In fact, owing to the geometric
relation $width \sim length^2$, the velocity field cannot
significantly bend our curvelet. To see this fact, let us model the
action of the velocity field as a warping $g(x)$. Consider now the
warped curvelet $\vf_\mu(g(x))$. Then the warping leaves the curvelet
almost invariant; i.e. our curvelets approximately remain needles
essentially fitting in boxes of side-lengths $2^{-j/2}$ and $2^{-j}$.
See figure \ref{fig:harmless_warping}.  This behavior is of course
very different from that of classical waves of the form $e^{i x \cdot
  \xi}$.  Because such waves have elongated support, the wavefronts
are geometrically distorted.

\subsection{The parabolic scaling is special} 

There is something truly remarkable about the parabolic scaling.
Instead of curvelets obeying the law
\[
width \approx length^2
\]
one may want to consider the family of element obeying the more
general power-law
\[
width \approx length^\alpha, \quad 1 \le \alpha \le \infty.
\]
However, we have seen that to exhibit the wave-like behavior, one
needs enough frequency localization. In our setup, we need that in the
frequency domain, our waveform lie in a box of width $\ell$ and length
$L$ obeying $\ell^2/L = O(1)$. In the spatial domain, this minimum
anisotropy imposes the relationship
\[
width \le  length^2.
\]
But we cannot afford too much anisotropy as otherwise curvelets would
be distorted. In the spatial domain, this says that curvelets have to
be supported near boxes obeying $L^2 \le \ell$ or in other words
\[
width \ge length^2.
\]
In conclusion, if we want both these properties to hold
simultaneously, namely, wave-like behavior and particle-like behavior
(no essential distortion by applying a smooth warping), only one
scaling may possibly work: the {\em parabolic scaling}.

\section{Representation of Linear Hyperbolic Systems}
\label{sec:lax}
\setcounter{theorem}{0}
\setcounter{equation}{0}

We now return to the main theme of this paper and consider linear
initial-value problems of the form
\begin{equation}\label{eq:hyp}
\frac{\d u}{\d t} + 
\sum_{k = 1}^m A_k(x) \frac{\d u}{\d x_k} + B(x) u = 0, 
\qquad u(0,x) = u_0(x),
\end{equation}
where in addition to the properties listed in the introduction, $A^k$
and $B$ together with all their partial derivatives are uniformly
bounded for $x \in \R^n$.  As explained in section
\ref{sec:decoupling}, we need to make the technical assumption that
for every set of real parameters $\xi_k$, the (real) eigenvalues of
the matrix $\sum_k A_k(x) \xi_k$ have constant multiplicity in $x$ and
$\xi$.

Our goal is to construct a concrete `basis' of $L^2(\R^n, \C^m)$ in
which the evolution is as simple/sparse as possible. We present a
solution based on the newly developed curvelets---which were
introduced in section \ref{sec:curvelets}---and choose to specialize our
discussion to $n = 2$ spatial dimensions. The reason is twofold:
first, this setting is indeed that in which the exposition of
curvelets is the most convenient; and second, this is not a
restriction as similar results would hold in arbitrary dimensions.

\subsection{Main result}

We need to prove
\begin{equation}\label{eq:main_result}
|E(t;\mu,\nu;\mu',\nu')| \leq C_{t,N} \cdot \sum_{\nu''} 
\omega(\mu, \mu'_{\nu''}(t))^{-N},
\end{equation}
for some constant $C_{t,N} > 0$ growing at most like $C_{N}e^{K_{N}t}$
for some $C_N, K_N > 0$. The sum over $\nu''$ indexes the different
flows and takes on as many values as there are distinct eigenvalues
$\lambda^0_{\nu''}$.

It is instructive to notice that the estimate (\ref{eq:main_result}) for $t
= 0$ is already the strongest of its sort on the off-diagonal decay of
the Gram matrix elements for a tight frame of curvelets. For $t > 0$,
equation (\ref{eq:main_result}) states that the strong phase-space
localization of every curvelet is preserved by the hyperbolic system,
thus yielding a sparse and well-organized structure for the curvelet
matrix. These warped and displaced curvelets are `curvelet molecules'
as introduced in section \ref{sec:molecules} because, as we will show,
they obey the estimates (\ref{eq:localization}) and
(\ref{eq:bandpass}).

The choice of the curvelet family being complex-valued in the above
theorem is not essential. $E(t)$ acting on real-valued curvelets would
yield \emph{two} molecules per flow (upstream and downstream).
Keeping track of this fact in subsequent discussions would be
unnecessarily heavy. In the real case, it is clear that the structure
and the sparsity of the curvelet matrix can be recovered by expressing
each real curvelet as a superposition of two complex curvelets.

\subsection{Architecture of the proof of the main result}
\label{sec:architecure}

Much of the remainder of this paper is devoted to the justification of
Theorem \ref{teo:main} or equivalently \eqref{eq:main_result}. However, before
engaging in the formal proof, we would like to present the overall
architecture of the argument.
\begin{itemize}
\item {\em Decoupling into polarized components}. The first step is to
  decouple the wavefield $u(t,x)$ into $m$ one-way components
  $f_\nu(t,x)$ 
\[
u(t,x) = \sum_{\nu = 1}^m R_\nu f_\nu(t,x),
\]
where the $R_\nu$ are operators mapping scalars to $m$-dimensional
vectors, and independent of time. The $f_\nu$ will also be called
`polarized' components. This allows a separate study of the $m$ flows
corresponding to the $m$ eigenvalues of the matrix $\sum_{k=1}^m
A_k(x) \xi_k$. In the event these eigenvalues are simple, the
evolution operator $E(t)$ can be decomposed as
\begin{equation}\label{eq:one-way-decomp}
E(t) = \sum_{\nu=1}^m R_\nu e^{-it \Lambda_\nu} L_\nu + 
\text{ negligible},
\end{equation}
where the $L_\nu$'s are operators mapping $m$-dimensional vectors to
scalars and the $\Lambda_\nu$'s are one-way wave operators acting on
scalar functions. In effect, each operator $E_\nu(t) = e^{-it
  \Lambda_\nu}$ convects wave-fronts and other singularities along a
separate flow. The `negligible' contribution is a smoothing
operator---not necessarily small. The composition operators $R_\nu$
and decomposition operators $L_\nu$ are provably pseudo-differential
operators, see section \ref{sec:decoupling}.

\item {\em Fourier Integral Operator parametrix.} We then approximate
  for small times $t > 0$ each $e^{-it \Lambda_\nu}$, $\nu = 1,
  \ldots, m$, by an oscillatory integral or Fourier Integral Operator  
  (FIO) $F_\nu(t)$. Such operators take the form
\[
F_\nu(t)f(x) = \int e^{i\Phi_\nu(t,x,\xi)} \sigma_\nu(t,x,\xi)
\hat{f}(\xi) \, d\xi,
\]
under suitable conditions on the phase function $\Phi_\nu(t,x,\xi)$
and the amplitude $\sigma_\nu(t,x,\xi)$. Again, the identification of
the evolution operator $E_\nu(t) = e^{-it \Lambda_\nu}$ with $F_\nu$
is valid up to a smoothing and localized additive remainder. The
construction of the so-called {\em parametrix} $F_\nu(t)$ and its
properties are detailed in section \ref{sec:FIOhyp}.

Historically \cite{Lax}, the construction of an oscillatory integral
parametrix did not involve the decoupling into polarized components as
a preliminary step. When applied directly to the system
(\ref{eq:hyp}), the construction of the parametrix gives rise to a
matrix-valued amplitude $\sigma(t,x,\xi)$ where all the couplings are
present. This somewhat simpler setting, however, is not adequate for
our purpose.  The reason is that we want to bootstrap the construction
of a parametrix to \emph{large} times by composing the small time FIO
parametrix with itself, $F(nt) = [F(t)]^n$. Without decoupling of the
propagation modes, each $E(t)$ or $F(t)$ involves convection of
singularities along $m$ families of characteristics or flows. Applying
$F(t)$ again, each flow would artificially split into $m$ flows again,
yielding $m^2$ fronts to keep track of. At time $T = nt$, that would
be at most $m^n$ fronts. This flow-splitting situation is not physical
and can be avoided by isolating one-way components \emph{before}
constructing the parametrix. The correct large-time argument is to
consider $E_\nu(nt)$ for small $t > 0$ and large integer $n$ as
$[E_\nu(t)]^n$. This expression involves one single flow, indexed by
$\nu$.

\item {\em Sparsity of Fourier Integral Operators.} The core of the
  proof is found in section \ref{sec:FIO} and consists in showing that
  very general FIO's $F(t)$, including the parametrices $F_\nu(t)$,
  are sparse and well-structured when represented in tight frames of
  (scalar) curvelets $\vf_\mu$. The scalar analog of Theorem
  \ref{teo:main} for FIO's is Theorem \ref{teo:main_FIO}---a
  statement of independent interest. Observe that pseudo-differential
  operators are a special class of FIO's and, therefore, are equally
  sparse in a curvelet frame.
\end{itemize}
Section \ref{sec:proof_main} assembles key intermediate results and
proves Theorem \ref{teo:main}.

\subsection{Decoupling into polarized components}\label{sec:decoupling}

How to disentangle the vector wavefield into $m$ independent
components is perhaps best understood in the special case of constant
coefficients, $A_k(x) = A_k$, and with $B(x) = 0$. In this case,
applying the 2-dimensional Fourier transform on both sides of
(\ref{eq:hyp}) gives a system of ordinary differential equations
\[
\frac{d \hat{u}}{dt}(t,\xi) + i a(\xi) \hat{u}(t,\xi) = 0, \quad
a(\xi) = \sum_k A_k \xi_k.
\]
(Note that $a(\xi)$ is a symmetric matrix with real entries.) It
follows from our assumptions that one can find $m$ real eigenvalues
$\lambda_\nu(\xi)$ and orthonormal eigenvectors $r_\nu(\xi)$, so that
\[
a(\xi) r_{\nu}(\xi) = \lambda_{\nu}(\xi) r_\nu(\xi).
\]
Put $f_\nu(t,\xi) = r_\nu(\xi) \cdot \hat{u}(t,\xi)$. Then our system
of equations is of course equivalent to the system of independent
scalar equations
\[
\frac{d f_\nu}{dt}(t,\xi) + i \lambda_\nu(\xi) f(t,\xi) = 0, 
\]
which can then be solved for explicitly;
\[
f_\nu(t,\xi) = e^{-it\lambda_\nu(\xi)}f_\nu(0,\xi).
\]
Hence, the diagonalization of $a(\xi)$ decouples the original equation
(\ref{eq:hyp}) into $m$ \emph{polarized} components; these can be
interpreted as waves going in definite directions, for example `up and
down' or `outgoing and incoming' depending on the geometry of the
problem. This is the reason why $f_{\nu}$ is also referred to as being a `one-way' wavefield.

The situation is more complicated when $A_k(x)$ is non-uniform since
Fourier techniques break down. A useful tool in the variable
coefficient setting is the calculus of pseudo-differential operators.
An operator $T$ is said to be {\em pseudo-differential} with symbol
$\sigma$ if it can be represented as
\begin{equation}\label{eq:defpseudo}
Tf(x) = \sigma(x,D) f = \frac{1}{(2 \pi)^2} \int_{\R^2} e^{ix \cdot \xi} 
\sigma(x,\xi) \hat{f}(\xi) \, d\xi,
\end{equation}
with the convention that $D = -i \nabla$. It is of type $(1,0)$ and
order $m$ if $\sigma$ obeys the estimate
\[
|\partial_{\xi}^\alpha \partial_x^\beta \sigma(x,\xi)| \leq
C_{\alpha,\beta} \cdot (1 + |\xi|)^{m-|\alpha|}
\]
for every multi-indices $\alpha$ and $\beta$. Unless otherwise stated,
all pseudo-differential operators in this paper are of type $(1,0)$.
An operator is said to be smoothing of order $-\infty$, or simply
smoothing if its symbol satisfies the above inequality for every $m <
0$. Observe that this is equivalent to the property that $T$ maps
boundedly distributions in the Sobolev space $H^{-s}$ to functions in
$H^s$ for every $s > 0$, in addition to a strong localization property
of its kernel $G(x,y)$ which says that for each $N > 0$, there is a
constant $C_N > 0$ such that $G$ obeys
\begin{equation}\label{eq:localkernel}
|G(x,y)| \leq C_N \cdot (1 + |x-y|)^{-N}
\end{equation}
as in \cite{Ste}[Chapter 6].

Now set 
\[
a(x,D) = \sum_{k = 1}^m A_k(x) D_k - i B(x),
\]
and its principal part 
\[
a^0(x,D) = \sum_{k = 1}^m A_k(x) D_k,
\]
so that equation (\ref{eq:hyp}) becomes $\d_t u + i a(x,D)u = 0$. The
matrices $a(x,\xi)$ (resp.~$a^0(x,\xi)$) are called the symbol of the operator
$a(x,D)$ (resp.~$a^0(x,D)$). Note that $a^0(x,\xi)$ is homogeneous of
degree one in $\xi$; $a^0$ also goes by the name of {\em dispersion
  matrix.}

It follows from the symmetry of $A_k$ and $B$ that for every set of
real parameters $\xi_1, \ldots, \xi_m$, the matrix $a^0(x,\xi) =
\sum_k A_k(x) \xi_k$ is also symmetric and thus admits real
eigenvalues $\lambda^0_\nu(x,\xi)$ and an orthonormal basis of
eigenvectors $r^0_\nu(x,\xi)$,
\begin{equation}\label{eq:orderone}
a^0(x,\xi) r^0_{\nu}(x,\xi) = \lambda^0_{\nu}(x,\xi) r^0_\nu(x,\xi).
\end{equation}
The eigenvalues being real and the set of eigenvectors complete is a
\emph{hyperbolicity} condition and ensures that equation
(\ref{eq:hyp}) will admit wave-like solutions. We assume
throughout this paper that the multiplicity of each
$\lambda^0_{\nu}(x,\xi)$ is constant in $x$ and $\xi$.

By analogy with the special case of constant coefficients, a first
impulse may be to introduce the components $r^0_\nu(x,D) \cdot u$,
where $r^0_\nu(x,D)$ is the operator associated to the eigenvector
$r^0_\nu(x,\xi)$ by the standard rule (\ref{eq:defpseudo}). In
particular this is how we defined hyper-curvelets from curvelets in
section \ref{sec:curvgeomoptics}.  Unfortunately, this does not
perfectly decouple the system into $m$ polarized modes---it only
approximately decouples. Instead, we would achieve perfect decoupling
if we could solve the eigenvalue problem
\begin{equation}\label{eq:twisted_egv}
a(x,D) r_\nu(x,D) = r_\nu(x,D) \lambda_\nu(x,D).
\end{equation}
Here, each $\Lambda_\nu = \lambda_\nu(x,D)$ is a scalar operator and
$R_\nu = r_\nu(x,D)$ is an $m$-by-1 vector of operators. Equation
(\ref{eq:twisted_egv}) must be understood in the sense of composition
of operators. Now let $f_\nu$ be the polarized components obeying the
scalar equation
\begin{equation}\label{eq:evol_fnu}
\frac{\d f_\nu}{\d t} + i \Lambda_\nu f_\nu = 0,
\end{equation}
with initial condition $f_\nu(0,x)$ and consider the superposition
\[
u = \sum_\nu u_\nu, \quad u_\nu = R_\nu f_\nu. 
\]
Then $u$ is a solution to our initial-value problem (\ref{eq:hyp}).
(We will make this rigorous later, and detail the dependence between
the initial values $u_0$ and the $f_\nu(0,\cdot)$'s.)


The following result and its proof show how in some cases,
(\ref{eq:twisted_egv}) can be solved up to a smoothing remainder of
order $-\infty$. When all the eigenvalues $\lambda^0_\nu(x,\xi)$ are
simple, the exact diagonalization is, in fact, possible. The situation
is more complicated when some of the eigenvalues are degenerate. This
complication does not compromise, however, any of our results.
\begin{lemma}\label{teo:twisted_egv}
  Suppose our hyperbolic system satisfies all the assumptions
  stated below (\ref{eq:hyp}). Then there exists an $m$-by-$m$
  block-diagonal matrix of operators $\Lambda$ and two $m$-by-$m$
  matrices of operators $R$ and $S$ such that
\[
a(x,D) R = R \Lambda + S,
\]
where $\Lambda, R$ and $S$ are componentwise pseudo-differential with
$\Lambda$ of order one, $R$ of order zero, and $S$ of order $-\infty$.
Each block of $\Lambda$ corresponds to a distinct eigenvalue
$\lambda^0_\nu$ whose size equals the multiplicity of that eigenvalue. The
principal symbol of $\Lambda$ is diagonal with the eigenvalues
$\lambda^0_\nu(x,\xi)$ as entries.
\end{lemma}
The lemma says that polarized components corresponding to distinct
eigenvalues can be decoupled modulo a smoothing operator. However,
further decoupling within the eigenspaces is in general not possible.
We would like to point out that a similar result was already given in
\cite{Tay} although in a different context and with a different proof,
and was used in $\cite{Stolk_deHoop}$.

\begin{proof}
  We already argued (\ref{eq:twisted_egv}) is not just the eigenvalue
  problem for the symbol $a(x,\xi)$ for the composition of two
  operators does not reduce to a multiplication of their respective
  symbols. Instead, it is common practice \cite{Fol} to define the
  \emph{twisted product} of two symbols $\sigma$ and $\tau$ as
\[
(\sigma \, \sharp \, \tau)(x,D) = \sigma(x,D) \tau(x,D),
\]
so that (\ref{eq:twisted_egv}) becomes the symbol equation $a \,
\sharp \, r_\nu = r_\nu \, \sharp \, \lambda_\nu$. Note that $D = -i
\nabla$. The explicit formula for the twisted product is, in
multi-index notation\footnote{All the pseudo-differential operators
  considered in this paper are of type $(1,0)$ therefore all such
  polyhomogeneous expansions are valid.},
\[
\sigma \, \sharp \, \tau = \sum_{|\alpha| \geq 0} \frac{1}{\alpha !}
\d^\alpha_\xi \sigma D^\alpha_x \tau.
\]
We can see that $\sigma \, \sharp \, \tau$ is the product $\sigma
\tau$ up to terms that are at least one order lower (because of the
differentiations in $\xi$).

Recall the decomposition of the symbol $a(x,\xi)$ into a principal
part $a^0(x,\xi) = \sum_k A_k(x) \xi_k$, homogeneous of degree one in
$\xi$, and a remainder $B(x)$ homogeneous of order zero. It follows
that the eigenvalues $\lambda^0_\nu(x,\xi)$ of $a^0(x,\xi)$ are
homogeneous of degree one, and the corresponding eigenvectors
$r^0_\nu(x,\xi)$ may be selected as homogeneous of degree zero (and
orthonormal). Up to terms of lower order in $\xi$, the original
problem (\ref{eq:twisted_egv}) therefore reduces to the eigenvalue
problem $a^0(x,\xi) r_\nu^0(x,\xi) = r_\nu^0(x,\xi)
\lambda_\nu^0(x,\xi)$ for the symbol $a^0$. It is then natural to look
for a solution $r_\nu$, $\lambda_\nu$ of (\ref{eq:twisted_egv}) as a
perturbation of $r^0_\nu$, $\lambda^0_\nu$ by lower-order terms.

Consider first the case in which each eigenvalue $\lambda^0_\nu$ is
simple and define the expansions
\[
r_\nu \sim r_\nu^0 + r_\nu^1 + r_\nu^2 + \ldots, \qquad \lambda_\nu
\sim \lambda_\nu^0 + \lambda_\nu^1 + \lambda_\nu^2 + \ldots
\]
so that $r_\nu^n$ is of order $-n$ in $\xi$ and $\lambda_\nu^n$ of
order $-n+1$ i.e.,
\[
|\d^\alpha_\xi \d^\beta_x r_\nu^n(x,\xi)| \leq C_{\alpha, \beta} (1 +
|\xi|)^{-n-|\alpha|},
\]
and similarly for $\lambda_\nu^n$. We plug these expansions in the
twisted product, or equivalently in (\ref{eq:twisted_egv}), and
isolate terms of identical degree.

The contribution at the leading order is, of course, $a^0 r^0_\nu =
\lambda^0_\nu r^0_\nu$ and the remainder is of the form $a \, \sharp
\, r_\nu^0 - r_\nu^0 \, \sharp \, \lambda_\nu^0$; put $e^0_\nu$ as its
principal symbol. The zero-order equation reads 
\begin{equation}
\label{eq:appendix_orderzero}
(a^0 - \lambda^0_\nu \,I) r^1_\nu = 
-e^0_\nu + r^0_\nu \lambda^1_\nu, 
\end{equation}
which admits a solution if and only if the right hand side has a zero
component in the eigenspace spanned by $r^0_\nu$. This is possible if
$\lambda^1_\nu$ is selected so that 
\[
-e^0_\nu + r^0_\nu \lambda^1_\nu \,\, \perp \,\, r^0_\nu \quad
\Leftrightarrow \quad \lambda^1_\nu = r^0_\nu \cdot e^0_\nu.
\]
It follows that equation (\ref{eq:appendix_orderzero}) admits the
family of solutions
\[
r^1_\nu = (a^0 - \lambda^0_\nu \, I)^{-1}(-e^0_\nu + r^0_\nu
\lambda^1_\nu) + f_1 r^0_\nu,
\]
where $f_1$ is actually a scalar function of $x$ and $\xi$, and
homogeneous of degree -1 in $\xi$. Our proof does not exploit this
degree of freedom.

It is clear that one can successively determine all the
$\lambda^n_\nu$'s and $r^n_\nu$'s in a similar fashion. Let $e^n_\nu$
be the principal symbol of $a \, \sharp \,(r_\nu^0 + \ldots + r_\nu^n)
- (r_\nu^0 + \ldots + r_\nu^n) \, \sharp \, (\lambda_\nu^0 + \ldots +
\lambda_\nu^n)$, then the equation at the order $-n$ is
\[
(a^0 - \lambda^0_\nu \, I) r^{n+1}_\nu = -e^n_\nu + r^0_\nu
\lambda^{n+1}_\nu, 
\]
and is solved exactly like (\ref{eq:appendix_orderzero}).

Suitable cut-offs of the low frequencies guarantee convergence of the
series for $r_\nu$ and $\lambda_\nu$. As is standard in the theory of
pseudo-differential operators \cite{Sog,Tre}, one selects a sequence
of $C^{\infty}$ cut-off functions $\chi_n(\xi) = \chi(\epsilon^{n} \xi)$ for some $\chi$ vanishing inside a
compact neighborhood of the origin, and identically equal to
one outside a larger neighborhood. Then $\epsilon$ is taken small enough so that
\[
r_\nu(x,\xi) = \sum_{n = 0}^{\infty} r^n_\nu(x,\xi) \chi_n(\xi), \qquad
\lambda_\nu(x,\xi) = \sum_{n = 0}^{\infty} \lambda^n_\nu(x,\xi)
\chi_n(\xi)
\]
are converging expansions in the topology of $C^\infty$. As a result,
the remainder $s_\nu = a \, \sharp \, r_\nu - r_\nu \, \sharp \,
\lambda_\nu$ also converges to a valid symbol which, by construction,
is of order $-\infty$, i.e. obeys
\[
|\d^\alpha_\xi \d^\beta_x s_\nu(x,\xi)| \leq C_{\alpha,\beta,N} \cdot
(1 + |\xi|)^{-N-|\alpha|}
\]
for every $N > 0$. The lemma is proved in the case when all
eigenvalues of the principal symbol are simple.

Consider now the case of a multiple eigenvalue $\lambda^0$, say.
Suppose the corresponding eigenspace is of dimension $p$ and spanned by
$r^0_{1}, \ldots, r^0_{p}$. The reasoning for simple eigenvalues does
not apply because the $p$ solvability conditions are too many for
purely diagonal lower-order corrections. Instead, the block
corresponding to $\lambda^0$ is now perturbed as
\[
\begin{pmatrix} \lambda^0 & \cdots & 0 \\ \vdots & \ddots & \vdots \\ 
0 & \cdots & \lambda_0 \end{pmatrix} + 
\begin{pmatrix} \lambda^1_{11} & \cdots & \lambda^1_{1p} \\ 
\vdots & \ddots & \vdots \\ \lambda^1_{p1} & \cdots & \lambda^1_{pp} 
\end{pmatrix} + \begin{pmatrix} \lambda^2_{11} & \cdots & \lambda^2_{1p}\\ 
\vdots & \ddots & \vdots \\ \lambda^2_{p1} & \cdots & \lambda^2_{pp} 
\end{pmatrix} + \ldots
\]
where each $\lambda^n_{ij}$ is homogeneous of degree $-n+1$ in $\xi$.
At the leading order, The $p$ equations relative to $\lambda_0$ are
\begin{equation}
\label{eq:appendix_orderzero_multiple}
(a^0 - \lambda^0 \, I) r^1_j = 
-e^0_j + \sum_{i = 1}^p r^0_i \lambda^1_{ij},
\end{equation}
where $e^0_j$ is the principal symbol of $a \, \sharp \, r_j^0 - r_j^0
\, \sharp \, \lambda^0$. Solvability requires that the projection of
the right-hand side on each of the $r^0_i$, $i = 1, \ldots, p$
vanishes. This unambiguously determines all the components of the
$p$-by-$p$ block $\lambda^1$ as
\[
\lambda^1_{ij} = r^0_i \cdot e^0_j.
\]
All blocks relative to other eigenvalues are solved for in a similar
way, yielding a block-diagonal structure for the zeroth order
correction $\lambda^1$. Each block should have dimension equal to the
multiplicity of the corresponding eigenvalue in order to meet the
solvability requirements.

The perturbed eigenvectors $r^1_1, \ldots, r^1_p$ are determined as
previously once the $\lambda^1_{ij}$ are known. The same reasoning
applies at all orders and thereby determines $\Lambda$ and $R$.
Convergence issues are addressed using cut-off windows just as
before.
\end{proof}

The above construction indeed provides efficient decoupling of the
original problem (\ref{eq:hyp}) into polarized modes.
\begin{lemma}\label{teo:decoupling}
  In the setting of Lemma \ref{teo:twisted_egv}, the solution operator
  $E(t)$ for (\ref{eq:hyp}) may be decomposed for all times $t > 0$ as
\[
E(t) = R e^{-it \Lambda} L + \tilde{S}(t),
\]
where the matrices of operators $\Lambda$ and $R$ are defined in Lemma
\ref{teo:twisted_egv} and $\tilde{S}(t)$ is (another) matrix of
smoothing operators of order $-\infty$. In addition,
\begin{enumerate}
\item $L$ is an approximate inverse of $R$, i.e. $RL = I$ and $LR = I$
  (mod smoothing). 
\item $L$ is a pseudo-differential of order zero (componentwise).
\end{enumerate}
Observe that $e^{-it \Lambda}$ inherits the block structure from
$\Lambda$, and is diagonal in the case where all the eigenvalues
$\lambda^0_\nu$ are simple.
\end{lemma}
\begin{proof}
  Begin by observing that $R = r(x,D)$---as an operator acting on
  $L^2(\R^2,\C^m)$---is invertible modulo a smoothing additive term.
  This means that one can construct a parametrix $L$ so that $LR = I$
  and $RL = I$ with both equations holding modulo a smoothing
  operator. To see why this is true, note that the matrix $r(x,\xi)$
  is a lower-order perturbation from the unitary matrix $r^0(x,\xi)$
  of eigenvectors of the principal symbol $a^0(x,\xi)$. The inverse of
  $r^0(x,\xi)$ is explicitly given by $\ell^0(x,\xi) = r^0(x,\xi)^*$.
  The symbol of $L$ can now be built as an expansion $\ell^0 + \ell^1
  + \ldots$, where each $\ell^n(x,\xi)$ is homogeneous of degree $-n$
  in $\xi$ and chosen to suppress the $O(|\xi|^{-j})$ contribution in
  $RL - I$ as well as $LR - I$. This construction implies that $L$ is
  pseudo-differential of order zero (componentwise). All of this is
  routine and detailed in \cite{Fol}[page 117].  An open question is
  whether $L$ can simply be chosen as the adjoint of $R$.
  
  In the sequel, $S$, $S_1$ and $S_2$ will denote a generic smoothing
  operator whose value may change from line to line. The composition
  of a pseudo-differential operator and a smoothing operator is
  obviously still smoothing. Set $f = Lu$ and let $A = a(x,D)$, so
  that $\d_t u = -iAu$. On the one hand, $u = Rf - Su$ and
\begin{equation}\label{eq:decoupling_dudt}
\d_t u = R \d_t f - S \d_t u  
= R \d_t f - S A u.
\end{equation}
On the other hand, Lemma \ref{teo:twisted_egv} gives 
\begin{equation}\label{eq:decoupling_Au}
Au = A R f - A S u = R \Lambda f + S_1 f + S_2 u 
= R \Lambda f + S u
\end{equation}
Comparing (\ref{eq:decoupling_dudt}) and (\ref{eq:decoupling_Au}), and
applying $L$ gives
\begin{equation}\label{eq:decoupling_dfdt}
\d_t f = -i \Lambda f + S u.
\end{equation}
This can be solved by Duhamel's formula,
\begin{equation}\label{eq:decoupling_f}
f(t) = e^{-it \Lambda} f(0) + 
\int_0^t e^{-i(t-\tau) \Lambda} Su(\tau) \, d\tau.
\end{equation}
We now argue that the integral term is, indeed, a smoothing operator
applied to the initial value $u_0$.
\begin{itemize}
\item First, the evolution operator $E(t) = e^{-itA}$ has a kernel $K(t,x,y)$ supported inside a neighborhood of the diagonal $y = x$ and for each $s \ge 0$, is well-known to map $H^s(\R^2,\C^m)$
  boundedly onto itself \cite{Lax}. Therefore, $S E(\tau)$ maps
  $H^{-s}$ to $H^s$ boundedly for every $s > 0$ and has a
  well-localized kernel in the sense of (\ref{eq:localkernel}). This
   implies that $S E(\tau)$ is a smoothing operator.
  
 \item Second, section \ref{sec:FIOhyp} shows that $e^{-it \Lambda}$
   is, for small $t$, a FIO of type $(1,0)$ and order zero, modulo a
   smoothing remainder. The composition of a FIO and a smoothing
   operator is smoothing. For larger $t$, think about $e^{-it
     \Lambda}$ as the product $(e^{-i\frac{t}{n} \Lambda})^n$ for
   appropriately large $n$.
  
\item And third, the integral extends over a finite interval $[0,t]$
  and may be thought as an average of smoothing operator---hence
  smoothing.
\end{itemize}
In short, $f(t) = e^{-it \Lambda} f(0) + S u_0$. Applying $R$ on both
sides of (\ref{eq:decoupling_dfdt}) finally gives
\[
u = R e^{-it \Lambda} L u_0 + S_1 u_0 + S_2 u = (R e^{-it \Lambda} L +
S) u_0
\]
which is what we set to establish.

It remains to see that the evolution operator $e^{-it \Lambda}$ for
the polarized components has the same block-diagonal structure as
$\Lambda$ itself. This is gleaned from equation
(\ref{eq:decoupling_dfdt}): evolution equations for two components
$f_{\nu_1}$, $f_{\nu_2}$ (corresponding to distinct eigenvalues) are
completely decoupled.
\end{proof}

\subsection{The Fourier Integral Operator parametrix}\label{sec:FIOhyp}

Lemma \ref{teo:decoupling} explained how to turn the evolution
operator $E(t)$ into the block-diagonal representation $e^{-it
  \Lambda}$. In this section, we describe how each of these blocks can
be approximated by a Fourier Integral Operator. The ideas here are
standard and our exposition is essentially taken from \cite{Gar} and
\cite{Sog}. The original construction is due to Lax \cite{Lax}.

Let us first assume that all eigenvalues of the principal symbol
$a^0(x,\xi)$ are simple. This is the situation where the matrix of
operators $\Lambda$ (Lemma \ref{teo:twisted_egv}) is diagonal with
diagonal elements $\Lambda_\nu$. Put $E_\nu(t) = e^{-it \Lambda_\nu}$,
the (scalar) evolution operator relative to the $\nu$th polarized
mode. We seek a parametrix $F_\nu(t)$ such that $S_\nu(t) = E_\nu(t) -
F_\nu(t)$ is smoothing of order $-\infty$.

Formally,
\[
f(t,x) = \int E_\nu(t) (e^{i x \cdot \xi}) \; \widehat{f_0}(\xi) \; d\xi. 
\]
Our objective is to build a high-frequency asymptotic expansion for 
$E_\nu(t)(e^{i x \cdot \xi})$ of the form 
\begin{equation}\label{eq:expansion}
e^{i \Phi_\nu(t,x,\xi)} \sigma_\nu(t,x,\xi),
\end{equation}
where $\sigma_\nu \sim \sigma^0_\nu + \sigma^1_\nu + \ldots$ and
$\sigma^n_\nu$ is homogeneous of degree $-n$ in $\xi$. As usual, one
obtains converging expansions through the use of appropriate
low-frequency cutoffs as in the proof of Lemma \ref{teo:twisted_egv}.

As is classical in asymptotic analysis, we proceed by applying $M_\nu
= \d_t + i \Lambda_\nu$ to the expansion (\ref{eq:expansion}) and
successively equate all the coefficients of the negative powers of
$|\xi|$ to zero, hence mimicking the relation $M_\nu E_\nu(t) (e^{i x
  \cdot \xi}) = 0$ which holds by definition. For obvious reasons, we
also impose that (\ref{eq:expansion}) evaluated at $t = 0$ be $e^{ix
  \cdot \xi}$.

We assume that $\Phi_\nu$ together with its partial derivatives in $t$
and $x$ be homogeneous of degree one in $\xi$ (note that
$\Phi_\nu(0,x,\xi) = x \cdot \xi$). Next, $\Lambda_\nu$ is also
expressed as a polyhomogeneous expansion $\Lambda_\nu \sim \sum_{j
  \geq 0} \lambda^{j}_\nu(x,D)$, where each symbol $\lambda_\nu^j(x,\xi)$
is homogeneous of degree $-j+1$ in $\xi$, compare with the proof of
Lemma \ref{teo:twisted_egv}. In passing, note that
$\lambda^0_\nu(x,\xi)$ is the $\nu$th eigenvalue of the principal
symbol $a^0(x,\xi)$, hence our consistent use of notations.

On the one hand,
\begin{equation}\label{eq:expansion_ddt}
e^{-i \Phi_\nu} \frac{\d}{\d t} [e^{i \Phi_\nu} \sigma_\nu] 
= i \frac{\d \Phi_\nu}{\d t} \sigma_\nu + \frac{\d \sigma_\nu}{\d t}.
\end{equation}
On the other hand, the action of $\Lambda_\nu$ on the oscillatory
product $e^{i \Phi_\nu} \sigma_\nu$ is more complicated but was
studied in \cite{Gar, Sog}. The following formula is available:
\begin{equation}\label{eq:expansion_lambdanu}
e^{-i \Phi_\nu} \Lambda_\nu [e^{i \Phi_\nu} \sigma_\nu] 
\sim \sum_\alpha \frac{1}{\alpha !} \d_\xi^\alpha 
\lambda_\nu(x, \nabla_x \Phi_\nu(t,x,\xi)) 
D_y^\alpha[e^{i q(t,x,y,\xi)} \sigma_\nu(t,y,\xi)] \vert_{y=x},
\end{equation}
where
\[
q(t,x,y,\xi) = \Phi_\nu(t,y,\xi) - \Phi_\nu(t,x,\xi) - 
\nabla_x \Phi_\nu(t,x,\xi) \cdot (y-x).
\]
Both equations (\ref{eq:expansion_ddt}) and
(\ref{eq:expansion_lambdanu}) are polyhomogeneous expansions. It is
rather tedious to show that this property indeed holds for
(\ref{eq:expansion_lambdanu}) and the discussion can be found in
\cite{Sog}.

Consider now the leading order relation (coefficient of $|\xi|$),
\begin{equation*}
i [\frac{\d \Phi_\nu}{\d t} + 
\lambda^{0}_\nu(x,\nabla_x \Phi_\nu)] \sigma_\nu^{0} = 0.
\end{equation*}
For nonzero $\sigma_\nu^{0}$, the term within brackets must vanish. For each
$\nu$, we are thus in the presence of a Hamilton-Jacobi equation
\begin{equation}\label{eq:HJ2}
\frac{\d \Phi_\nu}{\d t} + \lambda^0_{\nu}(x,\nabla_x \Phi_\nu) = 0,
\end{equation}
which are the well-known equations for the characteristic surfaces of
the original problem (\ref{eq:hyp}). As observed earlier, the initial
value is chosen as $\Phi_\nu(0,x,\xi) = x \cdot \xi$.

Local existence and uniqueness for (\ref{eq:HJ2}) with the given
initial condition is ensured as soon as $a^0(x,\xi)$ is $C^2$, from
the knowledge that $\Phi_\nu$ is constant along the bicharacteristics. See also \cite{Sog}.
The $\lambda^0_{\nu}$'s inherit the boundedness and smoothness
properties of $a^0(x,\xi)$. The solution to (\ref{eq:HJ2}) is not
expected to be global in time, because $\Phi_\nu$ would become
multi-valued when rays originating from the same point $x_0$ cross again later. This typically happens at
cusp points, when caustics start developing. We refer the reader to
\cite{Gar, Whi}.

Because the equation (\ref{eq:HJ2}) is homogeneous of degree one in
$\Phi_\nu$, the degree of homogeneity of $\Phi_\nu$ in $\xi$ must
match that of the initial condition at all times. This validates our
assumption that the phase is of degree one.

Consider the relation at the next order (coefficient of $|\xi|^0$),
\[
i [\frac{\d \Phi_\nu}{\d t}  + \lambda^0_\nu(x,\nabla_x
\Phi_\nu)] \sigma_\nu^{1} + \frac{\d \sigma_\nu^0}{\d t} + \nabla_\xi \lambda_\nu^0(x,\nabla_x \Phi_\nu) \cdot \nabla_x \sigma_\nu^0  = - i \lambda^1_\nu(x,\nabla_x \Phi_\nu) \sigma_\nu^0.
\]
The first term vanishes from (\ref{eq:HJ2}). What remains is a
transport equation for $\sigma_\nu^0$ along the bicharacteristic
vector field $\d_t + \nabla_\xi \lambda_\nu^0(x,\nabla_x \Phi_\nu)
\cdot \nabla_x$. This determines $\sigma_\nu^0$ uniquely from its
initial value $\sigma_\nu^0(0,x,\xi) = 1$. It is a smooth function of
$x$ and $\xi$, homogeneous of degree zero in $\xi$.

The other $\sigma_\nu^n$ are determined successively in a very similar
fashion. Each relation corresponding to a negative power $-n$ of
$|\xi|$ is a transport equation along the same vector field:
\begin{equation}\label{eq:expansion_general}
\frac{\d \sigma_\nu^n}{\d t} + 
\nabla_\xi \lambda_\nu^0(x,\nabla_x \Phi_\nu) \cdot \nabla_x \sigma_\nu^n  
= P_n(\sigma_\nu^0, \ldots, \sigma_\nu^n),
\end{equation}
where $P_n$ is a known differential operator applied to $\sigma_\nu^0,
\ldots, \sigma_\nu^n$. Because of (\ref{eq:HJ2}), the unknown
$\sigma_\nu^{n+1}$ will always disappear from the relation at the
degree $-n$. The initial condition is $\sigma_\nu^n(0,x,\xi) = 0$ for
$n > 0$. Since by construction, the right-hand side of
(\ref{eq:expansion_general}) is homogeneous of degree $-n$ in $\xi$,
the same holds for $\sigma_\nu^n$. The smoothness of each
$\sigma_\nu^n$ in $x$ and $\xi$ is inherited from that of the
principal symbol $a^0(x,\xi)$.

In the case where some eigenvalue $\lambda^0_\nu$ has multiplicity $p
> 1$, the construction of a FIO parametrix roughly goes the same way.
Let us use the same notation $\Lambda_{\nu}$ to denote the $p$-by-$p$
block corresponding to $\lambda^0_\nu$ in the matrix $\Lambda$ of
Lemma \ref{teo:twisted_egv}. The corresponding one-way evolution
operator is still denoted by $E_\nu(t) = e^{-it \Lambda_\nu}$. We seek
a $p$-by-$p$ matrix of operators $F_\nu(t)$ approximating $E_\nu(t)$.
The correct asymptotic expansion for $E_\nu(t)(e^{i x \cdot \xi})$ is
now
\[
e^{i \Phi_\nu(t,x,\xi)} \sum_{n = 0}^{\infty}
\bm{\sigma}_\nu^{n}(t,x,\xi),
\]
where $\Phi_\nu$ is a \emph{single} scalar phase function and each
$\bm{\sigma}_\nu^n$ is a $p$-by-$p$ matrix, homogeneous of degree $-n$
in $\xi$. Again, we apply $M_\nu = I \d_t + i \Lambda_\nu$ to this
expansion and successively equate identical powers of $|\xi|$ to zero.
The leading-order relation is satisfied if and only if $\Phi_\nu$
solves (\ref{eq:HJ2}) with initial condition $\Phi_\nu(0,x,\xi)
= x \cdot \xi$. The next equation, corresponding to $|\xi|^0$, is a
system of $p$ scalar evolutions equations which determines each column
of $\bm{\sigma}_\nu^0$ uniquely from the initial condition
$\bm{\sigma}_\nu^0(0,x,\xi) = I$. This situation repeats itself at all
orders and completely determines the expansion.

In all cases the parametrix $F_{\nu}$ thereby defined takes the form
of a FIO
\begin{equation}\label{eq:parametrix_FIO}
F_{\nu}(t) f_0(x) = \int e^{i \Phi_{\nu}(t,x,\xi)} 
\sigma_{\nu}(t,x,\xi) \hat{f}_0(\xi) \, d\xi,
\end{equation}
with $\Phi_{\nu}$ the phase function, homogeneous of degree one in
$\xi$, and $\sigma_{\nu}$ the amplitude, polyhomogeneous of degree 0
in $\xi$. The amplitude is scalar-valued or matrix-valued depending on
the multiplicity of $\lambda^0_\nu$. A suitable choice of cutoff
functions in $\xi$ will remove the singular behavior of each
$|\xi|^{-n}$ at the origin and guarantees convergence of the series
(\ref{eq:expansion}). This technicality, however, is the reason why
$F_\nu(t)$ is not exactly equal to $E_\nu(t)$ but is only an
approximation.

We emphasized that $\Phi_\nu$ may be defined only for small times. Put
$2t^*$ the infimum time for which a solution to (\ref{eq:HJ2}) ceases
to exist, uniformly in $\nu$. In order to establish results for large
times, we will simply compose evolution operators; e.g. $E_\nu(nt^*) =
E_\nu(t^*) \ldots E_\nu(t^*)$.
\begin{lemma}\label{teo:remainder_FIO}
  Define $t^*$ as above. In the setting of lemma \ref{teo:twisted_egv}
  denote by $\Lambda_\nu$ a block of $\Lambda$ and $E_\nu(t) = e^{-it
    \Lambda_\nu}$. Then for every $0 < t \leq t^*$, there exists a
  parametrix $F_\nu(t)$ for the evolution problem $\d_t f + i
  \Lambda_\nu f = 0$, that takes the form of a Fourier Integral
  Operator. For every such $t$ the phase function $\Phi_\nu$ is
  positive-homogeneous of degree one in $\xi$ and smooth in $x$ and
  $\xi$; the amplitude $\sigma_\nu$ is a symbol of type $(1,0)$ and
  order zero. The remainder $S_\nu(t) = E_\nu(t) - F_\nu(t)$ is a
  smoothing operator of order $-\infty$.
\end{lemma}
\begin{proof}
  The proof is a slight variation on the argument presented in
  \cite{Sog}[Pages 120 and below] and is omitted.
\end{proof}


\subsection{Sparsity of smoothing terms}
\label{sec:smooth}

The specialist will immediately recognize that a smoothing operator of
order $-\infty$ is very sparse in a curvelet frame. This is the
content of the following lemma.
\begin{lemma}
\label{teo:smoothing}
The curvelet entries of a smoothing operator $S$ obey the following
estimate: for each $N > 0$, there is a constant $C_N$ such that
\begin{equation}
\label{eq:smooth-estimate}
|\< \vf_{\mu}, S \vf_{\mu'} \>| \leq C_N \cdot 2^{-|j + j'|N}(1
+ |x_\mu - x_{\mu'}|)^{-N}. 
\end{equation}
\end{lemma}
Note that (\ref{eq:smooth-estimate}) is a stronger estimate than that
of Theorem \ref{teo:main}. Indeed,  our lemma implies that 
\[
|\< \vf_{\mu}, S \vf_{\mu'} \>| \leq 
C_{N} \cdot \omega(\mu, h_{\nu}(t,\mu'))^{-N}
\]
which is valid for each $N > 0$ and regardless of the value of $\nu$. 

\begin{proof}
 We know that $S$ maps $H^{-s}$ to
$H^s$ for arbitrary large $s$, so does its adjoint $S^{*}$. As a
result,
\begin{align*}
  | \< \vf_{\mu}, S \vf_{\mu'} \> | &
\leq |\< S^*\vf_{\mu}, \vf_{\mu'} \> |^{1/2} | 
\< \vf_{\mu}, S \vf_{\mu'} \> |^{1/2} \\
  &\leq \| S^*\vf_{\mu} \|_{H^s}^{1/2} 
\| \vf_{\mu'} \|_{H^{-s}}^{1/2}  \| \vf_{\mu} \|_{H^{-s}}^{1/2} 
\| S \vf_{\mu'} \|_{H^s}^{1/2} \\
  &\leq C \cdot \| \vf_{\mu'} \|_{H^{-s}} \| \vf_{\mu}
  \|_{H^{-s}} \leq C \cdot 2^{-(j+j')s}.
\end{align*}
On the other hand remember that curvelets have an essential spatial
support of size at most $O(1) \times O(1)$, which is the case when $j$
is small (wavelet regime $j = 0$). The action of $S$ is local on that
range of distances, so that
\[
| \< \vf_{\mu}, S \vf_{\mu'} \> | \leq C_N \cdot (1 + |x_{\mu} -
x_{\mu'}|)^{-N}.
\]
for arbitrary large $N > 0$. These two bounds can be combined to
conclude that the matrix elements of $S$ are negligible in the sense
defined above.
\end{proof}

\subsection{Proof of Theorem \ref{teo:main}}\label{sec:proof_main}

\begin{itemize}
\item Let us first show how the first assertion on the
  near-exponential decay of the curvelet matrix elements follows
  immediately from the second one, equation (\ref{eq:main2}). Let $a$
  be either a row or a column of the curvelet matrix and let
  $|a|_{(n)}$ be the $n$-largest entry of the sequence $|a|$. We have
$$
n \cdot |a|_{(n)}^{1/M} \le \|a\|_{\ell_{1/M}}
$$
and, therefore, it is sufficient to prove that the matrix $E$ has
rows and columns bounded in $\ell_p$ for every $p>0$. Consider the
columns. We need to establish
\begin{equation}\label{eq:ell_p_estimate}
\sup_{\mu',\nu'} \sum_{\mu, \nu} |E(t;\mu,\nu;\mu',\nu')|^p \leq C_{t,p},
\end{equation}
for some constant $C_{t,p} > 0$ growing at most like $C_{p}e^{K_{p}t}$
for some $C_p, K_p > 0$.

The sum over $\nu$ and the sup over $\nu'$ do not come
in the way since these subscripts take on a finite number of values.
The fine decoupling between the $m$ one-way components, crucial for
equation (\ref{eq:main2}), does not play any role here.

Let us now show that there exists $N$, possibly very large, so that
\[
\sum_\mu \omega(\mu,\mu')^{-Np} \leq C_{N,p},
\]
uniformly in $\mu'$. For the sum in $k$ and $\ell$ we can use the
bound (\ref{eq:sumkl}) with $Np$ in place of $N$, provided $Np$ obeys
$Np \geq 2$. What remains is
\[
C_{N,p} \cdot \sum_{j \geq 0} 2^{-|j - j'|Np} \cdot 2^{2|j - j'|},
\]
which is bounded by a constant depending on $N$ and $p$ provided again
that $Np \geq 2$.

Hence we proved the property for the columns. The same holds for the
rows because the same conclusion is true for the adjoint $E(t)^*$;
indeed, the adjoint solves the backward initial-value problem for the
adjoint equation $u_t = A^* u$, and $A^*$ satisfies the same
hyperbolicity conditions as $A$. We can therefore interchange the role
of the two curvelets and obtain
\[
\sup_{\mu',\nu'} \sum_{\mu,\nu} |E(t;\mu',\nu';\mu,\nu)|^p \le C_{t,p}. 
\]
Note that the classical interpolation inequality shows that $E(t)$ is
a bounded operator from $\ell_p$ to $\ell_p$ for every $0 < p \leq
\infty$.

\item We now turn to (\ref{eq:main2}). Let us assume first that all
  eigenvalues $\lambda^0_\nu$ of the principal symbol $a^0$ are
  simple.  According to Lemmas \ref{teo:decoupling} and
  \ref{teo:remainder_FIO}, each matrix element $E(t;\nu;\nu') =
  \mathbf{e_\nu} \cdot E(t)\mathbf{e}_{\nu'}$ of $E(t)$ can for fixed
  (possibly large) time $t>0$ be written as
\begin{equation}
  \label{eq:fio+smooth}
E(t;\nu;\nu') = \sum_{\nu''=1}^m R_{\nu, \nu''} (e^{-i\frac{t}{n}
  \Lambda_{\nu''}})^n L_{\nu'',\nu'} + S_{\nu,\nu'}(t).
\end{equation}
We have taken $n$ large enough---proportional to $t$---so that
$e^{-i\frac{t}{n} \Lambda_\nu}$ is a Fourier Integral Operator (mod
smoothing) for every $\nu$. Each $R_{\nu,\nu'}$ and $L_{\nu,\nu'}$ is
pseudo-differential of order zero and $S_{\nu,\nu'}(t)$ is smoothing.

Thanks to Lemma \ref{teo:smoothing}, it is sufficient to establish the
theorem about the the first term of (\ref{eq:fio+smooth}). To do this,
we invoke Theorem \ref{teo:main_FIO} from section \ref{sec:FIO} about
the sparsity of FIO's in a curvelet tight-frame.

It is an interesting exercise to notice that the ray dynamics is
equivalently expressed in terms of Hamiltonian flows,
\begin{equation}\label{eq:appendix_HS}
\left\{ \begin{array}{lll}
\dot x(t) & = \nabla_\xi \lambda^0(x(t),\xi(t)), 
\qquad & x(0) = x_0,\\
\dot \xi(t) & = - \nabla_x \lambda^0(x(t),\xi(t)), 
\qquad &\xi(0)  = \xi_0.
           \end{array}
      \right.
\end{equation}
 or canonical
transformations generated by the phase functions $\Phi_\nu$,
\begin{equation}\label{eq:canonical_phi}
\left\{ \begin{array}{ll}
x_0 &= \nabla_{\xi} \Phi(t,x(t),\xi_0), \\
\xi(t) &= \nabla_x \Phi(t,x(t),\xi_0),
           \end{array}
      \right.
\end{equation}
provided $\Phi(t,x,\xi)$ satisfies the Hamilton-Jacobi equation
$\frac{\d \Phi}{\d t} + \lambda^0(x,\nabla_x \Phi) = 0$ with the
initial condition $\Phi(0,x,\xi_0) = x \cdot \xi_0$. We obviously need
this property to ensure that the geometry of FIO's is the same as that
of hyperbolic equations.

Pseudo-differential operators are a special instance of Fourier
integral operators so the theorem equally applies to them. For
$E(t;\mu,\nu;\mu',\nu') = \< \vf_\mu, E(t;\nu;\nu') \vf_\mu' \>$
we get
\begin{align*}
|E(t;\mu,\nu;\mu',\nu')| \leq 
C_N \sum_{\nu'' = 1}^m \sum_{\mu_0} \cdots \sum_{\mu_{n}} \, 
&\omega(\mu,\mu_0)^{-N} \omega(\mu_0,{\mu_1}_{\nu''}(\frac{t}{n}))^{-N} 
\cdots \\
&\omega(\mu_{n-1},{\mu_n}_{\nu''}(\frac{t}{n}))^{-N}
\omega(\mu_{n},\mu')^{-N},
\end{align*}
for all $N > 0$. Inequality (\ref{eq:main_result}) then follows from
repeated applications of properties 3 and 4 of the distance $\omega$,
see proposition \ref{teo:omega}. The power growth in $t$ of the
overall multiplicative constant comes from the number of intermediate
sums over $\mu_0, \ldots, \mu_{n}$. There are $n+1 \sim t$ such sums
and they each introduce the same multiplicative constant $C_N$.

The reasoning is the same when at least some eigenvalues
$\lambda^0_\nu$ are degenerate. The subscript $\nu''$ now denotes the
flows i.e., the eigenvalues $\lambda^0_{\nu''}$ \emph{not} counting
their multiplicity. Each $R_{\nu,\nu''}$ is a row vector,
$e^{-i\frac{t}{n} \Lambda_{\nu''}}$ a matrix and $L_{\nu'',\nu'}$ a
column vector. The FIO parametrix for $e^{-i\frac{t}{n}
  \Lambda_{\nu''}}$ was constructed in such a way that only one flow
$h_{\nu''}$ appears in the majoration of its curvelet elements
(componentwise).  There is no intermediate sum over $\nu_0, \ldots
,\nu_n$ and that's the whole point of decoupling the polarized
components \emph{before} constructing the FIO parametrix.

\end{itemize}

\subsection{Relation to hyper-curvelets}

In section \ref{sec:curvgeomoptics} we introduced hyper-curvelets as
`polarized' curvelets which would not split into $m$ molecules along
the $m$ different flows. In light of section \ref{sec:decoupling}, it
is interesting to reformulate our main result (\ref{eq:main_result})
in terms of hyper-curvelets. We recall that
\[
\bm{\vf}^{(0)}_{\mu\nu} = r^0_\nu(x,D) \vf_\mu(x) = \frac{1}{(2
  \pi)^2} \int e^{ix\cdot \xi} r^0_\nu(x,\xi) \hat{\vf}_\mu(\xi)
\, d\xi.
\]

\begin{corollary}
  Define $E^{(0)}(t;\mu,\nu; \mu',\nu') = \< \bm{\vf}^{(0)}_{\mu\nu},
  E(t) \bm{\vf}^{(0)}_{\mu'\nu'} \>$. Then under the same assumptions
  as that of theorem \ref{teo:main} we have for all $N>0$
\beq\label{eq:main_hyper}
|E^{(0)}(t;\mu,\nu; \mu',\nu')| \leq C_{tN} \cdot
[\, \omega(\mu,\mu'_{\nu'}(t))^{-N} + 2^{-j'} \sum_{\nu'' \ne \nu'}
\omega(\mu,\mu'_{\nu''}(t))^{-N} \,].
\eeq
\end{corollary}

The main contribution to the right-hand side of this inequality is the
due to the $\nu'$th flow. All other flows are weighted by the small
factor $2^{-j'} = O(|\xi|^{-1})$ on the support of $\hat{\vf}_{\mu'}$.
In other words, the leaking of energy from one component $\nu$ to all
others when following hyper-curvelets $\bm{\vf}^{(0)}_{\mu\nu}$ is
smoothing of order $-1$, hence small at small scales.

\begin{proof}
  Equation (\ref{eq:main_hyper}) follows from Theorem \ref{teo:main}
  and the fact that the adjoint of the matrix operator $R^0$ whose
  columns are the $R^0_\nu = r^0_\nu(x,D)$ is an approximate left
  inverse for $R^0$---up to an error smoothing of order $-1$. Indeed,
  by the standard rules for composition and computation of the adjoint
  of pseudo-differential operators,
\begin{align*}
(R^0_\nu)^* R^0_{\nu'} &= ((r^0_\nu)^* \, \sharp \, r^0_{\nu'})(x,D) \\
&= ((r^0_\nu)^*r^0_{\nu'})(x,D) + \mbox{ order} - \! 1 \\
&= \delta_{\nu\nu'} I + \mbox{ order} - \! 1.
\end{align*}
We have used the fact that the dispersion matrix $a^0(x,\xi)$ is
assumed to be symmetric, hence admits an orthobasis of eigenvectors
$r^0_\nu(x,\xi)$. We then conclude from Theorem \ref{teo:main_FIO} applied
to pseudo-differential operators of order $-1$.
\end{proof}

Alternatively, we could have defined hyper-curvelets as
\[
\bm{\vf}^{(\infty)}_{\mu\nu} = r_\nu(x,D) \vf_\mu(x) = \frac{1}{(2
  \pi)^2} \int e^{ix\cdot \xi} r_\nu(x,\xi) \hat{\vf}_\mu(\xi)
\, d\xi.
\]
This would have given the same result.\footnote{We can only conjecture
  that the decoupling should be better if we use the improved
  $\bm{\vf}^{(\infty)}_{\mu\nu}$.} The reason why we did not use
hyper-curvelets in the preceding sections is that they do not
necessarily constitute a suitable practical basis to decompose
wavefields onto. We do not even know if they always constitute a
frame. Digital implementation would also seem less obvious.

\section{Representation of FIOs}
\label{sec:FIO}
\setcounter{theorem}{0}
\setcounter{equation}{0}

The purpose of this section is to show that Fourier Integral Operators
admit a sparse and well-organized structure in a curvelet frame. The
main result, Theorem \ref{teo:main_FIO}, is a key step in completing
the discussion of the previous section. (Observe that by construction,
the FIO's encountered in the previous section satisfy all the
assumptions stated in section \ref{sec:strategy} right below
(\ref{eq:FIO}).) As in the previous section, we will restrict the
discussion to $x \in \R^2$ which is no loss of generality, see section
\ref{sec:discussion}.

\subsection{Main results}

In the introduction section, we detailed a notion of Hamiltonian
correspondence for hyperbolic equations. This correspondence also
exists for FIO's and is `encoded' in the phase function $\Phi$ of the
FIO. It is called the canonical transformation associated to $\Phi$,
and is defined as the mapping $(x,\xi) \to (y,\eta)$ of phase-space
\begin{equation}\label{eq:canonical_transformation}
x = \nabla_\xi \Phi(y,\xi), \qquad \qquad \eta = \nabla_x \Phi(y,\xi).
\end{equation}
As suggested in section \ref{sec:proof_main}, this formulation is
equivalent to that involving trajectories along the bicharacteristic
flow as in equation (\ref{eq:HS}), provided the phase function solves
an appropriate Hamilton-Jacobi equation. This canonical transformation
induces a mapping of curvelet subscripts, denoted by $\mu' = h(\mu)$.

The main result for this section reads as follows.
\begin{theorem}
\label{teo:main_FIO} 
Let $T$ be a Fourier Integral Operator of order $m$ acting on
functions of $\R^2$, with the assumptions stated above, and
$T(\mu;\mu')$ denote its matrix elements in the complex curvelet
tight frame. Then with $h$ the curvelet index mapping and $\omega$ the
distance defined in (\ref{eq:omega}), the elements $T(\mu;\mu')$
obey for each $N > 0$
\[
|T(\mu;\mu')| \leq C_N \cdot  2^{m j'}\omega(\mu,h(\mu'))^{-N},
\]
for some $C_N > 0$. Moreover, for every $0 < p \leq \infty$,
$(T(\mu,\mu'))$ is bounded from $\ell^p$ to $\ell^p$.
\end{theorem}

The interpretation of Theorem \ref{teo:main_FIO} is in strong analogy
with that of Theorem \ref{teo:main}. Namely, a FIO has the
property of transporting and warping a curvelet into another
curvelet-like molecule. (Again, the choice of using complex-valued
curvelets is not essential, as a real curvelet would be mapped onto
two molecules.)

The proof of Theorem \ref{teo:main_FIO} relies on the factorization of
$T$ on the space-frequency support of $\vf_\mu$ as a nice
pseudo-local operator $T_{1,\mu}$ followed by a smooth change of
variables, or warping $T_{2,\mu}$. This decomposition goes as follows.

Let $\vf_\mu$ be a fixed curvelet centered around the lattice point
$(x_\mu, \xi_\mu)$ in phase-space. The phase of our FIO can be
decomposed as
\begin{equation}
  \label{eq:phase-linear}
  \Phi(x, \xi) = \Phi_\xi(x, \xi_\mu) \cdot \xi + \delta(x, \xi),
\qquad \phi_\mu(x) = \Phi_\xi(x, \xi_\mu). 
\end{equation}
In effect, the above decomposition `linearizes' the frequency variable
and is classical, see \cite{SSS,Ste}. With these notations, we may
rewrite the action of $T$ on a curvelet $\vf_\mu$ as
\begin{equation}
  \label{eq:FIO-curvelet}
  (T\vf_\mu)(x) = \int e^{i \phi_\mu(x) \cdot \xi} e^{
  i \delta(x, \xi)} \sigma(x, \xi) \hat{\vf}_\mu(\xi) \, d\xi.
\end{equation}
Now for a fixed value of the parameter $\mu$, we introduce the
decomposition
\[
T  = T_{2,\mu} T_{1,\mu}, 
\]
where 
\begin{equation}
  \label{eq:T1}
  (T_{1,\mu} f)(x) = \int e^{i x \cdot \xi} b_\mu(x, \xi) 
\hat{f}(\xi) \, d\xi, \quad (T_{2,\mu} f)(x) = f (\phi_\mu(x)), 
\end{equation}
with $b_\mu(x,\xi) = e^{i \delta(\phi_\mu^{-1}(x), \xi)}
\sigma(\phi_\mu^{-1}(x), \xi))$. This decomposition allows the
separate study of the nonlinearities in frequency $\xi$ and space $x$
in the phase function $\Phi$. The point is that both $T_{1,\mu}$ and
$T_{2,\mu}$ are sparse in a curvelet tight frame---only for very
different reasons.

\begin{theorem} 
\label{teo:T1}
Let $(\vf_{\mu})_\mu$ be a tight frame of curvelets compactly
supported in frequency. For each $\mu$, $T_{1,\mu}$ maps $\vf_\mu$
into a curvelet molecule $m_\mu$ with arbitrary regularity $R$,
uniformly over $\mu$ in the sense that the constants in estimates
(\ref{eq:localization}) and (\ref{eq:bandpass}) do not depend on $\mu$.
\end{theorem}

As we shall see, the proof of Theorem \ref{teo:T1}, presented in
section \ref{sec:proof_T1}, relies on the property of compact support
in frequency of the $\vf_\mu$. In contrast the corresponding result
for the operators $T_{2,\mu}$ which we present next, is extraordinarily
simplified if one uses curvelets compactly supported in space.
Although well localized in space, the tight frame introduced in
section \ref{sec:curvelets} does not meet this requirement. In order
to circumvent this technical difficulty, we introduce compactly
supported \emph{curvelet atoms} in section \ref{sec:atoms}. They are
built on the model of atomic decompositions, standard in approximation
theory \cite{FJW}.

\begin{theorem}
\label{teo:T2} 
Let $(\rho_{\mu})_{\mu}$ be a family of complex-valued curvelet atoms,
compactly supported in space, with regularity $R$. Denote by $h$ the
canonical index correspondence associated to $\Phi$, as defined above.
For each $\mu$, $T_{2,\mu}$ maps $\rho_{\mu}$ into a molecule
$m_{h(\mu)}$ of the same regularity $R$, uniformly over $\mu$.
\end{theorem}

The latter theorem says that the `warped' atom $\rho_\mu \circ
\phi_\mu$ is still an atom, only its scale, orientation, and location
may have been changed. That a smooth warping preserves the sparsity of
curvelet expansions is a result of independent interest.

The remaining 3 sections are devoted to the proofs of Theorems
\ref{teo:T1}, \ref{teo:T2} and \ref{teo:main_FIO}. The dependence of
$\phi_\mu$ upon $\mu$ is not essential in proving Theorems
\ref{teo:T1}, \ref{teo:T2} as the only property of interest is that
the derivatives of $\phi_\mu$ are bounded from above and below
uniformly over $\mu$ (which follows from our assumptions about
$\Phi$). This is the reason why in the next sections we will drop the
explicit dependence on $\mu$ and work with a generic warping $\phi$.

\subsection{Proof of Theorem \ref{teo:T1}}\label{sec:proof_T1}
\setcounter{theorem}{0}
\setcounter{equation}{0}

We will assume without loss of generality that our curvelet
$\vf_\mu$ is centered near zero ($k = 0$) and is nearly vertical
($\theta_\mu = 0$).

Set $m_\mu = T_1 \vf_{\mu}$. We first show that $m_\mu$ obeys the
smoothness and spatial localization estimate of a molecule
(\ref{eq:localization}).  With the same notations as before, recall
that $m_\mu$ is given by
\begin{equation}
\label{eq:defbmu}
m_{\mu}(x) = \int e^{ i x \cdot \xi} 
b_{\mu}(x,\xi) \hat{\vf}_{\mu}(\xi) \, d\xi, \quad 
b_{\mu}(x,\xi) = e^{ i \delta(\phi^{-1}(x),\xi)}
\sigma(\phi^{-1}(x),\xi).
\end{equation}
To study the spatial decay of $m_{\mu}(x)$, we introduce the
differential operator
\[
L_{\xi} = I - 2^{2j} \frac{\partial^2}{\partial \xi_1^2} - 2^{j}
\frac{\partial^2}{\partial \xi_2^2},  
\]
and evaluate the integral (\ref{eq:defbmu}) using an integration by
parts argument. First, observe that 
\[
L_{\xi}^N e^{ i x \cdot \xi} = \left(1 +  |2^j x_1|^2 + 
  |2^{j/2} x_2|^2\right)^{N} \, e^{ i x \cdot \xi}.
\]
Second, we claim that for every integer $N \geq 0$, 
\begin{equation}\label{eq:bigOone}
|L_{\xi}^N [b_{\mu}(x,\xi) \hat{\vf}_{\mu}(\xi)]| 
\le C \cdot 2^{-3j/4}. 
\end{equation}
(The factor $2^{-3j/4}$ comes from the $L^2$ normalization of
$\hat{\vf}_{\mu}$.)  This inequality is proved in appendix
\ref{sec:appendix_FIO}.  Hence,
\[
m_{\mu}(x) = \left(1 +  |2^j x_1|^2 +  |2^{j/2}
x_2|^2\right)^{-N} \, \int L_{\xi}^N [b_{\mu}(x,\xi)
\hat{\vf}_{\mu}(\xi)] \, e^{ i x \cdot \xi}.
\]
Since $|L_{\xi}^N [b_{\mu}(x,\xi) \hat{\vf}_{\mu}(\xi)]| \le C
\cdot 2^{-3j/4}$ and is supported on a dyadic
rectangle $R_\mu$, of length about $2^{j}$ and width $2^{j/2}$, we then
established that
\[
|m_{\mu}(x)| \le C \cdot \frac{2^{3j/4}}{\left(1 + 
    |2^j x_1|^2 +  |2^{j/2} x_2|^2\right)^{N}}.
\]
The derivatives of $m_{\mu}$ are essentially treated in the
same way. Begin with
\begin{eqnarray*}
\partial_x^\alpha (e^{ i x \cdot \xi} 
b_{\mu}(x,\xi)) & = & \sum_{\beta + \vf \le \alpha} 
\partial^{\beta} (e^{ i x \cdot \xi}) \,  
\partial^{\vf} (b_{\mu}(x,\xi)) \\
& = & \sum_{\beta + \vf \le \alpha}    
\partial^{\vf} (b_{\mu}(x,\xi)) \, \xi^\beta e^{ i x \cdot \xi} 
\end{eqnarray*}
Therefore, the partial derivatives of $m_\mu$ are given by  
\begin{equation}
\label{eq:derbmu}
(\partial_x^\alpha m_{\mu})(x) = 
\sum_{\beta + \vf \le \alpha} I_{\beta, \vf}(x), 
\end{equation}
where 
\begin{equation}
  \label{eq:Ibeta}
  I_{\beta,\vf}(x) = \int e^{ i x \cdot \xi}
\partial_x^{\vf}(b_{\mu}(x,\xi)) \xi^{\beta} \hat{\vf}_{\mu}(\xi) \,
d\xi.
\end{equation}
First, observe that on the support of $\hat{\vf}_\mu$,
$|\xi|^\beta$ obeys $|\xi|^\beta \le C \cdot 2^{j\beta_1} \cdot
2^{j\beta_2/2}$. Second, the term $\partial_x^\vf b(x,\xi)$ is of
the same nature as $b_\mu(x,\xi)$ in the sense that it obeys all the
same estimates as before. In particular, we claim that for every
integer $N \geq 0$,
\begin{equation}\label{eq:bigOone2}
|L_{\xi}^N [\partial_x^\vf b_{\mu}(x,\xi) \xi^\beta 
\hat{\vf}_{\mu}(\xi)]| 
\le C \cdot 2^{-3j/4} \cdot 2^{j\beta_1} \cdot 2^{j\beta_2/2}. 
\end{equation}
Hence, the same argument as before gives  
\[
|I_{\beta,\vf}(x)| \le C \cdot \frac{2^{3j/4} \cdot 2^{j\beta_1}
  \cdot 2^{j\beta_2/2}}{\left(1 +  |2^j x_1|^2 + 
    |2^{j/2} x_2|^2\right)^{N}}.
\]
Now since $\beta \le \alpha$, we may conclude that 
\[
|(\partial^\alpha_x m_\mu)(x)| \le 
C \cdot \frac{2^{3j/4} \cdot 2^{j\alpha_1}
  \cdot 2^{j\alpha_2/2}}{\left(1 +  |2^j x_1|^2 + 
    |2^{j/2} x_2|^2\right)^{N}}. 
\]
This establishes the smoothness and localization property. 

The above analysis shows that $m_\mu$ is a ``ridge'' of effective
length $2^{-j/2}$ and width $2^{-j}$; to prove that $m_\mu$ is a
molecule, we now need to evidence its oscillatory behavior across the
ridge. In other words, we are interested in the size of the Fourier
transform at low frequencies (\ref{eq:bandpass})--(\ref{eq:bp2}).

Formally, the Fourier transform of $m_\mu$ is given by 
\begin{equation}
  \label{eq:fourier-T1}
   \hat{m}_{\mu}(\xi) = \int \! \!  \int 
  e^{ i x \cdot (\eta-\xi)} b_{\mu}(x,\eta)
  \hat{\vf}_{\mu} (\eta)\, dx d\eta.  
\end{equation}
We should point out that because the amplitude $b$ is not of compact
support in $x$, the sense in which (\ref{eq:fourier-T1}) holds is not
obvious. This is a well-known phenomenon in Fourier analysis and a
classical technique to circumvent such difficulties would be to
multiply $m_\mu$ (or equivalently $b_\mu$) by a smooth and compactly
supported cut-off function $\chi(\epsilon x)$ and let $\epsilon$ tend
to zero.  We omit those details as they are standard.

Set $D_1 = -i \frac{\d}{\d x_1}$. 
To develop bounds on $|\hat{m}_{\mu}(\xi)|$, observe that 
\[
D_1^N e^{ i x \cdot \eta} = (\eta_1)^N.  
\]
An integration by parts then gives 
\[
\hat{m}_{\mu}(\xi) = \int \! \!  \int e^{ i x \cdot \eta} D_1^N
\left(e^{-  i x \cdot \xi} \, b_\mu(x, \xi)\right) \,
\eta_1^{-N} \hat{\vf}_\mu(\eta) \, dx d\eta. 
\]
Hence, 
\[
\hat{m}_{\mu}(\xi) = \sum_{m = 0}^N c_m \, \xi_1^m \hat{F}_m(\xi),
\]
where 
\[
F_m(x) = \int e^{ i x \cdot \eta}  (\partial^{n-m}_{x_1}
b(x,\xi)) \, \eta_1^{-n} \hat{\vf}_\mu(\eta) \, d\eta. 
\]
Note that $F_m$ is exactly of the same form as (\ref{eq:Ibeta})---but
with $\eta_1^{-n}$ instead of $\eta^\beta$---and therefore, the
exact same argument as before gives
\[
|F_m (x)| \le C \cdot \frac{2^{3j/4} \cdot 2^{-jn}}
{\left(1 +  |2^j x_1|^2 +  |2^{j/2} x_2|^2\right)^{N}}.
\]
We then established
\[
\|\hat{F}_m\|_{L_\infty} \le \|F\|_{L_1} \le C_m \cdot 2^{-3j/4} \cdot 
2^{-jn}, 
\]
which gives
\[
|m_\mu(\xi)| \le C \cdot 2^{-3j/4} \cdot 2^{-jn} \cdot (1 + |\xi|^n), 
\]
as required. This finishes the proof of Theorem \ref{teo:T1}. 

The careful reader will object that we did not study the case of
coarse scale curvelets; it is obvious that coarse scale elements are
mapped into coarse scale molecules and, here, the argument would 
not require the deployment of the sophisticated tools we exposed
above. We omit the proof.

\subsection{Atomic decompositions}\label{sec:atoms}

As we will see later, to prove our main result and especially Theorem
\ref{teo:T2}, it would be most helpful to work with tight frames of
curvelet compactly supported in space.  Unfortunately, it is unclear
at this point how to construct such tight frames with nice frequency
localization properties. However, there exist useful atomic
decompositions with compactly supported curvelet-like atoms. We now
explore such decompositions.

In this section, the notation $f_{a,\theta}$ refers to the function
obtained from $f$ after applying a parabolic scaling and a rotation  
\[
f_{a, \theta}(x) = a^{-3/4} f\left( D_a R_\theta x \right), \quad 
D_a = \begin{pmatrix} 1/a & 0 \\ 0 &  1/\sqrt{a} \end{pmatrix}, 
\]
and where $R_\theta$ is the rotation matrix which maps the vector
$(1,0)$ into $(\cos \theta, - \sin \theta)$. 
Note that this is an isometry as
\[
\|f_{a, \theta}\|_{L_2} = \|f\|_{L_2}.
\]

In \cite{Sm1}, Smith proved the following result: let $\tilde{\psi}$
be a Schwartz function obeying $\hat{\tilde{\psi}}(1,0) \neq 0$; then
one can find another Schwartz function $\psi$, and a function $q(\xi)$
such that the following formula holds
\begin{equation}
  \label{eq:smith-calderon}
  q(\xi) \int_{a \le 1} \hat{\tilde{\psi}}_{a,\theta}(\xi) 
\hat{\psi}_{a,\theta}(\xi) \, a da d\theta = r(\xi);  
\end{equation}
here $r$ is a smooth cut-off function obeying 
\[
r(\xi) =  \left\{ \begin{array}{ll}
                          1 & \quad |\xi| \ge 2 \\
                     0 & \quad |\xi| \le 1
           \end{array}
      \right. , 
\]
and $q$ is a standard Fourier multiplier of order zero; that is, for
each multiindex $\alpha$, there exists a constant $C_\alpha$ such that
\[
|\d^\alpha_\xi q(\xi)| \le C_\alpha (1 + |\xi|)^{-|\alpha|}. 
\] 

This formula is useful because it allows us to express any object
whose Fourier transform vanishes on $\{|\xi| \le 2\}$ as a continuous
superposition of curvelet-like elements.  We now make some specific
choices for $\vf$. In the remainder of this section, we 
will take $\tilde{\psi}(x) = \psi(-x)$ and the function $\psi$ of the form 
\begin{equation}
\label{eq:psi-choice}
\psi(x_1, x_2) = \psi^D(x_1) \, \vf(x_2), 
\end{equation}
where both $\vf$ and $\psi^D$ are compactly supported and obey 
\[
\text{Supp}\, \vf \subset [0,1], 
\quad \text{Supp}\, \psi^D \subset [0,1]. 
\]
We will assume that $\vf$ and $\psi^D$ are $C^\infty$ and that the
function $\psi^D$ has vanishing moments up to order $D$, i.e.
\begin{equation}
\label{eq:psiD-vanish}
\int \psi^D(x_1) x_1^k \, dx_1 = 0, \quad k = 0, 1, \ldots, D. 
\end{equation}

For each $a \le 1$, each $b
\in \R^2$ and each $\theta \in [0,2\pi)$, introduce
\begin{equation}
\label{eq:continuous-curvelet}
       \psi_{a,\theta, b} (x) := \psi_{a,\theta}(x - b) = 
  a^{-3/4} \psi\left(D_a R_\theta(x - b)\right);  
\end{equation}
and given an object $f$, define coefficients by
\begin{equation}
\label{eq:coeff}
\cR(f)(a,b,\theta) = \int \overline{\psi_{a,\theta,b}(x)} f(x) dx.
\end{equation}
Now, suppose for instance that $\hat{f}$ vanishes over $|\xi| \le 2$,
then (\ref{eq:smith-calderon}) gives the exact reconstruction formula
\begin{equation}
\label{eq:CCT}
f(x) = \int_{a \ge 1} \cR(q(D) f))(a,b,\theta) 
{\psi}_{a,\theta,b}(x) \mu(da d\theta db), 
\end{equation}
with $\mu(da d\theta db) = a da d\theta db$. In the remainder of this
section, we will use the shorter notation $d\mu$ for $\mu(da d\theta
db)$. 

As is now well-established, the reproducing formula may be turned into
a so-called `atomic decomposition'. Not surprisingly, our atomic
decomposition will just mimic the discretization of the curvelet frame
as introduced in section \ref{sec:curvelets}. With the
notations of that section, we introduce the cells $Q_\mu$ defined as
follows: for $j \ge 0$, $\ell = 0, 1, \ldots 2^{\lfloor j/2 \rfloor} -
1$ and $k = (k_1, k_2) \in \bZ^2$, the cell $Q_\mu$ is the collections
of triples $(a, \theta, b)$ for which
\[
2^{-(j+1)} \le a < 2^{-j}, \quad 
|\theta - \theta_\mu| \le \frac{\pi}{2} 2^{-\lfloor j/2 \rfloor}
\] 
and 
\[
D_{2^{-j}} R_{\theta_\mu} b \in [k_1, k_1 + 1) \times [k_2, k_2 + 1). 
\]
Note that $\int_{Q_\mu} d\mu = 3\pi/2$ for $j$ even, and $3 \pi$ for
$j$ odd.  We may then break the integral (\ref{eq:CCT}) into a sum of
terms arising from different cells, namely,
\begin{equation}
\label{eq:atomic-decomposition}
f(x) = \sum_\mu \alpha_\mu \rho_\mu(x)
\end{equation}
where 
\begin{equation}
\label{eq:atoms}
\alpha_\mu = \|\cR(q(D) f)\|_{L_2(Q_\mu)}, \quad \rho_\mu(x) = 
\frac{1}{\alpha_\mu} \int_{Q_\mu}  \cR(q(D) f))(a,b,\theta) 
{\psi}_{a,\theta,b}(x)\,  d\mu. 
\end{equation}

Of course, the decomposition (\ref{eq:atomic-decomposition}) greatly
resembles the tight frame expansion, compare
(\ref{eq:tight-frame}). In particular, the atoms $\rho_\mu$ are
curvelet-like in the sense that they share all the properties of the
tight frame $( \vf_\mu)_\mu $ -- only they are compactly supported in
space. In the remainder of the paper, we will call these elements {\em
  atoms.} Below are some crucial properties of these atoms.
\begin{lemma}
\label{teo:atoms}
Rewrite the atoms $\rho_\mu$ as $\rho_\mu(x) = 2^{3j/4} a^{(\mu)}\left(D_{2^{-j}}
  R_{\theta_\mu} x - k \right)$. In other words, $\rho_\mu$ is obtained from
$a^{(\mu)}$ after parabolic scaling, rotation, and translation. For
all $\mu$, the $a^{(\mu)}$'s verify the following properties.
  \begin{itemize}
  \item Compact support;
  \begin{equation}
    \label{eq:compact-support}
    \text{Supp}\,\, a^{(\mu)} \subset c Q.
  \end{equation}
\item Nearly vanishing moment along the horizontal axis; let $m =
  D/2$. Then for each $k = 0, 1, \ldots, m$, there is a constant $C_m$
  such that
  \begin{equation}
    \label{eq:nearly-vanishing-moments}
    \int a^{(\mu)}(x_1,x_2) x_1^k \, dx_1 \le C_m \cdot 2^{-j(m+1)}.
  \end{equation}
\item Regularity; for every multiindex $\alpha$
\begin{equation}
  \label{eq:regularity}
  |\d^\alpha_x a^{(\mu)}(x)| \le C_\alpha. 
\end{equation}
\end{itemize}
In (\ref{eq:nearly-vanishing-moments}) and (\ref{eq:regularity}), the
constants may be chosen independently of $\mu$ and $f$.
\end{lemma}
\begin{proof}
See appendix \ref{sec:appendix_FIO} 
\end{proof}

Needless to say that curvelet atoms are molecules with spatial compact
support, compare lemma \ref{teo:atoms} with the definition of a
molecule. Finally, observe (and this is important) that it is of
course possible to decompose a molecule into a series of atoms
\[
m_\mu = \sum_{\mu'} \alpha_{\mu\mu'} \rho_{\mu'}. 
\]
The coefficients would then obey the same estimate as in lemma
\ref{teo:almost_orthogonality}
\begin{equation}\label{eq:molec-atom}
|\alpha_{\mu\mu'}| \leq C_N \cdot \omega(\mu, \mu')^{-N},
\end{equation}
and in particular, for each  $p > 0$,
\[
\sup_{\mu} \sum_{\mu'} |\alpha_{\mu\mu'}|^p < C_p. 
\]
This is briefly justified in appendix \ref{sec:appendix_FIO}.

\subsection{Proof of Theorem \ref{teo:T2}}\label{sec:proof_T2}

As mentioned earlier, curvelet atoms depend in a nonessential way
upon the object $f$ we wish to analyze and we shall drop this
dependence in our notations. To prove Theorem \ref{teo:T2}, recall
that we need to show that for each curvelet atom $\rho_\mu$ with
regularity $R$, the 'warped' atom $\rho_\mu \circ \phi$ is also a
curvelet atom, with the same regularity.

As in Section \ref{sec:atoms}, we suppose our curvelet atom is of the
form
\[
\rho_\mu(x) = 2^{3j/4} a^{(\mu)}(D_{2^{-j}} R_{\theta_\mu} (x - x_\mu)), 
\]
where $a^{(\mu)}$ obeys the conditions of Lemma \ref{teo:atoms}.
(Here, the location $x_\mu$ may be formally defined by $x_\mu = (D_{2^{-j}}
R_{\theta_\mu})^{-1}k_\delta$.)  Define $y_\mu$ and $A_\mu$ by
\begin{equation}
  \label{eq:def-microlocal}
  y_\mu = \phi^{-1}(x_\mu), \quad \text{and} \quad 
A_\mu = (\nabla \phi)(y_\mu)
\end{equation}
so that
\[
\phi(y) = x_\mu + A_\mu(y - y_\mu) + g(y - y_\mu). 
\]
With these notations, it is clear that the warped atom $\rho_\mu \circ
\phi$ will be centered near the point $y_\mu$; that is,
\[
\rho_\mu(\phi(y)) = 2^{3j/4} a^{(\mu)}\left(D_{2^{-j}} 
R_{\theta_\mu} (A_\mu (y - y_\mu) + g(y - y_\mu))\right). 
\]
To simplify matters, we first assume that $A_\mu$ is the identity and
show that $\rho_\mu \circ \phi$ is a curvelet atom with the same scale
and orientation as $\rho_\mu$. Later, we will see that in general,
$\rho_\mu \circ \phi$ is an atom whose orientation depends upon
$A_\mu$, and whose scale may be taken to be the same as that of
$\rho_\mu$.  Assume without loss of generality that $\theta_\mu = 0$ and
$y_\mu = 0$ (statements for arbitrary orientations and locations are
obtained in an obvious fashion) so that
\begin{equation}
  \label{eq:A-equal-I}
  \rho_\mu(\phi(y)) = 2^{3j/4} a^{(\mu)}\left(D_{2^{-j}} (y + g(y))\right) =
2^{3j/4} b^{(\mu)}(D_{2^{-j}} y),
\end{equation}
with
\[
b^{(\mu)}(y) =  a^{(\mu)}\left(y + D_{2^{-j}} g(D_{2^j}y)\right).
\]
The atom $a^{(\mu)}$ is supported over a square of sidelength about 1;
likewise, $b^{(\mu)}$ is also compactly supported in a box of roughly
the same size---uniformly over $\mu$.  We then need to derive
smoothness estimates and show that $b^{(\mu)}$ obeys
\begin{equation}
\label{eq:warped-smooth}
|\partial^\alpha b^{(\mu)}(y)| \le C_\alpha, \quad |\alpha| \le R. 
\end{equation}
Over the support of $\rho_\mu \circ \phi$, $g = (g_1,g_2)$ deviates
little from zero and for each $k = 1, 2$, $g_k$ obeys
\[
|g_k(y)| \le C \cdot 2^{-j}, \quad |\partial^\alpha g_k(y)| \le C
\cdot 2^{-j/2}, \,\, |\alpha| = 1.
\]
Similarly, for each $\alpha$, $|\alpha| > 1$,  
\begin{equation}
\label{eq:deriv-h}
|\partial^\alpha g_k(y)| \le C_\alpha. 
\end{equation}
These estimates hold uniformly over $\mu$.  It follows that for
$|y_1|, |y_2| \le C$ and each $\alpha$, the perturbation $g$ obeys
\begin{equation}
\label{eq:deriv-h2}
2^j \cdot |\partial^\alpha g_1(2^{-j} y_1, 2^{-j/2} y_2)| \le C_\alpha, \quad 
2^{j/2} \cdot |\partial^\alpha g_2(2^{-j} y_1, 2^{-j/2} y_2)| \le C_\alpha.
\end{equation}
The bound (\ref{eq:warped-smooth}) is then a simple consequence of
(\ref{eq:deriv-h2}) together with the fact that all the derivatives of
$a^{(\mu)}$ up to order $R$ are bounded, uniformly over $\mu$.

We now show that $\rho_{\mu} \circ \phi$ exhibits the appropriate
behavior at low frequencies.
\begin{align*}
  \widehat{\rho_{\mu} \circ \phi}(\xi) &= \int e^{-i x \cdot \xi} 
\rho_{\mu}(\phi(x)) \, dx \\
  &= \int e^{-i \phi^{-1}(x) \cdot \xi} \rho_{\mu}(x)
  \frac{dx}{|\det \nabla \phi|(\phi^{-1}(x))}.
\end{align*}
We will use the nearly vanishing moment property of $\rho_\mu$. Set
\[
S_\xi(x) = e^{-i \phi^{-1}(x) \cdot \xi}/|\det \nabla
\phi|(\phi^{-1}(x));  
\]
note that over the support of $\rho_\mu$ and for each $N \le R$, we
have available the following upper bound on the partial derivative of
$S_\xi$
\[
|\partial_1^N S_\xi(x)| \le C_N \cdot (1+|\xi|)^N.
\]
Classical arguments give  
\begin{equation}
\label{eq:final-estimate-1}
\widehat{\rho_{\mu} \circ \phi}(\xi) = \sum_{k = 0}^{n -1} 
\int \frac{\partial^k_1
S_\xi(0,x_2)}{k!} \, dx_2 \int \rho_{\mu}(x_1,x_2) x_1^k \, dx_1 dx_2 + E, 
\end{equation}
where $E$ is a remainder term obeying 
\begin{equation}
\label{eq:final-estimate}
|E| \le C_n \cdot 2^{-3j/4} \cdot 2^{-jn} \cdot 
\sup |\partial_1^n S_\xi(x)| \le  C_n \cdot 2^{-3j/4} 
\cdot 2^{-jn} (1 + |\xi|^n). 
\end{equation}
The near-vanishing moment property gives that each term in the
right-hand side of (\ref{eq:final-estimate-1}) obeys the estimate in
(\ref{eq:final-estimate}).  This proves that the Fourier transform of
$\rho_{\mu} \circ \phi$ obeys
\[
|\widehat{\rho_{\mu} \circ \phi}(\xi)| \le   C_n \cdot 2^{-jn} (1 + |\xi|^n)
\] 
as required.

We now discuss the case where the matrix $A_\mu$ is not the identity. 
In this case, (\ref{eq:A-equal-I}) becomes
\[
\rho_\mu(\phi(y)) = m_\mu(A_\mu y), 
\]
with 
\[
m_\mu(y) = 2^{3j/4} a^{(\mu)}\left(D_{2^j} (y + \tilde g(y))\right),
\quad \text{and} \quad \tilde g(y) = g(A_\mu^{-1} y).
\]
Our assumptions about FIOs guarantee that $|A_\mu^{-1}|$ is uniformly
bounded and, therefore, it follows from the previous analysis that
$m_\mu$ is a curvelet atom. As a consequence $\rho_\mu \circ \phi$ is a
curvelet atom with the same regularity $R$ since it is clear that
bounded linear transformations of the plane map curvelet atoms into
curvelet atoms.

\subsection{Proof of Theorem \ref{teo:main_FIO}}
\label{sec:together}
\setcounter{theorem}{0}
\setcounter{equation}{0}

Let $\vf_{\mu_0}$ be a fixed curvelet and decompose $T$ as
$T_{2,\mu_0} \circ T_{1,\mu_0}$.  First, Theorem \ref{teo:T1} proved
that $T_{1,\mu_0} \vf_{\mu_0}$ is a curvelet molecule $m_{\mu_0}$
which we will express as a superposition of curvelet atoms
$\rho_{\mu_1}$
\[
T_{1,\mu_0} \vf_{\mu_0} = m_{\mu_0} = \sum_{\mu_1} \beta_0(\mu_1,\mu_0) 
\rho_{\mu_1}. 
\]
Second, for each $\mu_1$, Theorem \ref{teo:T2} shows that $T_{2,\mu_1}
\rho_{\mu_1}$ is a molecule $m_{h(\mu_1)}$ at the location $h(\mu_1)$.
We are not exactly in that setting since in $T_{2,\mu_0}
\rho_{\mu_1}$, the subscripts do not in general match. This does not
pose any difficulty since Theorem \ref{teo:T2} can be understood as a
statement concerning general warpings $\phi$. We can define the map
$h_{\mu_0}$ as induced by the transformation $(x,\xi) \to (y,\eta)$
given by
\[
x = \nabla_x \Phi(y,\xi_{\mu_0}), \qquad \qquad \eta = \nabla_\xi
(y,\xi_{\mu_0})
\]
(compare this with equation (\ref{eq:canonical_transformation})).
Then, according to Theorem \ref{teo:T2}, \, $T_{2,\mu_0} \rho_{\mu_1}$
is a molecule $m_{h_{\mu_0}(\mu_1)}$ at the location
$h_{\mu_0}(\mu_1)$. So
\[
\<\vf_{\mu_2}, T_{2,\mu_0} \rho_{\mu_1}\> = \beta_1(\mu_2,h_{\mu_0}(\mu_1)).
\]
Hence, 
\[
\< \vf_{\mu_2}, T \vf_{\mu_0} \> = \sum_{\mu_1}
\beta_1(\mu_2,h_{\mu_0}(\mu_1)) \beta_0(\mu_1,\mu_0). 
\]
Of course, both $\beta_0$ and $\beta_1$ obey very special decay properties.
\begin{itemize}
\item By Theorem \ref{teo:T1} and Lemma \ref{teo:almost_orthogonality},
  $|\beta_0(\mu_1, \mu_0)| \leq C_n \cdot \omega(\mu_1,\mu_0)^{-N}$ for
  arbitrarily large $N > 0$, provided that the selected atoms are
  regular enough.
\item By Theorem \ref{teo:T2} and Lemma \ref{teo:almost_orthogonality},
  $|\beta_1(\mu_2, h_{\mu_0}(\mu_1))| \leq C_N \cdot \omega(\mu_2,
  h_{\mu_0}(\mu_1))^{-N}$ for arbitrarily large $N > 0$, provided that the
  selected atoms are regular enough.
\end{itemize}


Theorem \ref{teo:main_FIO} now follows from the observation that
\begin{equation}\label{eq:microlocal-multiply-omega}
\sum_{\mu_1} \omega(\mu_2,h_{\mu_0}(\mu_1))^{-N} \cdot 
\omega(\mu_1,\mu_0)^{-N} \leq 
C_n \cdot \omega(\mu_2,h_{\mu_0}(\mu_0))^{-(N-1)},
\end{equation}
This is an immediate consequence of properties 3 and 4 of the
pseudo-distance $\omega$, see proposition \ref{teo:omega}.

Cases involving coarse scale elements are treated similarly and we
omit the proof. The boundedness from $\ell^p$ to $\ell^p$ for every $p
> 0$ follows from the same argument as in the proof of Theorem
\ref{teo:main}.

\section{Discussion}
\label{sec:discussion}
\setcounter{theorem}{0}
\setcounter{equation}{0}

All along we specialized our discussion to the special case where the
dimension of the spatial variable is $n = 2$.  It is clear that
nothing in our arguments depends upon this specific assumption.
Indeed, we could just as well construct tight-frames of curvelets in
arbitrary dimensions by smoothly partitioning the frequency plane into
dyadic coronae, which would then be angularly localized near regions
of sidelength length $2^j$ in the radial direction and $2^{j/2}$ in
all the other directions; in order to this, we would use smooth
partitions of the unit sphere of $\R^n$ into spherical caps of radius
about $2^{-j/2}$. All of our analysis would apply as is, and would
prove versions of Theorem \ref{teo:main} in arbitrary dimensions.

Our main result assumes that the coefficients of the equation
\eqref{eq:hyperbolic-system} be smooth. In many applications of
interest, however, the coefficients may be smooth away from singular
smooth surfaces.  In geophysics for example, we typically have
different layers with very different physical properties.  A very
important question would be to know how our analysis would adapt to
this situation. In fact, it seems natural to believe that sparsity
would continue to hold in this more general setting. Intuitively, the
wave group would still be approximated by rigid motion along the
Hamiltonian flow. Only, one would need to account for possible
reflections/refractions. A curvelet hitting a singularity at a small angle of incidence would
typically produce two curvelets, a reflected and a refracted curvelet.
This is merely an intuition which one would need to justify by a
careful analysis quantifying the behavior of a curvelet near the
interface (here, the singular surface). We regard this type of
question as an important extension to this work.

The estimates proven in this paper are motivated by efforts towards
applications as sparser expansions theoretically lead to better and
faster algorithms. Our goal is to transform our theoretical insights
into effective algorithms, and derive fast accurate solvers to certain
classes of partial differential equations for which our methods have a
comparative advantage. In this direction, suppose that a discretized
version $E_N(t)$ of the curvelet matrix $E(t)$ is known in the sense
that it has been precomputed once for all.  Assuming that discretized
version inherits the sparsity from its continuous analog, then for
each initial value problem, we would have available a fast algorithm
for calculating the solution of the full wave equation. Indeed,
Corollary \ref{teo:maincoro} asserts that the truncated matrix
$\text{trunc}(E_N(t))$ obtained by keeping only $O(\epsilon^{1/m})$
terms per row obeys
\[
\|E_N(t) - \text{trunc}(E_N(t))\| \le \epsilon.
\]
Therefore, ignoring the cost of the digital curvelet transform, the
total cost of an algorithm calculating the solution of a full wave
equation to within accuracy $\epsilon$ would be linear in the number
of unknowns (number of voxels) and scale at most like $C(\epsilon)
\cdot N$. This is not the most practically relevant situation,
however, as in general, one would not know $E_N(t)$. Work in progress
shows that one can still design effective algorithms to `build' the
matrix $E_N(t)$, and for all practical purposes compute the solution
of the full wave equation for a given accuracy $\epsilon$ in $O(N
(\log N)^2)$.

Finally, we would like to conclude by pointing out that there is now
considerable evidence that the curvelet transform is a very useful
mathematical architecture; curvelets can do things that other
classical systems simply cannot do. This paper proved that
tight-frames of curvelets provide optimally sparse representations of
large classes of linear systems of hyperbolic differential equations.
But they also allow for optimally sparse representations of wavefront
phenomena \cite{CurveEdge}.  Further, they can also be useful in many
different settings. For example, they have useful microlocal features
which make them especially suitable for deployment in many inverse
problems and especially limited-angle tomography \cite{LAT}. In short,
curvelets addresses a new range of problems, going beyond what
traditional multiscale systems offer.

\section{Appendix}
\label{sec:appendix}
\setcounter{theorem}{0}
\setcounter{equation}{0}

\subsection{Additional proofs for section \ref{sec:curvelets}}
\label{sec:appendix_curvelets}

{\em Proof of proposition \ref{teo:omega}}

These four properties were already formulated in \cite{SmithWave}, although with a slightly weaker definition of pseudo-distance. Properties 1 and 2 are not proved in that reference, and property 3 is not extensively documented. We give the justification for these three results for completeness.

\begin{enumerate}
\item We are to show that $d(\mu,\mu') \asymp d(\mu',\mu)$. With $e_\mu = \xi_\mu/|\xi_\mu|$, this is
\[
|\< e_\mu, \Delta x \>| + | \Delta x |^2 + |\Delta \theta|^2 \asymp |\< e_{\mu'}, \Delta x \>| + | \Delta x |^2 + |\Delta
\theta|^2.
\]
It is sufficient to notice that
\[
|\< e_\mu, \Delta x \>| + | \Delta x |^2 + |\Delta \theta|^2 \asymp |\<
e_\mu, \Delta x \>| + |\< e_{\mu'}, \Delta x \>| + | \Delta x |^2 + |\Delta
\theta|^2.
\]
In order to justify the nontrivial inequality, use the law of cosines
illustrated in Figure \ref{fig:geometry}:
\begin{align*}
  |\< e_\mu, \Delta x \>|^2 + |\< e_{\mu'}, \Delta x \>|^2 &=
  \sin^2|\Delta \theta| \, (d_{\mu}^2 + d_{\mu'}^2) \\
  &= \sin^2|\Delta \theta| \, |\Delta x|^2 \pm 2 |\< e_\mu, \Delta x \>|
  \, |\< e_{\mu'}, \Delta x \>| \cos |\Delta \theta| \\
  &\leq \sin^2 |\Delta \theta| \, |\Delta x|^2 + 2 |\< e_\mu, \Delta x
  \>| \, |\< e_{\mu'}, \Delta x \>|.
\end{align*}
It follows that $||\< e_\mu, \Delta x \>| - |\< e_{\mu'}, \Delta x
\>|| \leq C \cdot |\Delta \theta| |\Delta x| \leq C \cdot (|\Delta
\theta|^2 + |\Delta x|^2)$ and, therefore,
\[
|\< e_\mu, \Delta x \>| + |\< e_{\mu'}, \Delta x \>| \leq C \cdot (2
|\< e_\mu, \Delta x \>| + |\Delta \theta|^2 + |\Delta x|^2).
\]
\begin{figure}
\centering
\includegraphics[width = 14cm]{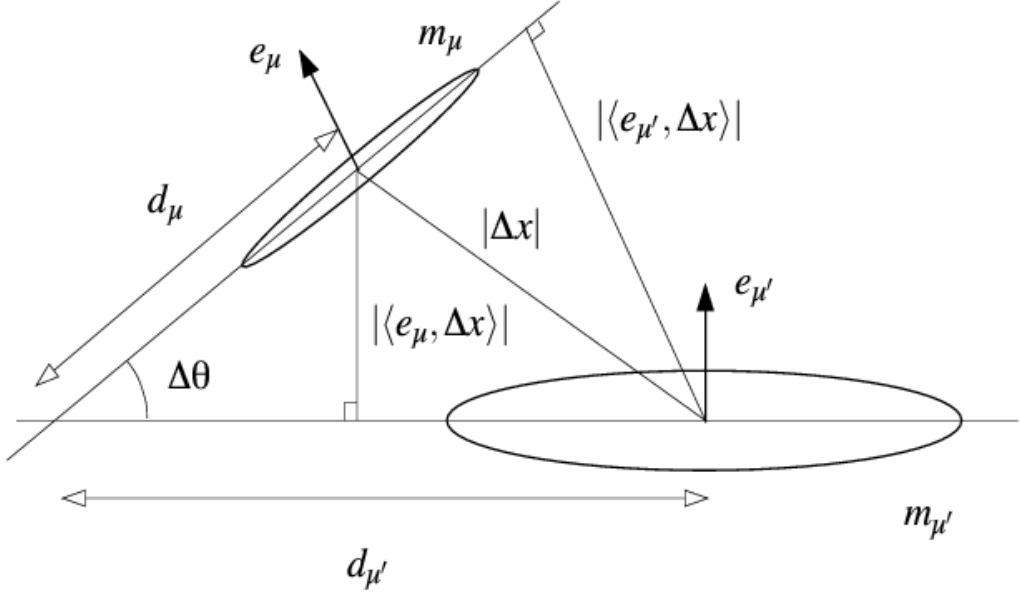}
\caption{Relative position and orientation of two curvelet molecules 
  in $x$-space. The ellipses indicate their essential support.}
\label{fig:geometry}
\end{figure}
 
\bigskip

\item Recall that $\omega(\mu,\mu') = 2^{|j-j'|}(1 + 2^{\min(j,j')}
  d(\mu,\mu'))$. Let us show that $d(\mu,\mu') \leq C \cdot
  (d(\mu,\mu'') + d(\mu'',\mu'))$. To simplify notations, set in the
  coordinates defined by $\{ e_\mu, e_\mu^\bot \}$,
\begin{align*}
  &x_{\mu} = (0,0) \qquad &x_{\mu'} = (x_1, x_2) \qquad
  &x_{\mu''} = (y_1, y_2) \\
  &e_\mu = (1,0) \qquad &e_{\mu'} = (\cos \alpha, \sin{\alpha}) \qquad 
&e_{\mu''} = (\cos \beta, \sin \beta) \\
  &|\theta_{l} - \theta_{l''}| = |\beta| \qquad &|\theta_{l'} -
  \theta_{l''}| = |\alpha - \beta|
\end{align*}
It is enough to show that there exists $\epsilon > 0$ such that
\[
\epsilon |x_1| \leq |y_1| + |\cos \alpha (x_1-y_1) + \sin \beta
(x_2-y_2)| + (|\beta| + |\alpha - \beta|)(|y_1| + |x_1 - y_1| + |y_2|
+ |x_2 - y_2|),
\]
because then $(|\beta| + |\alpha - \beta|)(|y_1| + |x_1 - y_1| + |y_2|
+ |x_2 - y_2|) \leq C \cdot (|\beta|^2 + |\alpha - \beta|^2 + |y_1|^2 +
|x_1 - y_1|^2 + |y_2|^2 + |x_2 - y_2|^2$). By contradiction let us
assume that the inequality fails. Then we must have $|y_1| < \epsilon
|x_1|$. It is always true that $|x_1 - y_1| + |y_1| \geq |x_1|$ so it
is necessary that $|\beta| + |\alpha - \beta| < \epsilon$. But then
$|\alpha| < 2 \epsilon$ thus $\cos \alpha > 1-4 \epsilon^2$ and $|\sin
\alpha| < 2 \epsilon$. The term $|\cos \alpha (x_1-y_1) + \sin \beta
(x_2-y_2)|$ is therefore always greater than $(1-4 \epsilon^2) |x_1 -
y_1| - \epsilon |x_2-y_2|$. But this quantity must also be less than
$\epsilon |x_1 - y_1|$, otherwise its sum with $|y_1|$ would exceed
$\epsilon |x_1|$. So we must have $|x_2 - y_2| > \frac{1-\epsilon - 4
  \epsilon^2}{\epsilon} |x_1 - y_1|$. But then the sum $|y_1| + |x_1 -
y_1| + |x_2 - y_2|$ must dominate $\frac{|x_1|}{2 \epsilon}$, which
implies $|\beta| + |\alpha - \beta| \leq 2 \epsilon^2$. By induction,
$\alpha = \beta = 0$ and $|y_1| + |x_1 - y_1| \geq |x_1|$ yields a
contradiction.

\bigskip

\item Let us show that $\sum_{\mu_1} \omega(\mu_0,\mu_1)^{-N} \cdot
  \omega(\mu_1, \mu_2)^{-N} \leq C_N \cdot \omega(\mu_0,
  \mu_2)^{-(N-1)}$. We closely follow and expand the argument in
  \cite{SmithWave}. We will need to use $d(\mu_0,\mu_1) \asymp
  d(\mu_1,\mu_0)$, as we have just showed.

Define $I_{\mu_1}$ by 
\begin{eqnarray}
\label{eq:sum_mu1}
  I_{\mu_1} & := & \omega(\mu_2,\mu_1)^{-N} \cdot \omega(\mu_1,\mu_0)^{-N} \\
\nonumber & = & \left(2^{|j_2 - j_1| + |j_1 - j_0|}
(1+2^{\min(j_2,j_1)} d(\mu_2,\mu_1))(1+2^{\min(j_0,j_1)}
  d(\mu_0,\mu_1))\right)^{-N}.
\end{eqnarray}
To ease notations, put temporarily $a_0 = 2^{\min(j_0,j_1)}$, $a_2 =
2^{\min(j_2,j_1)}$, $d_{01} = d(\mu_0,\mu_1)$, and $d_{12} =
d(\mu_2,\mu_1)$. We develop a lower bound on
$(1+a_2d_{12})(1+a_0d_{01}) =
1+a_2d_{12}+a_0d_{01}+a_2d_{12}a_0d_{01}$. 
We make three simple observations:
first, 
\[
a_2d_{12}+a_0d_{01} \geq \min ( a_2,a_0 ) (d_{12}+d_{01}) = A_0, \quad
\text{and} \quad d_{12}+d_{01} \geq C \cdot d(\mu_0,\mu_2);
\]
second, 
\[
a_2d_{12} + a_0 d_{01} \geq \max ( a_2d_{12},a_0d_{01} ) \geq \max (
a_2, a_0) \min ( d_{12}, d_{01} ) = B_0; 
\]
and third
\begin{eqnarray*}
a_2d_{12}a_0d_{01} &  = & 
\max ( a_2, a_0 ) \min ( a_2, a_0 ) \max ( d_{12},
  d_{01} ) \min ( d_{12}, d_{01} ) \\
& \geq & \max ( a_2, a_0 ) \min ( a_2, a_0 )
  \min ( d_{12}, d_{01} ) \frac{d_{12}+d_{01}}{2} = A_0 B_0/2.
\end{eqnarray*}
This gives
\begin{equation*}
1 + a_2 d_{12} + a_0d_{01} + a_2 d_{12} 
a_0 d_{01} \ge  \frac{1}{2}\left(1 + A_0 + B_0 
+ A_0 B_0\right) \ge \frac{1}{2} (1 + A_0) (1 + B_0).  
\end{equation*}
We replace the values of $A_0$, $B_0$ by their expression, use the
relation $A_0 \ge d(\mu_0,\mu_2)$ and obtain
\begin{equation}\label{eq:I1}
I_{\mu_1} \le C \cdot 2^{-(|j_2 - j_1|+|j_0-j_1|)N} \cdot
\left(1+2^{\min(j_2,j_0,j_1)} d(\mu_2,\mu_0)\right)^{-N} \cdot  (L_1)^{-N}  
\end{equation}
with 
\[
L_1 = 1+ \max(2^{\min(j_2,j_1)},2^{\min(j_0,j_1)}) \min(d_{01},
d_{12}).
\]
Note that
\begin{eqnarray*}
L_1 & = & \min \left( 1+ \max(2^{\min(j_2,j_1)},2^{\min(j_0,j_1)})
d_{12}, 1+ \max(2^{\min(j_2,j_1)},2^{\min(j_0,j_1)})
d_{01} \right) \\
& \ge & \min \left( 1+ 2^{\min(j_2,j_1)}
d_{12}, 1+ 2^{\min(j_0,j_1)} d_{01} \right) 
\end{eqnarray*}
and, therefore, 
\begin{eqnarray*}
(L_1)^{-N} & \le & 
\max\left((1+ 2^{\min(j_2,j_1)} d_{12})^{-N}, (1+ 2^{\min(j_0,j_1)} 
d_{01})^{-N}\right)\\
& \le & (1+ 2^{\min(j_2,j_1)} d_{12})^{-N} + 
(1+ 2^{\min(j_0,j_1)} 
d_{01})^{-N}.
\end{eqnarray*}

In the sequel we will repeatedly make use of the bound
\begin{equation}\label{eq:sumkl}
\sum_{k,\ell} (1+2^{q} d(\mu,\mu'))^{-N} \leq C \cdot 2^{2(j-q)_+},
\end{equation}
valid for $N \geq 2$, any real $q$ and where the subscript $+$ denotes
the positive part. This is justified as follows. Without loss of
generality, assume that $\mu' = (j',0,0)$ so that the curvelet
$\gamma_{\mu'}$ is nearly vertical and centered near the origin. We
recall that $\Delta \theta = \pi \cdot \ell \cdot 2^{-\lfloor j/2
  \rfloor}$, $\ell = 0, 1, \ldots 2^{\lfloor j/2 \rfloor} - 1$, and
$x_{\mu} = R_{\theta_{\mu}} D_{j}^{-1} k$, say. Then the left-hand side
is
\begin{equation}\label{eq:lhssumkl}
\sum_{\ell = 0}^{2^{\lfloor j/2 \rfloor} - 1} \sum_{k \in \Z^2} \left(1
    + 2^{q}(|2^{-j/2}\ell|^2 + |2^{-j/2}k_2|^2 + |2^{-j} k_1|)\right)^{-N}.
\end{equation}
For $j \geq q$ this can be seen as a Riemann sum and bounded---up to a
numerical multiplicative constant---by the corresponding integral
\begin{equation*}
  \int_{\R^2} \frac{dx}{2^{-3j/2}} \int_{\R} \frac{dy}{2^{-j/2}} [1 + 2^{q}(y^2 + x_2^2 + |x_1|)]^{-N}
\end{equation*}
which in turn is less than $C \cdot 2^{2(j-q)}$ provided $N \geq 2$.
For $j \leq q$, the sum (\ref{eq:lhssumkl}) essentially consists of a
few terms, giving a $O(1)$ contribution. This gives the bound $C \cdot
2^{2(j-q)_+}$.

By symmetry, we can now assume $j_0 \leq j_2$. Let us consider three
cases.
\begin{itemize}
\item $\mathbf{0 \leq j_2 \leq j_1}$. In that case we have the bound
\[
(L_1)^{-N} \leq C \cdot [(1+2^{j_2} d_{01})^{-N} + (1+2^{j_2}d_{12})^{-N}]. 
\]
Summing this quantity over $k_1$ and $\ell_1$ i.e., over all $\mu_1$
that correspond to a given $j_1$, and using (\ref{eq:sumkl}), we
obtain for $j_1 \geq j_2$
\begin{align*}
  \sum_{\mu_1} I_{\mu_1} &\leq C \cdot (1 + 2^{j_0} d_{02})^{-N}  \sum_{j_1 \geq j_2} 2^{-(2j_1 - j_0 - j_2)N} \cdot 2^{2(j_1-j_2)} \\
  &\leq C \cdot 2^{-(j_2-j_0)N} (1 + 2^{j_0} d_{02})^{-N} = C \cdot
  \omega(\mu_0,\mu_2)^{-N}.
\end{align*}
\item $\mathbf{0 \leq j_1 \leq j_0}$. We now have
\[
(L_1)^{-N} \leq C \cdot [(1+2^{j_1} d_{01})^{-N} + (1+2^{j_1}d_{12})^{-N}]. 
\]
According to (\ref{eq:sumkl}), the sum over $k_1$ and $\ell_1$ of
$(L_1)^{-N}$ is bounded by a constant independent of $j_1$. The
remaining sum is
\[
\sum_{\mu_1} I_{\mu_1} \leq C \cdot 2^{-(j_0 + j_2)N} \sum_{j_1 \leq
  j_0} 2^{2j_1 N} \cdot (1 + 2^{j_1} d_{02})^{-N}.
\]
Observe that $2^{j_1 N} (1 + 2^{j_1} d_{02})^{-N} \leq 2^{j_0 N} (1 +
2^{j_0} d_{02})^{-N}$, therefore
\[
\sum_{\mu_1} I_{\mu_1} \leq C \cdot 2^{-(j_2-j_0)N} (1 + 2^{j_0}
d_{02})^{-N} = C \cdot \omega(\mu_0,\mu_2)^{-N}.
\]
\item $\mathbf{j_0 \leq j_1 \leq j_2}$. In that case we still have
\[
(L_1)^{-N} \leq C \cdot [(1+2^{j_1} d_{01})^{-N} + (1+2^{j_1}d_{12})^{-N}]. 
\]
summed over $k_1$ and $\ell_1$ into a $O(1)$ contribution. What
remains is
\begin{align*}
  \sum_{\mu_1} I_{\mu_1} \leq C \cdot 2^{-(j_2 - j_0)N} (1 + 2^{j_0} d_{02})^{-N} \sum_{j_0 \leq j_1 \leq j_2} 1 \\
  \leq C \cdot \omega(\mu_0,\mu_2)^{-(N-1)}.
\end{align*}
\end{itemize}
We conclude by collecting the estimates corresponding to the three
different cases. Remark that the loss of one (fractional) power of
$\omega$ in the third case is unavoidable unless one modifies its
definition in the spirit of \cite{SmithWave}. This would however make
notations unnecessarily heavy.

\bigskip

\item See \cite{SmithWave} p. 804.

\end{enumerate}

\bigskip

{\em Proof of the inequality (\ref{eq:freq-loc})}.  Assume without
loss of generality that $\mu = \mu_0$.  We may express $S_{\mu'}(\xi)$ as
$S_{\mu'_0}(R_{\Delta \theta} \xi)$, with $\Delta \theta = \theta_\mu -
\theta_{\mu'}$. We begin by expressing the integral in polar
coordinates,
\begin{align*}
  \xi_1 &= r \cos \theta \qquad \qquad \qquad & (R_{\Delta
    \theta}\xi)_1 = r
  \cos (\theta + \Delta \theta), \\
  \xi_2 &= r \sin \theta \qquad \qquad \qquad & (R_{\Delta
    \theta}\xi)_2 = r \sin (\theta + \Delta \theta).
\end{align*}
As we can see, the cosine factor is not crucial and we may just as
well drop it. Consequently,
\begin{align*}
\int | S_{\mu}(\xi) \, S_{\mu'}(\xi)|^n \, d\xi \leq C \cdot
&\int_0^{\infty} r dr \, \frac{1}{[1+2^{-j}r]^N}
\frac{1}{[1+2^{-j'}r]^N} \\
\times &\int_0^{2 \pi} d\theta [1+a|\sin
\theta|]^{-N} [1+a'|\sin (\theta + \Delta \theta)|]^{-N}, 
\end{align*}
where $a = \frac{2^{-j/2} r}{1+2^{-j} r}$ and $a' = \frac{2^{-j'/2}
  r}{1+2^{-j'} r}$. This decoupling makes the problem of bounding the
inner integral on the variable $\theta$ tractable. For example when $a
> a' > 1$, following \cite{Mey2} p.56,
\[
\int_{-\infty}^{\infty} d\theta [1+a|\theta|]^{-N} [1+a'|\theta +
\Delta \theta|]^{-N} \leq C \cdot \frac{1}{a} \frac{1}{[1 + a' |\Delta
  \theta|]^N}.
\]
We get other estimates for other values and orderings of $a$ and $a'$.
The integral on $r$ is then broken up into several pieces according to
the values of $a$, $a'$, $j$ and $j'$. It is straightforward to show
that each of these contributions satisfies the inequality
(\ref{eq:freq-loc}).

\subsection{Additional proofs for section \ref{sec:FIO}}
\label{sec:appendix_FIO}

{\em Proof of lemma \ref{teo:atoms}}

By definition $a^{(\mu)}(x) = 2^{-3j/4} m_\mu\left(D_{2^{-j}} R_{\theta_\mu} x - k
\right)$ and, therefore,
  \begin{eqnarray}
    \label{eq:atom-expression}
    \nonumber 
a^{(\mu)}(x) & = & \frac{1}{\alpha_\mu} \int (Rf)(a,\theta, b) a^{3/4} 
 2^{-3j/4} \psi(D_a R_\theta (R^{-1}_{\theta_\mu} D_{2^j}(x + k) - b) \, 
  d\mu \\
 & = &  \frac{1}{\alpha_\mu} \int (Rf)(a,\theta, b) |A|^{1/2} 
 \psi\left(A(x - (\beta - k)\right)  \, d\mu, 
  \end{eqnarray}
where $A = D_a R_\delta D_{2^j}$ with $\delta = \theta - \theta_\mu$ and 
$\beta = D_{2^-j} R_{\theta_\mu} b$.  

Let us first verify the assertion about the support of $a^{(\mu)}$.
Recall that over a cell $Q_\mu$, $\beta \in [k_1, k_1 + 1) \times
[k_2, k_2 + 1)$, and hence for all $b \in Q_\mu$, we have 
\[
\text{Supp}\, \psi\left(A(x - (\beta - k))\right) \subset \text{Supp} \, 
\psi(Ax) + [0,1]^2.
\]
Next $\text{Supp} \, \psi(Ax) \subset A^{-1}[0,1]^2$ with $A^{-1} =
D_{2^{-j}} R_{-\delta} D_a^{-1}$. It is not difficult to check that
$A^{-1}[0,1]^2 \subset [c_1, c_2) \times [d_1, d_2)$ which then gives
(\ref{eq:compact-support}).

There are several ways to prove the property about nearly vanishing
moments. A possibility is to show that the Fourier transform of
$a^{(\mu)}$ is appropriately small in a neighborhood of the axis $\xi_1
= 0$. We choose a more direct strategy and show that 
\begin{equation}
  \label{eq:nearly-vanishing-moments2}
\left| \int \psi(A(x - \beta)) x_1^k \, dx_1 \right| \le C_m 
\cdot 2^{-j(m+1)}.  
\end{equation} 
uniformly over the $(a, \theta, b) \in Q_\mu$. The property
(\ref{eq:nearly-vanishing-moments}) follows from this fact. Indeed,
\[
\int a^{(\mu)}(x_1, x_2) \, x_1^k \, dx_1  = \frac{1}{\alpha_\mu} 
\int_{Q_\mu} R f(a,\theta, b) d\mu \int |A|^{1/2}  
\psi(A(x - \beta)) x_1^k \, dx_1, 
\]
and the Cauchy-Schwarz inequality gives 
\begin{eqnarray*}
\left| \int a^{(\mu)}(x_1, x_2) \, x_1^k \, dx_1 \right| & \le & 
\frac{1}{\alpha_\mu} \cdot \|Rf\|_{L_2(Q_\mu)} \cdot  \left(
\int_{Q_\mu} \left|  \int |A|^{1/2}  
\psi(A(x - \beta)) x_1^k \, dx_1 \right|^2 \, 
d\mu\right)^{1/2} \\
& =  & \left(\int_{Q_\mu} \left|  \int |A|^{1/2}  
\psi(A(x - \beta)) x_1^k \, dx_1 \right|^2 \, 
d\mu\right)^{1/2}.
\end{eqnarray*}
The uniform bound (\ref{eq:nearly-vanishing-moments2}) together with
the fact that $\int_{Q_\mu} d\mu$ is either $3\pi$ or $3\pi/2$ gives
(\ref{eq:nearly-vanishing-moments}).

We then need to establish (\ref{eq:nearly-vanishing-moments2}).  Let
$\d_2$ be $\partial/\partial x_2$, recall that by assumption
(\ref{eq:psi-choice})--(\ref{eq:psiD-vanish}), we have that for all
$x_2 \in \R$,
\[
\int \d_2^n \psi(x_1, x_2) x_1^k \, dx_1 = 0, 
\quad k = 0, 1,\ldots, R,  
\]
and more generally, for each $\alpha \neq 0$ and $\beta$
\begin{equation}
\label{eq:vanish}
\int \d_2^n \psi(\alpha x_1 + \beta, x_2) x_1^k \, dx_1 = 0, 
\quad k = 0, 1,\ldots, R. 
\end{equation}

We shall use (\ref{eq:vanish}) to prove
(\ref{eq:nearly-vanishing-moments2}).  Letting
\[
A = \begin{pmatrix} a_{11} & a_{12}\\
a_{21} & a_{22} \end{pmatrix}
\]
and with the same notations as before, a simple calculation shows that
$a_{21} = - \frac{2^{-j} \sin \delta}{\sqrt{a}}$. As $a \ge 2^{-(j+1)}$ and
$|\delta| \le \pi/2 \cdot 2^{-\lfloor j/2 \rfloor}$, we have 
\begin{equation}
\label{eq:epsilon}
|a_{21}| \le c \cdot 2^{-j}. 
\end{equation}
We then write 
\begin{eqnarray*}
\psi(A x) & = & \psi(a_{11} x_1 + a_{12} x_2, a_{21} x_1 + a_{22} x_2)\\
          & = & \sum_{n = 0}^{N-1} D^{n} \psi(a_{11} x_1 + a_{12} x_2,
          a_{22} x_2) \frac{(a_{21} x_1)^n}{n!} + O((a_{21} x_1)^N)
\end{eqnarray*}
and, therefore, 
\[
\int \psi(Ax) x_1^k  \, dx_1 = \sum_{n = 0}^{N-1}  
\frac{a_{21}^n}{n!}  \int D^{n} \psi(a_{11} x_1 + a_{12} x_2,
      a_{22} x_2) x_1^{n + k} \, dx_1  + O(a_{21}^N) 
\]
Fix $k \le D$ and pick $N = D - k + 1$ so that for $n = 0, 1,
\ldots, N -1$,  $n + k \le D$. By virtue of (\ref{eq:vanish}) 
all the integrals in the sum vanish and the only remaining term is 
$O(a_{21}^N)$ which because of (\ref{eq:epsilon}) is $O(2^{-jN})$. 
As a consequence, setting $m = D/2$, we conclude that 
\[
\left| \int \psi(Ax) x_1^k \, dx_1 \right| \le C_m \cdot 2^{-j(m+1)},
\quad k = 0,1, \ldots m;
\]
this is the content of (\ref{eq:nearly-vanishing-moments2}).

The careful reader will notice that inequality
(\ref{eq:nearly-vanishing-moments2}) or equivalently
(\ref{eq:nearly-vanishing-moments}) is a weaker statement than
inequality (\ref{eq:bandpass}) for the definition of nearly vanishing
moments. There is no doubt that the stronger estimate
(\ref{eq:bandpass}) also holds for curvelet atoms. The proof of this
fact uses standard arguments and we choose not to reproduce it here.

Last, the regularity property is a simple consequence of the Cauchy
Schwarz inequality;
\begin{eqnarray*}
\left| a^{(\mu)}(x_1, x_2) \right| & \le & 
\frac{1}{\alpha_\mu} \int |Rf(a, \theta, b) | |A|^{1/2}
\|\psi\|_{L_\infty} \, d\mu\\
& \le  &  
\|\psi\|_{L_\infty} \cdot \frac{1}{\alpha_\mu} 
\|Rf\|_{L_2(Q_\mu)} \cdot \left( \int_{Q_\mu} \, |A| d\mu\right)^{1/2} \\
& = & 2\sqrt{3\pi} \cdot \|\psi\|_{L_\infty}.  
\end{eqnarray*}
These last inequalities used the facts that $|A| \le 4$ for $(a,
\theta, b) \in Q_\mu$ and $ \int_{Q_\mu} \, d\mu \le 3\pi$.  Estimates
for higher derivatives are obtained in exactly the same fashion--after
differentiation of the integrand.  This finishes the proof of the
lemma.

\bigskip

{\em Proof of (\ref{eq:molec-atom})}.  Recall that
\[
\alpha_{\mu\mu'} = \left( \int_{Q_{\mu'}} |
  \mathcal{R}(q(D)\gamma_{\mu})(a,b,\theta) |^2 \, d\mu \right)^{1/2}.
\]
The first thing to notice is that $q(D) \gamma_{\mu}$ is still a
family of curvelet molecules, because $q(\xi)$ is a multiplier of
order zero. Since $\psi_{a,\theta,b}$ also obeys the molecule
properties, lemma \ref{teo:almost_orthogonality} implies the corresponding
almost-orthogonality condition. Integrating over $Q_{\mu'}$ does not
compromise this estimate, as can be seen by applying the
Cauchy-Schwarz inequality.

\bigskip

{\em Proof of inequality (\ref{eq:bigOone})}.  Derivatives of
$\hat{\gamma}_{\mu}$ and $\sigma$ are treated using the following
estimates.
\begin{align*}
  &|\partial_{\xi}^{\alpha} \hat{\gamma}_{\mu}(\xi)| \leq C_{\alpha}
  \cdot 2^{-3j/4}
  2^{-\alpha_1 j} 2^{-\alpha_2 j/2} \\
  &|\partial_{\xi}^{\alpha} \sigma(\phi^{-1}(x),\xi)| \leq C_{\alpha} \cdot
  2^{-|\alpha| j} \qquad \text{on } W_{\mu} =
  \mbox{supp}(\hat{\gamma}_{\mu}).
\end{align*}
We now develop size estimates for the phase perturbation $\delta$.
Following closely the discussion in \cite{Ste}, p.407, we claim 
that on $W_{\mu}$,
\begin{equation}\label{eq:esth1}
| \partial_{\xi}^{\alpha} \partial_{x}^{\beta} \delta(x,\xi)| \leq 
C_{\alpha \beta} \cdot 2^{-\alpha_1 j} 2^{-\alpha_2 j/2}.
\end{equation}
The derivations in $x$ add no complications. Hence, assume that $\beta
= 0$.  As the above result (\ref{eq:esth1}) relies upon the homogeneity
of the phase with respect to $\xi$, we recall a few useful facts about
homogeneous functions of degree one:
\begin{itemize}
\item[] $\Phi = \Phi_{\xi} \cdot \xi \;$ (Euler's theorem), 
\item[] $\Phi_{\xi \xi} \cdot \xi = 0 \;$ (differentiate the above
  relation), 
\item[] $\partial_{\xi}^{\alpha} \Phi = O(|\xi|^{1-|\alpha|}) \;$. 
\end{itemize}
It follows from the definition that $\delta(x,\xi_1,0) = 0$ and
likewise $\frac{\partial \delta}{\partial \xi_2}(x,\xi_1,0) = 0$. Thus
for every $n$, $\frac{\partial^n \delta}{\partial \xi_1^n} (x,\xi_1,0)
= 0$ and $ \frac{\partial}{\partial \xi_2} \frac{\partial^n
  \delta}{\partial \xi_1^n}(x,\xi_1,0) = 0$. Recall that the support
conditions are $|\xi_1| \leq C \cdot 2^j$ and $|\xi_2| \leq C \cdot 2^{j/2}$.
Taylor series expansions about $\xi_2 = 0$ together with homogeneity
assumptions give 
\begin{align*}
  \frac{\partial^{\alpha_1}}{\partial \xi_1^{\alpha_1}}\delta(x,\xi) &=
  O(|\xi_2|^2 |\xi|^{-1-\alpha_1}) = O(2^{-\alpha_1 j}), \\
  \frac{\partial}{\partial \xi_2} \frac{\partial^{\alpha_1}} {\partial
    \xi_1^{\alpha_1}}\delta(x,\xi) &= O(|\xi_2| |\xi|^{-1-\alpha_1})
  = O(2^{-j/2} 2^{-\alpha_1 j}), \\
  \frac{\partial^{\alpha_2}}{\partial \xi_2^{\alpha_2}}
  \frac{\partial^{\alpha_1}}{\partial \xi_1^{\alpha_1}}\delta(x,\xi) &=
  O(|\xi|^{1-\alpha_1-\alpha_2}) = O(2^{-\alpha_1 j} 2^{-\alpha_2
    j/2}) \qquad \text{when } \alpha_2 \geq 2,
\end{align*}
as claimed. The point about these estimates is that they exhibit
exactly the parabolic scaling of curvelets. We conclude
\[
|\partial_{\xi}^{\alpha} e^{i \delta(\phi^{-1}(x),\xi)}| \leq C_{\alpha} \cdot
2^{-\alpha_1 j} 2^{-\alpha_2 j/2} \qquad \text{on } W_{\mu}
\]
and therefore (\ref{eq:bigOone}).


\begin{thebibliography}{99}
  
\bibitem{AM} J.P. Antoine, R. Murenzi, Two-dimensional directional
  wavelets and the scale-angle representation, \emph{Sig. Process.},
  \textbf{52} (1996) 259-281
  
\bibitem{BCR} G. Beylkin, R. Coifman, V. Rokhlin, Fast wavelet
  transforms and numerical algorithms, \emph{Comm. on Pure and Appl.
    Math.}, \textbf{44} (1991), 141-183
  
\bibitem{CurveFIO} E. J. Cand\`{e}s, L. Demanet, Curvelets and Fourier
  Integral Operators, \emph{C. R. Acad. Sci. Paris, Ser. I},
  \textbf{336} (2003), 395-398

\bibitem{Curvelets-StMalo} E. J. Cand\`{e}s, D. L. Donoho, 
  Curvelets - a suprisingly
  effective nonadaptive representation for objects with edges, in
  \emph{Curves and Surface Fitting : Saint-Malo 1999}, A. Cohen, C.
  Rabut, L. Schumaker (eds.), Vanderbilt University Press, Nashville,
  2000, 105-120
  
\bibitem{CurveEdge} E. J. Cand\`{e}s, D. L. Donoho, New Tight Frames
  of Curvelets and Optimal Representations of Objects with Piecewise
  $C^2$ Singularities, \emph{Comm. on Pure and Appl.  Math.},
  \textbf{57} (2004), 219-266.
  
\bibitem{DWTPSI} E. J. Cand\`{e}s, D. L. Donoho, Directional Wavelet
  Transform via Parabolic Scaling: I. Resolution of the Wavefront Set,
  Technical Report, Stanford University, 2003.  
  
\bibitem{DCTvUSFFT} E. J. Cand\`{e}s, D. L. Donoho, DCTvUSFFT: Digital
  Curvelet Transforms via Unequispaced Fast Fourier Transforms,
  Technical Report, California Institute of Technology, 2004.
  
\bibitem{LAT} E. J. Cand\`{e}s, D. L. Donoho, Curvelets: New Tools for
  Limited-Angle Tomography, Technical Report, California Institute of
  Technology, 2004.

\bibitem{TVSynthesis} E. J. Cand\`{e}s, F. Guo, New multiscale
  transforms, minimum total variation synthesis : Application to
  edge-preserving image reconstruction, \emph{Sig. Process.}, special
  issue on Image and Video Coding Beyond Standards, \textbf{82}
  (2002), 1519--1543.
  
\bibitem{CF} A. C\'{o}rdoba, C. Fefferman, Wave packets and Fourier
  integral operators, \emph{Comm. PDE's}, \textbf{3(11)} (1978),
  979-1005
  


\bibitem{Donoho-IEEE} D. L. Donoho, M. Vetterli, R. A. DeVore, I. Daubechies, Data compression and harmonic analysis, \emph{IEEE Trans. Inform. Theory}, \textbf{44}, (1998), 2435-2476

\bibitem{DonohoSoft} D. L. Donoho, De-noising by Soft-Thresholding, \emph{IEEE Trans. Inform. Theory}, \textbf{41}, (1995), 613-627

\bibitem{Dui} J. Duistermaat, \emph{Fourier Integral Operators},
  Birkhauser, Boston, 1996
  
  
\bibitem{Gar} L. Garding, \emph{Singularities in Linear Wave
    Propagation}, Lecture notes in Math 1241, Springer, 1987
 
\bibitem{Fef} C. Fefferman, A note on spherical summation multipliers,
  \emph{Israel J. Math.} \textbf{15} (1973) 44-52
 
\bibitem{Fol} G. Folland, \emph{Harmonic Analysis in Phase Space},
  Princeton University Press, 1989
 
\bibitem{FJW} M. Frazier, B. Jawerth, G. Weiss, \emph{Littlewood-Paley
    Theory and the Study of Function Spaces}, CBMS 79, AMS,
  Providence, 1991
  
\bibitem{Lax} P. Lax, Asymptotic solutions of oscillatory initial
  value problems, \emph{Duke Math J.}, \textbf{24} (1957) 627-646

\bibitem{Lemarie} P. G. Lemari\'e and Y. Meyer,  Ondelettes
  et bases {H}ilbertiennes.  \emph{Rev. Mat. Iberoamericana},
  {\bf 2} (1986), 1-18.
  
\bibitem{MeyerAlgo} Y. Meyer.  \emph{Wavelets: Algorithms and
    Applications}, SIAM, Philadelphia, 1993.

\bibitem{MallatBook} S. Mallat, \emph{A Wavelet Tour of Signal Processing}, 2nd ed., Academic Press, Orlando-San Diego, 1999.

\bibitem{Mey1} Y. Meyer, \emph{Ondelettes et Op\'{e}rateurs}, Hermann,
  Paris, 1990
  
\bibitem{Mey2} Y. Meyer, R. Coifman, \emph{Wavelets,
    Calder\'{o}n-Zygmund and Multilinear Operators}, Cambridge
  Univresity Press, Cambridge, 1997
  
\bibitem{SSS} A. Seeger, C. Sogge, E. Stein, Regularity properties of
  Fourier integral operators, \emph{Annals of Math.} \textbf{134}
  (1991), 231-251
  
\bibitem{Seeley} R. Seeley, Complex powers of an elliptic operator,
  \emph{Proc. Symp. Pure Math.} \textbf{10} (1968), 288-307

\bibitem{Sm1} H. Smith, A Hardy space for Fourier integral operators,
  \emph{J. Geom. Anal.} \textbf{7} (1997)
  
\bibitem{SmithWave} H. Smith, A parametrix construction for wave equations
  with $C^{1,1}$ coefficients, \emph{Ann. Inst. Fourier}, Grenoble,
  \textbf{48}, 3 (1998), 797-835
  
\bibitem{Sog} C. Sogge, \emph{Fourier Integrals in Classical
    Analysis}, Cambridge University Press, 1993
  
\bibitem{Ste} E. Stein, \emph{Harmonic Analysis}, Princeton University
  Press, Princeton NJ, 1993
  
\bibitem{Stolk_deHoop} C. Stolk, M. de Hoop, Microlocal Analysis of
  Seismic Inverse Scattering in Anisotropic Elastic Media, \emph{Comm.
    Pure and Appl. Math.}, \textbf{55} (2002), 261-301

  
\bibitem{Tay} M. Taylor, Reflection of Singularities of Solutions to
  Systems of Differential Equations, \emph{Comm. Pure and Appl.
    Math.}, \textbf{28} (1975), 457-478

 
  
\bibitem{Tre} F. Tr\`{e}ves, \emph{Introduction to Pseudo-Differential
    and Fourier Integral Operators}, Plenum press, 1982, 2 volumes.
  
\bibitem{Whi} G. Whitham, \emph{Linear and Nonlinear waves}, Wiley
  Interscience, 1999.

\end{thebibliography}

\end{document}